\documentclass[12pt]{article}

\usepackage{amsfonts,amssymb,amsmath,eucal,pstricks,pst-plot,pst-node}
\usepackage[dvips]{graphics}

\numberwithin{equation}{section}

\usepackage{graphicx}

\newcommand{\R}{{\mathbb R}}
\newcommand{\C}{{\mathbb C}}
\newcommand{\Q}{{\mathbb Q}}
\newcommand{\N}{{\mathbb N}}
\newcommand{\Z}{{\mathbb Z}}
\newcommand{\SL}{{\mathcal S} \ell}

\newcommand{\half}{\frac{1}{2}}

\newcommand{\s}{\bar{s} }
\newcommand{\y}{\bar{y}}

\newcommand{\bareta}{\bar{\eta}}
\newcommand{\bDelta}{\bar\Delta}
\newcommand{\bz}{\bar{z}}

\newtheorem{theo}{{\sc Theorem}}[section]
\newtheorem{cor}[theo]{{\sc Corollary}}

\newtheorem{lem}[theo]{{\sc Lemma}}

\newtheorem{prop}[theo]{{\sc Proposition}}
\newtheorem{defn}[theo]{{\sc Definition}}

\newtheorem{prob}{{\sc Problem}}[section]

\begin{document}

\author{Steve Zelditch}
{\thanks{Research partially supported by  NSF grants \#
DMS-0071358, DMS-0302518.}}

\title{The inverse spectral problem  }

\maketitle

\section{Introduction}

  The  inverse spectral problem on a Riemannian manifold
$(M, g)$, possibly with boundary,  is to determine as much as
possible of the geometry of  $(M, g)$ from the spectrum of its
Laplacian $\Delta_g$ (with some given boundary conditions).
 The special
inverse problem of Kac is to determine a Euclidean domain $\Omega
\subset \R^n$ up to isometry  from the spectrum $Spec_B(\Omega)$
of its Laplacian $\Delta_B$ with Dirichlet, Neumann or more
general boundary conditions $B$. The physical motivation is to
identify physical objects from the light or sound they emit, which
may be all that is observable of remote objects such as stars or
atoms.

The inverse spectral problem is just one among many kinds of
inverse problems whose goal is to determine a metric, domain or
scatterer from physically relevant invariants. A comparison with
other inverse problems shows just how small a set of invariants
the spectrum is. For instance, the  boundary  inverse problem asks
to determine the metric $g$ on a fixed  bounded domain $\Omega
\subset M$ of a Riemannian manifold $(M, g)$ from the spectrum of
the Dirichlet Laplacian on $L^2(\Omega)$, and from the Cauchy data
$\frac{\partial \phi_{\lambda}}{\partial \nu}|_{\partial \Omega}$
of its eigenfunctions, or equivalently from  its
Dirichlet-to-Neumann operator \cite{BK, KKL, LU}. The inverse
scattering problem seeks to determine an obstacle from its
scattering amplitude \cite{Ma, Ma2}, or from its scattering length
spectrum \cite{St, St2}, a set of lengths parametrized by $S^{n-1}
\times S^{n-1}$. By comparison, the inverse spectral problem has
to make do  with just a discrete `unformatted' set of eigenvalues
(or resonance poles in the open case).

It is probably a consequence of this relative poverty  of
invariants that most of the results in inverse spectral theory
over the last two decades (especially since the appearance of
Sunada's aricle \cite{Su}) are `negative results', showing that
one cannot determine metrics or domains by their spectra. The
collection of non-isometric isospectral pairs of Riemannian
manifolds would require a lengthy survey of its own and for that
we refer to the lectures of C. Gordon   (cf. \cite{Gor, Gor2}).
By comparison, the number of `positive results' showing that one
can indeed recover a domain or metric is rather small.

This  survey is devoted to the positive results. The emphasis is
on relations between the spectrum of the Laplacian and the
dynamics of the geodesic flow $G^t: S^*M \to S^*M$. Most of the
new material concerns wave trace invariants and their
 applications to solving concrete inverse spectral problems.  The
 wave group is the quantization of the geodesic flow, and so wave
 trace methods often `reduce' inverse spectral problems to
 inverse dynamical problems. Because of their relevance, we have attempted to describe inverse
 dynamical problems and results.

\subsection{Some basic inverse spectral problems}

 Let us introduce some basic terminology.  The spectrum of a
 compact Riemannian manifold defines a map
$$Spec: {\cal M} \to \R_+^{\N}, \;\;\; (g, B) \to  Spec(\Delta_{g, B}) =  \{\lambda_0  < \lambda_1^2 \leq \lambda_2^2
\leq \cdots \}$$
 from some class of metrics ${\mathcal M}$ on a manifold $M$  to
 the spectrum of its Laplacian,
 $$\left\{\begin{array}{l}  \Delta \phi_j =
\lambda_j^2 \phi_j,\;\; \langle \phi_i, \phi_j \rangle = \delta_{i
j}  \\ \\
B\phi_j = 0\;\; \mbox{on}\; \partial M, \end{array} \right.
$$
with boundary conditions $B: C^{\infty}(M) \to C^{\infty}(\partial
M)$ if  $\partial M \not= \emptyset$. Here, $\Delta$ denotes the
{\it positive} Laplacian $$\Delta = -
\frac{1}{\sqrt{g}}\sum_{i,j=1}^n \frac{\partial}{\partial
x_i}g^{ij} \sqrt{g} \frac{\partial}{\partial x_j}$$ of  a
Riemannian manifold $(M,g)$, where    $g_{ij} =
g(\frac{\partial}{\partial x_i},\frac{\partial}{\partial x_j}) $,
$[g^{ij}]$ is the inverse matrix to $[g_{ij}]$ and $g = {\rm det}
[g_{ij}].$
 We will only consider
Dirichlet $B u = u|_{\partial M}$ and Neumann $B u =
\partial_{\nu} u |_{\partial M}$.
 Eigenvalues are
repeated according to their multiplicities.

Two metrics or domains are called isospectral if they have the
same spectrum.  The main problem in inverse spectral theory is to
describe the possible spectra $\Lambda \subset \R^{{\bf N}}$ of
Laplacians and, for each possible spectrum, to describe the
metrics or domains in the  spectral class
\begin{equation} Spec^{-1} (\Lambda). \end{equation}
Somewhat simpler is to describe the possible smooth curves in the
isospectral class, since it apriori eliminates  irregular subsets.
An isospectral deformation of a Riemannian manifold (possibly with
boundary)  is one-parameter family of metrics satisfying $Spec(M,
g_t) = Spec(M, g_0)$ for each $t$. Similarly, an isospectral
deformation of a domain with a fixed background metric $g_0$ and
boundary conditions $B$ is a family $\Omega_t$ with
$Spec_B(\Omega_t) = Spec_B(\Omega)$. One could also pose the
inverse spectral problems for  boundary conditions (while holding
the other data fixed) as in \cite{GM2, PT}.

The  inverse spectral and isospectral deformation problems are
difficult because the map $Spec$ is highly nonlinear. The
linearization of the problem is to find infinitesimal isospectral
deformations, i.e. deformations for which the eigenvalue
variations vanish to first order.
 By first order perturbation theory, the variations of the
eigenvalues under a variation of the metric are given by
\begin{equation}\label{DOTLAMBDA} \dot{\lambda}_j\; = \;\frac{d}{dt} \lambda_j(t)|_{t = 0}\;  = \;\langle \dot{\Delta}
\phi_j, \phi_j \rangle, \end{equation} where $\dot{\Delta}$ is the
variation of the Laplacian and where  $\phi_j = \phi_j(0)$ is an
orthonormal basis of eigenfunctions which varies smoothly in $t$
(such a basis exists by the Kato-Rellich theory). As will be
recalled below, lengths $L_{\gamma}$ of closed geodesics $\gamma$
are spectral invariants (at least, when there is at most one
closed geodesic of each length), so infinitesimal
iso-length-spectral deformations are those where
\begin{equation} \label{LE} \dot{L}_{\gamma} = 0, \;\;\; \forall \gamma.
\end{equation}
The deformation of the metric is a symmetric $2$-tensor $\dot{g}$,
so the linearized problem is to determine the space of $\dot{g}
\in S^2 T^*M$ (modulo tensors arising from diffeomorphisms
$\phi_t^*(g)$) for which
\begin{equation} \label{INF}  \int_{\gamma} \dot{g} ds = 0 \;\forall \gamma
\iff  \langle \dot{\Delta} \phi_j, \phi_j \rangle = 0,\;\;
(\forall j).\end{equation} The operator $\dot{\Delta} =
Op(\dot{g})$ is the differential operator with symbol $\dot{g}.$
The linearized problem is still very difficult because it requires
a study of the asymptotic behavior of the expressions (\ref{INF})
as the lengths or eigenvalues tend to infinity. This is tantamount
to the study of the equidistribution theory of closed geodesics
and eigenfunctions.

The basic distinctions in inverse spectral theory are the
following. We say that
\begin{itemize}

\item a metric or domain is {\it spectrally determined} (within ${\cal
M}$)
if it is the unique element of ${\cal M}$ with its spectrum;

\item it is {\it locally spectrally determined} if there exists a
neighborhood of the metric or domain in ${\cal M}$ on which it is
spectrally determined;

\item a metric or domain is {\it spectrally rigid} in ${\cal M}$ if it
does not admit an isospectral deformation within the class;

\item the inverse spectral problem is solvable in ${\cal M}$
if $Spec|_{{\mathcal M}}$ is $1 -1$, i.e. if  any other metric or
domain in ${\cal M}$ with the same spectrum is isometric to it.
 If not, one has found a
counterexample.

\end{itemize}

There are analogous problems for Laplacians on non-compact
Riemannian manifolds, which often  have continuous spectra as well
as discretely occurring eigenvalues. In place of eigenvalues, one
considers the {\it resonances} $Res (\Delta)$  of $\Delta$, i.e.
the poles of the analytic continuation of its resolvent
$$R(z) = (\Delta + z^2)^{-1}.$$
Depending on whether the dimension is odd or even, $Res (\Delta)$
is a discrete subset of $\C$ or of the logarithmic plane. The
inverse spectral problems above  have natural analogues for
resonance poles. We refer to Zworski's expository articles
\cite{Zw2, Zw3}  for background.

To illustrate the current state of knowledge, we note that even
the simplest special metrics are not known to be spectrally
determined  at the present time (at least to the author's
knowledge). It is not known:

 \begin{itemize}

 \item if the standard metric $g_0$
on $S^n$ is determined by its spectrum (in dimensions $\geq 7$),
i.e. if $(M, g)$ (or even $(S^n, g)$) is isospectral to $(S^n,
g_0)$ then it is isometric to it. This has been proved in
dimensions $\leq 6$ \cite{T}.

\item  if ellipses in the plane are determined by their Dirichlet
spectra, or even if they are  spectrally rigid, i.e. if there
exist isospectral deformations of ellipses (with Dirichlet
boundary conditions).

\item if hyperbolic manifolds are determined by their spectra in dimensions $\geq
3$. I.e. if $(M_0, g_0)$ is hyperbolic  and $(M, g)$ is
isospectral to it, then is $(M, g)$  hyperbolic? This is of course
true in dimension $2$.  Is $(M, g)$ isometric to $(M_0, g_0)$? In
dimension $2$, this is known to be false for some hyperbolic
surfaces.

\item if flat metrics are determined by their spectra in the sense that
if $(M, g_0)$ is flat and $(M, g)$ is isospectral to it, then $(M,
g)$ is flat (it is known that this is true in dimensions $\leq 6$
or in all dimensions  if additionally $g$ is assumed to lie in a
sufficiently small neighborhood of $g_0$ \cite{Ku3}); it is also
classical that there are non-isometric flat tori with the same
spectra.

\end{itemize}

These special cases are tests of the strength of the known
methods. Another test is given by the two-dimensional inverse
spectral problem. One-dimensional problems are comparatively well
understood because the eigenvalue problems are ordinary
differential equations and the underlying dynamics consists of
just one orbit (an interval)!  Two-dimensional problems are
already rich in spectral and dynamical complexities, as
illustrated by the classical dynamics of twist maps or geodesic
flows on Riemannian surfaces. In general, the inverse spectral
problem grows rapidly in difficulty with the dimension, and is
already quite open for analytic surfaces and domains  in two
dimensions. This motivates our concentration on two-dimensional
problems for much of the survey. The following simple-sounding
problems are still apparently beyond the reach of known methods:

\begin{itemize}

\item Are convex analytic domains determined by their spectra among other such
domains? Are they spectrally rigid?

\item Are convex analytic surfaces of revolution determined by
their spectra among all metrics on $S^2$? Are they spectrally
rigid?

\end{itemize}

These problems are in some ways analogous to each other in  that
the unknown is a function of one variable (the boundary or the
profile curve), and that is limit  of what wave trace invariants
at one orbit can hope to recover. Surfaces of revolution are
simpler than plane domains since the geodesic flow is integrable,
while billiards on plane domains could have any dynamical type. On
the other hand, in the domain problem the class ${\cal M}$ only
consists of convex plane domains, whereas in the second problem we
allow any other metric. If we similarly restricted the class of
surfaces of the  first problem entirely to analytic convex
surfaces of revolution, then the answer is known to be `yes'
\cite{Z2}.

The inverse problem for surfaces of revolution might sound
reasonably simple since the geodesic flow is  completely
integrable, and the feeling arises that one should be able to
detect this property from the spectrum. This is one of many
problems which relate spectral invariants to dynamics of the
geodesic flow.

The relations between Laplace spectrum and dynamics have been at
the center of at least the positive results in inverse spectral
theory in the last thirty years, by comparison with the emphasis
on heat invariants in the earlier period. Before going into the
technical relations between spectral (in particular, wave)
invariants and dynamics, it might be helpful to give some
heuristic principles which suggest the relevance of dynamical
inverse spectral problems to Laplace inverse spectral problems.
The first is the relation of classical to quantum mechanics. Two
Laplacians are isospectral if
\begin{equation} \label{DIRAC} \Delta_{g_1} = U \Delta_{g_2} U^*, \end{equation} where $ U: L^2(M_1, g_1)
\to L^2(M_2, g_2)$ is a unitary operator.
 In the Dirac dictionary
of analogies, the classical analogue of this similarity is the
symplectic conjugacy
$$|\xi|_{g_1} = \chi^* |\xi|_{g_2} \iff G_{g_1}^t = \chi \circ
G_{g_2}^t \circ \chi^{-1} $$ of the corresponding geodesic flows
(cf.  (\ref{SYMP})). Here, $\chi: T^*M_1 \backslash 0 \to T^*M_2
\backslash 0$ is a homogeneous symplectic diffeomorphism.) This
analogy should not be taken too literally, but it is useful in
suggesting conjectures. In modern language, the analogue would
hold if $U$ were a unitary Fourier integral operator quantizing
$\chi$. It would clearly be  difficult to prove, even in special
cases, that isospectral Laplacians are conjugate by unitary
Fourier integral operators, though it was observed independently
by Uribe and the author (see \cite{Z6}) and by P. B\`erard
(\cite{Be, Be2}) that the Sunada counterexamples \cite{Su} have
this property (the resulting Fourier integral operators were
termed `translplantations' by B\`erard). It was also observed in
\cite{Z6} that such Fourier integral intertwining operators need
not be quantizations of symplectic diffeomorphisms, but could be
(and indeed they are, in the Sunada examples)  quantizations of
multi-valued symplectic correspondences.

But the analogy is suggestive and is fruitful  on a local (or more
accurately, formal local)  level. One of the main results in
inverse spectral theory in recent years is the theorem due to V.
Guillemin \cite{G, G3} (see also \cite{Z3, Z4}) that
isospectrality (with a simple length spectrum assumption)
 implies the  {\it formal local}  symplectic equivalence of the
 geodesic flows
around corresponding pairs of closed geodesics, i.e.  it implies
the equality of their {\it Birkhoff normal forms}. This implies
local symplectic equivalence around hyperbolic orbits, although
not around elliptic orbits (see Problem \ref{INVBNF}). It is
perhaps the closest that the above heuristic principle comes to
being valid in a general setting.

 A further
heuristic principal is that much of the dynamics of a flow is
encoded in the structure of the flow near closed orbits. To the
extent that this is true,  local equivalence would be  a powerful
fact and one should be able to obtain strong information about the
metric $g$ from studying the wave trace expansion around closed
geodesics. This  raises the possibility that the zeta functions of
the flows might determine their dynamical type (see Problem
\ref{ZETA}).

\begin{itemize}

\item Can one determine the dynamical type of the geodesic flow from the spectrum of $\Delta$, i.e.  whether the
geodesic flow is integrable, ergodic or of some other type? The
Ruelle zeta function of the geodesic flow is generically a Laplace
spectral invariant. Can one determine dynamical type from its
analytic properties.

\item In the above problems, we almost always assume in addition
that the length spectrum is  multiplicity-free: i.e. that the set
of closed geodesics of a fixed length, or more generally the fixed
point sets of the geodesic flow, should contain at most two
components interchanged by the time reversal involution $(x, \xi)
\to (x, - \xi)$. Are any of the known counterexamples, i.e.
non-isometric isospectral pairs, multiplicity free?

\end{itemize}

 The
motivation for the second problem is that wave trace methods
cannot get off the ground, and in particular the dynamical zeta
function need not be a spectral invariant,  unless the length
spectrum is simple. The issue is that there could exist
complicated cancellations among invariants of closed geodesics of
the same length.

\subsection{Strategies for solving the inverse spectral problem}

As mentioned above, we are concerned here primarily with the
positive results, ones which prove that certain geometric data is
determined by the spectra.  Roughly speaking, the strategy for
obtaining positive results has long consisted of the the following
steps:
\medskip

\noindent (A) Define  a lot of spectral invariants;

\noindent (B) Calculate them in terms of geometric or dynamical
invariants;

\noindent  (C) Try to determine the metric or domain from the
invariants.

\medskip

A crucial limitation arises in  step  (B), which accounts for the
relative paucity of positive results compared to negative results.
It is easy to define a complete set of spectral invariants, namely
the `special values' of any one of
\begin{equation}\label{TRACES}  \left\{ \begin{array}{ll} \mbox{The heat trace}, & Z(t)
= Tr e^{- t \Delta} = \sum_{j = 0}^{\infty} e^{- \lambda_j^2 t}\;\; (t > 0),\\&  \\
\mbox{The zeta function} &  \zeta(s) = Tr \Delta^{-s} = \sum_{j = 0}^{\infty}  \lambda_j^{-2 s} \;\; (\Re s > n)\\ & \\
\mbox{The wave trace} & S(t) = Tr e^{i t \sqrt{\Delta}}  = \sum_{j
= 0}^{\infty} e^{i \lambda_j t},\; \mbox{or} \; S_{ev}(t) = Tr
\cos t \sqrt{\Delta}.\end{array} \right.
\end{equation}

Of course, $\zeta(s)$ must be meromorphically continued to $\C$
and $S(t)$ is a distribution rather than a function. But  the
point we are making is that special values are rarely computable
in terms of the geometry and are therefore of limited use for
positive results on the inverse spectral problem. By comparison,
special values can be used to prove negative results by showing
that the traces of any of the above operators are the same for two
non-isometric $(M, g)$.  A key step in obtaining positive results
is to find computable invariants, and to give efficient algorithms
for computing them in terms of the simplest possible geometric
invariants. It should be mentioned that there exists spectral
invariants such as $\lambda_1$ and $\log \det \Delta = -
\zeta'(0)$ which are  useful although not computable in the above
sense. To maintain our theme of wave invariants and dynamics, we
will not discuss such invariants.

Oversimplifying a bit, the computable invariants arise from the
singularity (or related) asymptotic expansions of the traces
defined above or, in another language, from non-commutative
residues of functions of the Laplacian (we re refer to \cite{G,
G2, Z9} for discussion of non-commutative residues). In fact, all
of the computable invariants known to the author are wave
invariants, i.e. arising from the singularities of the
distribution trace or residues of the wave operator $U(t) = e^{i t
\sqrt{\Delta}}$ at times $t $ in the length spectrum of $(M, g)$
(including $t = 0$).

And it should not be forgotten that  the goal of inverse spectral
theory  is step (C). There now exist a number of rather abstract
results showing that various dynamical or quantum mechanical
invariants (e.g. Birkhoff normal forms) are spectral invariants
\cite{G, GM2, ISZ, Z1, Z3, Z4}. But this only trades one inverse
problem for another, and there is relatively little work on the
subsequent inverse problem of determining the domain or metric
from these invariants. For instance, it is not hard to see that
the classical Birkhoff normal form of the Poincare map of a
bouncing ball orbit {\it does not} determine all of the Taylor
coefficients of the boundary at the endpoints of the orbit, even
if the domain has one symmetry.  Often step (C) is the deepest,
requiring a separate study of inverse dynamical problems.

\subsection{Contents of the survey}

The focus of this survey is  on the use of wave trace formula to
derive information about the metric and geodesic flow around
closed geodesics and to determine metrics or domains from the
information.   We survey in some detail the relation between wave
invariants and Birkhoff normal form invariants on general
Riemannian manifolds, initiated by V. Guillemin, and developed by
the author and by Iantchenko-Sj\"ostrand-Zworski \cite{G, G3, ISZ,
SjZ, Z3, Z4}. In an appendix to this article \cite{SjZ2}, J.
Sj\"ostrand and M. Zworski describe in more detail how their
general results on quantum monodromy apply to the Laplacian.

But as mentioned above, to succeed with step (C) we need to be
able to determine a metric or domain from such invariants. In the
end, the crucial problem is to compute wave trace invariants in
the simplest and most efficient way,  to analyze them in detail,
and to reconstruct the domain or metric.  This is most feasible in
dimension two, so we review in some detail the articles which have
succeeded in determining special families of domains or metrics
from wave invariants, to wit, bounded plane domains \cite{CdV, P,
Z1, Z2, Z5, Z7, Z10, ISZ, GM, MM, S1} and surfaces of revolution
\cite{Z2}. One of our aims is to describe a new method for
calculating wave invariants from \cite{Z5} which so far has
achieved better results than the Birkhoff normal form approach.

In addition, we provide a fair amount of background that hopefully
puts the special problems in context.  There already exist a
number of surveys on the inverse spectral problem (e.g. \cite{Ber,
Be3, C3, Me, Gor, Gor2}) including our own expository articles
\cite{Z9, Z10}, and we have tried to avoid duplication of material
which already appears elsewhere. However, to make the survey more
self-contained we quote from a number or prior surveys, including
our own.  We also follow
 the lecture notes of Melrose \cite{Me} in our discussion of the
 Lifshits (Penrose mushroom) example of two domains with the same
 wave invariants, and also  the (much better known)
 examples of domains with the same heat invariants.

We also omit a number of topics as being too far from our focus on
wave invariants and dynamics. As mentioned above, we do not
discuss counterexamples and negative results (cf. \cite{Gor,
Gor2}). We also omit discussion of compactness results of
isospectral sets, of which there are many since the (unpublished)
work of Melrose and the work of Osgood-Phillips-Sarnak  \cite{OPS}
on isospectral sets of plane domains (for the resonance analogue,
see \cite{HZel2}). To avoid dissipation of energy, we do not
discuss the inverse resonance problem in detail, but only mention
some recent results closely related to the inverse spectral
problems covered in this survey. We refer to
 \cite{Z10} for further discussion of inverse resonance
 problems for exterior domains and to \cite{BJP, BP} for
geometric scattering settings.

The author would like to thank V. Baladi, G. Besson,  Y. Colin de
Verdi\`ere,  G. Courtois, R. Kuwabara, R. de La Llave, G. Lebeau,
M. Rouleux, K. F. Siburg and M. Zworski for informative remarks on
the contents of this survey.  Of course, errors and omissions are
the author's responsibility.

\tableofcontents

\section{Expansions at  $t = 0$}

The most classical spectral invariants are the heat invariants,
namely the coefficients of the expansion at $t = 0$ of the trace
of the heat kernel. They are closely related to the  coefficients
of the trace of the wave group at $t = 0$, although it should be
noted that the powers  $t^{-\frac{n}{2} + m}$ of the heat kernel
expansion at $t = 0$ are not singular if $m$ is even and if $m
\geq \frac{n}{2}$ (where $n = \dim M$). Hence, the heat kernel
expansion contains more information than the singularity expansion
of the wave trace at $t = 0$.

\subsection{Boundaryless case}

The earliest work in inverse spectral was based on calculations of
heat invariants in terms of curvature invariants, and the recovery
of special metrics or domains from these curvature invariants
\cite{Ber, Pa, T, T2,Ku2, Ku3}.

We recall that the heat trace expansion in dimension $n$ on a
boundaryless manifold has the asymptotic expansion,
\begin{equation} Tr e^{ t \Delta_g} \sim t^{-n/2} \sum_{j =
0}^{\infty} a_j \;  t^j. \end{equation} The coefficients $a_j$ are
the heat invariants.  We note that when $n$ is odd, the powers of
$t$ are singular and hence the expansion may be viewed as a
singularity expansion in which the terms become more regular. When
$n$ is even, the terms with $-n/2 + j < 0$ are singular but the
rest are smooth and hence are not residual. Rather one may view
the expansion as a Taylor expansion at $t = 0$ of $t^{n/2} Tr e^{
t \Delta_g} $.  But just like the singular terms, the coefficients
are spectral invariants given by integrals of curvature
invariants.

  The
first four heat invariants in the boundaryless case are given by
\cite{T, T2}
\begin{equation} \begin{array}{l} a_0 = Vol(M) = \int dVol_M \\ \\
a_1 = \frac{1}{6} \int S dVol_M \\ \\
a_2 = \frac{1}{360} \int \{2 |R|^2 \; - 2 |Ric|^2 \; + 5 S^2 ]
dVol_M \\ \\
a_3 = \frac{1}{6!} \int \{ - \frac{1}{9} |\nabla R|^2 -
\frac{26}{63} |\nabla Ric|^2 - \frac{143}{63} |\nabla S|^2 \\ \\
- \frac{8}{21} R^{ij}_{k \ell} R^{kl}_{rs} R^{rs}_{ijk\ell} -
\frac{8}{63} R^{rs} R_{r}^{j k \ell} R_{s j k \ell} + \frac{2}{3}
S |R|^2 \\ \\
- \frac{20}{63} R^{ik} R_{j \ell} R_{i j k \ell} - \frac{4}{7}
R^i_j R^j_k R^k_i - \frac{2}{3} S |Ric|^2 + \frac{5}{9} S^3\}
dVol_M. \end{array} \end{equation} Here, $S$ is the scalar
curvature, $Ric$ is the Ricci tensor and $R$ is the Riemann
tensor. In general the heat invariants are integrals of curvature
polynomials of various weights in the metric. We refer to
\cite{Ber} for background.

The heat invariants are  complicated and it is difficult to detect
meaningful patterns in the curvature polynomials. Nevertheless,
they have been successfully used to obtain inverse spectral
results.

\begin{itemize}

\item Spheres: Tanno \cite{T2} used $a_0, a_1, a_2, a_3$ to prove that the round metric $g_0$ on $S^n$ for $n \leq 6$
is determined among all Riemannian manifolds by its spectrum, i.e.
any isospectral metric $g$ is necessarily isometric to $g_0$. He
also used $a_3$ \cite{T2} to prove that canonical spheres are
locally spectrally determined (hence spectrally rigid) in all
dimensions. Patodi proved that round spheres are determined by the
spectra $\text{Spec}^0(M,g)$ and $\text{Spec}^1(M,g)$ on zero and
$1$ forms.

\item Complex projective space: Let $(M,g,J)$ be a compact
K\"ahler manifold and let $(CP^n(H),g_0,H_0)$ be a complex
$n$-dimensional projective space with the Fubini-Study metric of
constant holomorphic sectional curvature $H$. Tanno \cite{T2}
proves that  if the complex dimension $n\leq 6$ and if
$\text{Spec}(M,g,J)=\text{Spec}(CP^n(H),g_0,J_0)$, then $(M,g,J)$
is holomorphically isometric to $(CP^n(H),g_0,J_0)$. He also
proves that $(CP^n(H),g_0,H_0)$ is locally spectrally determined
in all dimensions \cite{T3}.

\item Flat manifolds: Patodi \cite{Pa} and Tanno \cite{T, T2}  used the heat invariants to prove in
dimension $\leq 5$  that if $(M, g)$ is isospectral to a flat
manifold, then it is flat. More precisely, they showed that if
$a_j = 0$ for $j \geq 1$, and if $n \leq 5$ then $(M, g)$ is flat.
The result is sharp, as  Patodi (loc. cit.) showed that $a_j = 0$
for $j \geq 1$ for the product of a 3-dimensional sphere with a
3-dimensional space of constant negative curvature. In fact, Tanno
showed that if $a_2 = a_3 = 0$, then $(M, g)$ is either
$E^6/\Gamma_1$, where $\Gamma_1$ is some discontinuous group of
translations of the Euclidean space $E^6$, or (2) $[S^3(C)\times
H^3(-C)]/\Gamma_2$, where $S^3(C)[H^3(-C)]$ is the 3-sphere
[hyperbolic 3-space] with constant curvature $C>0 [-C<0]$ and
$\Gamma_2$ is some discontinuous group of isometries of
$S^3(C)\times H^3(-C)$.
 Kuwabara \cite{Ku2, Ku3} used the invariants to prove that  flat manifolds  are locally spectrally determined,
 hence spectrally rigid.

\end{itemize}

\subsubsection{The boundary case}

When $\partial \Omega \not=0$, the heat trace has the form
\begin{equation} Tr e^{ t \Delta_g} \sim t^{-n/2} \sum_{j =
0}^{\infty} a_j \;  t^{j/2}. \end{equation} The coefficients have
been calculated for a variety of boundary conditions (see
\cite{BG} and its references).

The formulae are simplest for plane domains, where the only
invariant is the curvature $\kappa$ of the boundary. Using a
nicely adapted calculus of pseudodifferential operators, L. Smith
obtained the first five heat kernel coefficients in the case of
Dirichlet boundary conditions. They are given by :
\begin{equation} \left\{ \begin{array}{l} a_0 = \; \mbox{area of}
\; \Omega, \\ \\
a_1 = - \sqrt{\frac{\pi}{2}} \;|\partial \Omega| \; (\mbox{the
length of the boundary}), \\ \\
a_2 = \frac{1}{3} \int_{\partial \Omega} \kappa ds,\\ \\
a_3 = \frac{\sqrt{\pi}}{64} \int_{\partial \Omega} \kappa^2 ds, \\
\\
a_4 = \frac{4}{315} \int_{\partial \Omega} \kappa^3 ds, \\
\\
a_5 =  \frac{37\sqrt{\pi}}{2^{13}} \int_{\partial \Omega} \kappa^4
ds - \frac{\sqrt{\pi}}{2^{10}} \int_{\partial \Omega}
 (\kappa')^2 ds, \\
\\. \end{array} \right. \end{equation}
Here, $ds$ is arclength and $\kappa'$ is the derivative with
respect to arclength. Certain useful patterns in the heat
coefficients were used by R. B. Melrose to prove a compactness
result (later improved by Osgood-Phillips-Sarnak). We refer to
\cite{Me} for the details. In higher dimensions, one still has
\begin{equation} \label{HIGHERD} a_0 = C_n \; Vol_n(\Omega), \;\; a_1 =
C_n' \; Vol_{n-1}(\partial \Omega). \end{equation}

There exist a  few inverse  inverse spectral results using heat
invariants:

\begin{itemize}

\item Euclidean balls in all dimensions are spectrally determined among simply connected bounded Euclidean
 domains by their Dirichlet or Neumann spectra. This follows from
 (\ref{HIGHERD}) and from the fact that isoperimetric
 hypersurfaces in $\R^n$ are spheres.

 \item The exterior of the unit ball $B_3 \subset \R^3$ in dimension $3$ is uniquely determined among
 exterior domains of simply connected compact obstacles by its resonance poles. \cite{HZ}.

\end{itemize}

\subsection{Domains and metrics with the same heat invariants}

It was soon realized  that heat invariants are insufficient to
determine smooth metrics or domains. This is due to the fact that
they are integrals of local invariants of the metrics. Pairs of
non-isometric metrics with the same heat invariants can be
obtained by putting two isometric bumps. The bumped spheres will
not be isometric if the distances between the bumps are different,
but the heat invariants will be the same. There are many
variations on this well-known example. But heat invariants might
be quite useful for analytic metrics and domains, and  have also
been used in compactness results.

\section{Dynamics and dynamical inverse problems}

We now turn to the more dynamical theory of the wave group in
inverse spectral theory. The trace of the wave group expresses
spectral invariants in terms of the dynamics of the geodesic flow,
and often `reduces' inverse Laplace spectral problems to inverse
problems in dynamics. We therefore begin by recalling the relevant
dynamical notions and inverse problems.

\subsection{Geodesic flow on boundaryless manifolds}

We denote by $(T^*M, \sum_j dx_j \wedge d \xi_j)$ the cotangent
bundle of $M$ equipped with its natural symplectic form. Given a
metric $g$, we define the metric Hamiltonian in a standard
notation by
 \begin{equation} \label{MH} H(x,\xi)=|\xi|:= \sqrt{\sum_{ij=1}^{n+1} g^{ij}(x) \xi_i\xi_j} \end{equation} and
 define the energy surface to be the unit sphere bundle  $S_g^*M = \{(x, \xi)|: |\xi|_g = 1\}.$
 The spectral theorists geodesic flow is the Hamiltonian flow
 \begin{equation} G^t = \exp t \Xi_H: T^*M \backslash \to T^*M \backslash 0,\;\;\;
 \Xi_H = \; \mbox{the Hamiltonian vector field of }\; H.
 \end{equation}
It is homogeneous of degree $1$ with respect to the dilation $(x,
\xi) \to (x, r \xi), r > 0$, so nothing is lost by restricting
$G^t$ to $S^*_g M$.  We also denote its generator by $\Xi$, the
Hamiltonian and metric being understood.

The periodic orbits of $G^t$ of period $T$ of the geodesic flow
are the fixed points of $G^T$ on $S^*M$. Equivalently they are the
critical points of the length functional on the free loop space of
$M$ of length $T$.

\subsection{\label{BFLOW} Billard flow and billiard map on domains with boundary}

To define the geodesic or billiard flow $G^t$  on a domain
$\Omega$ with boundary $\partial \Omega$, we need to specify what
happens when a geodesic intersects the boundary. The definition is
dictated by the propagation of singularities theorem for solutions
of the wave equation (due to R. B. Melrose and J. Sj\"ostrand
\cite{MS}), and is therefore not  purely geometric or dynamical. A
billiard trajectory is defined to be the path along which a
singularity of a solution of the wave equation propagates. We give
a quick, informal review of the flow; for more details, the reader
might consult \cite{MS}, \cite{PS}.  For simplicity we assume that
there are no points of infinite order tangency.

We denote by $S^* \Omega$ the unit tangent vectors to the interior
of $\Omega$ and by $S^*_{in} \partial \Omega$  the manifold with
boundary of inward pointing unit tangent vectors to $\Omega$ with
footpoints in $\partial  \Omega$. The boundary consists of unit
vectors tangent to $\partial \Omega$. The billiard flow $G^t$ is a
flow on $S^* \Omega \cup S^*_{in} \partial \Omega$, defined as
follows: When an interior geodesic of $\Omega$ intersects the
boundary $\partial \Omega$  transversally, it is reflected by the
usual Snell law of equal angles. Such trajectories are called
(transveral)  reflecting rays.   The complications occur when a
geodesic intersects the boundary tangentially in $S^* \partial
\Omega$.

Convex domains are simpler than non-convex domains, since interior
rays cannot intersect the boundary tangentially. Dynamical studies
of billiards (see e.g. \cite{MF}) often restrict to convex
domains.  It should be noted that  geodesics of $\partial \Omega$
with the induced metric are important   billiard trajectories.
They are limits of `creeping rays', i.e. rays with many small
links (interior segments) which stay close to $\partial \Omega$.
In particular, the boundary of a convex plane domain is a closed
billiard trajectory. In higher dimensions, closed geodesics on the
boundary which are limits of interior creeping rays are closed
billiard trajectories. (Are all closed geodesics on the boundary
limits in this sense?)

On a non-convex domain,  a general billiard trajectory is divided
into segments which are either geodesic segments  in the interior
or geodesic segments of the boundary. Geodesic segments of the
boundary only occur where the boundary is convex. Intuitively, the
boundary segments  are limits of creeping rays along the convex
parts of the boundary. Trajectories enter and exist the boundary
at inflection points. In particular,  the boundary of a
non-convex plane domains  is not a closed billiard trajectory. If
a trajectory intersects the boundary tangentially at a
non-inflection point of a non-convex domain, it
 goes straight past the point of intersection into the interior.

When $(\Omega, g)$ is non-trapping (i.e. if there is no geodesic
ray which remains forever in the interior), the  set $S^*_{in}
\partial \Omega$ behaves like a global cross section to the
billiard flow. It is then natural to reduce the dimension by
defining the {\em billiard ball map}  $\beta: B^*(\partial \Omega)
\to B^*(\partial \Omega)$, where $B^*(\partial \Omega)$ is the
ball bundle of the boundary. We first identify $S^*_{in}
\partial \Omega \simeq B^*(\partial \Omega)$ by adding to  a tangent
(co)vector $\eta \in B^*_q \partial \Omega$  of length $<1$ a
multiple $c \nu_q$ of the inward point unit normal $\nu_q$  to
form a covector in $S^{in}_{\partial \Omega}\Omega$.  The image
$\beta(q,v)$ is then defined to be the tangential part of the
first intersection of $G^t(q, \eta + c \nu_q)$ with $\partial
\Omega.$ The billiard map is symplectic with respect to the
natural symplectic form on $B^*(\partial \Omega).$

An equivalent description of the billiard map of a plane domain is
as follows. Let $q \in \partial \Omega$ and let $\phi \in (0,
\pi)$. The point $(q, \phi)$ corresponds to an inward pointing
unit vector making an angle $\phi$ with the tangent line, with
$\phi = 0$ corresponding to a fixed orientation (say
counter-clockwise). The billiard map is then $\beta(q, \phi) =
(q', \phi')$ where $(q`, \phi')$ are the parameters of the
reflected ray at the first point of intersection with the
boundary. The map $\beta$ is then area preserving with respect to
$\sin \phi d s \wedge  d\phi$ (see e.g. \cite{MF}).

\subsection{\label{POB} Closed orbits and their Poincar\'e maps}

Closed orbits (or periodic orbits) $\gamma$ of flows are orbits of
points $(x, \xi) \in T^*M$ satisfying $G^T(x, \xi) = (x, \xi)$ for
some
 $T \not= 0$ (the period). They project to closed geodesics on the
 Riemannian manifold or domain.

We recall the definition of the nonlinear Poincare map ${\cal
P}_{\gamma}$: in $S^*M$ one forms a symplectic transversal
$S_{\gamma}$  to $\gamma$ at some point $m_0$. One then defines
the first return map, or nonlinear Poincar\"e map, $${\mathcal
P}_{\gamma}(\zeta): S_{\gamma} \to S_{\gamma}$$ by setting
${\mathcal P}_{\gamma}(\zeta) = G^{T(\zeta)}(\zeta)$, where
$T(\zeta)$ is the first return time of the trajectory to
$S_{\gamma}$. This map is well-defined and symplectic  from a
small neighborhood of $\gamma(0) = m_0$ to a larger neighborhood.
By definition, the linear Poincare map is its derivative,
$P_{\gamma} = d {\mathcal P}_{\gamma}(m_0).$

Closed geodesics are classified by the spectral properties of the
symplectic linear map $P_{\gamma}$. Its eigenvalues come in
$4$-tuples $\lambda, \bar{\lambda}, \lambda^{-1},
\bar{\lambda}^{-1}$. A closed geodesic  $\gamma$ is called:
\begin{itemize}

\item {\it non-degenerate} if $\det(I - P_{\gamma}) \not= 0$;

 \item {\it  elliptic} if all of its
eigenvalues  are of modulus one and not equal to $\pm 1$, in which
case they come in complex conjugate pairs $ e^{i \pm \alpha_j}$.

\item {\it hyperbolic} if all of its eigenvalues are real, in which case
they come in inverse pairs $\lambda_j \lambda_j^{-1}$

\item {\it loxodromic} or {\it complex hyperbolic} in the case where the
$4$-tuple consists of distinct eigenvlaues as above.

\end{itemize}

There are other possibilities (parabolic) in the degenerate case.
In the case of Euclidean domains, or more generally domains where
there is a unique geodesic between each pair of boundary points,
one can specify a billiard trajectory by its successive points of
contact $q_0, q_1, q_2, \dots$ with the boundary.  The {\em
$n$-link periodic reflecting rays} are the trajectories where $q_n
= q_0$ for some $n
> 1.$ The point $q_0, \dots, q_n$ is then a critical point of the
length functional
$${\cal L}(q_0, \dots, q_n) = \sum_{i = 0}^{n-1} |q_{i +1} -
q_i|$$
 on $(\partial \Omega)^n$.

Among the periodic orbits, bouncing ball orbits often have special
applications in inverse spectral theory.  By a {\em bouncing ball}
orbit $\gamma$ one means a periodic  $2$-link reflecting ray, i.e.
$q_0 = q_2.$  Convex domains always have at least two bouncing
ball orbits, one of which is its diameter. There exist non-convex
domains without any bouncing ball orbits. For geometric aspects of
bouncing ball orbits, we refer to \cite{Gh}.

The projection to $\Omega$ consists of a segment
$\overline{q_0q_1}$ which is orthogonal to the boundary at both
endpoints, i.e. $\overline{q_0q_1}$ is an extremal diameter. The
period is of course twice the length of the segment, which we
denote by $L$. We write $q_0 = A, q_1 = B.$ We  orient the domain
so that the links of $\gamma$ are vertical and so that the
midpoint is at the origin of $\R^2$. The top, resp. bottom, of the
boundary $\partial \Omega$ is then the graph of a function $y =
f_+(x)$ resp. $y = f_-(x)$ over the $x$-axis. It is clear that the
wave invariants depend only on the Taylor coefficients of
$f_{\pm}$ at $A, B$.

We denote by $R_A$, resp. $R_B$ the radius of curvature of the
boundary at the endpoints $A$, resp. $B$ of an extremal diameter.
The bouncing ball orbit is elliptic if $L < \min\{R_A, R_B\}$ or
if $\max\{R_A, R_B\} < L < R_A + R_B$, and is hyperbolic if $L >
R_A + R_B$ or if $\min\{R_A, R_B\} < L < \max\{R_A, R_B\}$ (see
\cite{KT}).  If $\bar{AB}$ is a local minimum diameter than $L <
R_A + R_B$ while if it is a local maximum diameter then $L > R_A +
R_B$. So a (non-degenerate) local maximum diameter must be
hyperbolic;
 a local minimum diameter is elliptic if it satisfies the additional inequalities above.

 When $\gamma$ is
elliptic, the eigenvalues of $P_{\gamma}$ are of the form
$\{e^{\pm i \alpha}\}$ while in the hyperbolic case they are of
the form $\{e^{\pm \lambda}\}.$  The explicit formulae for them
are:
\begin{equation} \label{COSALPHA}\begin{array}{l} \cos \alpha/2 =
 \sqrt{(1 - \frac{L}{R_A}) (1 - \frac{L}{R_B})}\;\;(\mbox{elliptic
case}),\\ \\ \cosh \lambda/2 = \sqrt{(1 - \frac{L}{R_A}) (1 -
\frac{L}{R_B})}\;\;(\mbox{hyperbolic case}).
\end{array} \end{equation}
We note that $f_{\pm}''(0) = \frac{1}{R_{\pm}} = \kappa_{\pm}$,
where $\kappa_{\pm}$ denotes the curvature at ($A = -, B = +)$.

\subsection{The length spectrum and the marked length spectrum}
The {\it length spectrum} of a boundaryless manifold $(M, g)$ is
the discrete set
 \begin{equation} Lsp(M,g) = \{ L_{\gamma_1} < L_{\gamma_2} < \cdots\} \end{equation}
 of lengths of closed geodesics $\gamma_j$.  In the boundary case, the length spectrum
 $Lsp(\Omega)$  is the set of lengths of closed billiard
 trajectories in the sense of \S \ref{BFLOW} and is no longer discrete, but rather has points of
 accumulation at lengths of trajectories which have intervals
 along the boundary. In the case of convex plane domains, e.g.,
 the length spectrum is the union of the lengths of periodic
 reflecting rays and multiples of $|\partial \Omega|.$
According to  the standard terminology, $Lsp(M, g)$ is  the  set
of
 distinct
 lengths, not including multiplicities, and one refers to the the
 length spectrum repeated according to multiplicity as the
 {\it extended length spectrum}.

 In the notation for $Lsp(M, g)$ we wrote $L_{\gamma_j}$ as if the
 closed geodesics of this length were isolated. But in many
 examples (e.g. spheres or flat tori), the geodesics come in
 families, and the associated length $T$ is the common length of
 closed geodesics in the family. In place of closed geodesics, one
 has
 components of the fixed point sets of $G^T$ at this time. The fixed point sets could be
 quite messy, so it is also common to assume that they are clean,
 i.e. that the fixed point sets are manifolds, and that their
 tangent spaces are fixed point sets of $d G^T$. It is equivalent
 that the length functional is Bott-Morse on the free loop space.

 The set of lengths is unformatted in the sense that one does not
 know which lengths in the list  correspond to which closed
 geodesics. A formatted notion in which lengths are assigned to
 topogically distinct types of closed geodesics is the {\it marked length spectrum}
  $ML_g$. On a manifold without boundary, it assigns to each
 free homotoply class of closed loops
 the length of the shortest closed geodesic in its class. On a
 convex plane domain with boundary, topologically distinct closed
 geodesics correspond to different rotation numbers $\frac{m}{n} =
 \frac{\mbox{winding number}}{\mbox{number of reflections}}.$ The marked length
 spectrum then associates to each rational rotation number the {\it maximal length}
 of the closed geodesics having $n$ reflection points and winding number $m$.    We
 refer to \S \ref{MLSR} (based on \cite{S1, S2}) for further  discussion and applications to
 inverse spectral theory.

\subsection{Birhoff normal forms}   Birkhoff normal forms are
approximations to Hamiltonians (or symplectic maps) near
equilibria by completely integrable Hamiltonians (or symplectic
maps). We now briefly consider the simplest kind of equilibrium
occurring for geodesic flows, a closed geodesic $\gamma$.

Let us first consider the Birkhoff normal form of the metric
Hamiltonian (\ref{MH})
 near a non- degenerate elliptic closed geodesic $\gamma$. To put
 $H$ into normal form is to conjugate it (approximately) to a
 function of locally defined action variables $(\sigma, I_1, \dots, I_n)$
 on the model space $T^*(S^1 \times \R^n)$, where $\sigma$ is the
 momentum coordinate in $T^*S^1$ and $I_j$ are transversal action
 variables which depend on the type (elliptic, hyperbolic, loxodromic)  of $\gamma$.
 More precisely, the
 normal form algorithm defines a sequence of canonical transformations $\chi_M$ at $\gamma$ which
 conjugate $H$    to the normal forms
\begin{equation} \chi_M^*H \equiv
  \sigma + \frac{1}{L} \sum_{j =1}^n \alpha_j I_j +\frac{p_1(I_1,...,I_n)}{\sigma}+...+\frac{p_M(I_1,...,I_n)}{\sigma^M}\;\;\;\;
mod\;\;\;\;\; O^1_{M + 1} \label{CBNF} \end{equation} where $p_k$
is homogeneous of order k+1 in $I_1,\dots, I_n$, and where
 $O^1_{M+ 1}$ is the space of germs
of functions homogeneous of degree 1 which vanish to order $M + 1$
along $\gamma.$ Note that all the terms in (\ref{CBNF}) are
homogenous of degree 1 in $(\sigma, I_1,...,I_n)$, and that the
order of vanishing at $|I|=0$ equals one plus the order of decay
in $\sigma$.  The coefficients of the monomials in the
$p_j(I_1,\dots,I_n)$ are known as the classical Birkhoff normal
form invariants.

The algorithm for putting a symplectic map (or Hamiltonian) into
Birkhoff normal form around a fixed point (or equilibrium point)
can be found in many places (see e.g. \cite{AKN}). For geodesic
flows, we need Birkhoff normal forms for  homogeneous Hamiltonians
around periodic orbits.
 An algorithm for putting a metric Hamiltonian into Birkhoff
normal form around a closed geodesic is described in \cite{G} and
also in  the appendix to \cite{Z3}, among other places. To carry
out the algorithm to infinite order, one needs to assume that no
eigenvalues of $\gamma$ are roots of unity, which is of course a
stronger condition than non-degeneracy.

There are also a Birkhoff normal forms for the geodesic flow and
for the Poincar\'e map. The normal form is simpler for  ${\mathcal
P}_{\gamma}$ since we have eliminated the directions along
$\gamma$ and therefore may express ${\mathcal P}_{\gamma}$ just in
terms of local action-angle variables $(I, \phi)$. The normal form
is as follows:
\begin{equation} {\mathcal P}(I, \phi) = (I, \phi + \nabla _I G_M(I)),\;\;\;mod\;\;\;\;\; O^1_{M + 1}, \label{PBNF}
\end{equation}
where $G_M(I)$ is a polynomial of degree $M$ in the $I$ variables.
Thus, to order $M + 1$, the Poincare map leaves invariant the
level sets of the actions (ellipses, hyperbolas etc. according to
the type of $\gamma$) and `rotates' the angle along them.

The well-known question arises whether the full (infinite series)
Birkhoff normal form converges and whether the formal  symplectic
map conjugating the Hamiltonian to its normal form converges. The
latter is sometimes called the Birkhoff transformation.  Clearly,
if the Birkhoff normal form diverges, then so must the Birkhoff
transformation. According to a recent article of Perez-Marco
\cite{PM}, there are no known examples of analytic Hamiltonians
having divergent Birkhoff normal forms, and we refer to that
article for further background and history on the problem. If both
the normal form and the transformation converge, the Hamiltonian
must be integrable in a neighborhood of the orbit. The Birkhoff
normal form is then the expression of the Hamiltonian in local
action-angle variables. In the generic non-resonant case, it was
proved by H. Ito that the Birkhoff transformation converges
\cite{I1} (see also \cite{I2} for the resonant case).

For non-integrable systems, the Birkhoff normal form (\ref{CBNF})
 and the
conjugating map an around elliptic orbit $\gamma$ are only
approximations in small neighborhoods of $\gamma$, which shrink at
$M $ increases. The hyperbolic case is simpler. In the case of
hyperbolic orbits of
 analytic symplectic maps in two degrees of freedom, it was proved by J. Moser \cite{Mo}
 that the Birkhoff normal form and transformation do
 converge. Convergence is more complicated in higher dimensional,
 and we refer to  Banyaga-de La Llave-Wayne
\cite{BLW}, Perez-Marco \cite{PM}  and to Rouleux \cite{R} for
recent results. In Rouleux, results of \cite{BLW} are used to
prove the existence of a local smooth symplectic  conjugacy
$\kappa^* H = q(I)$  of the metric Hamiltonian $H$ near a
hyperbolic fixed point (or orbits) to a smooth normal form $q(I)$,
obtained by Borel summation of the  formal Birkhoff  normal form
of $H$. (In the elliptic case, one only has a conjugacy $\kappa^*
H = q(I) + r$ up to a remainder $r$ which vanishes to infinite
order at the orbit.)  It follows that two metric Hamiltonians with
the same Birkhoff normal form at respective hyperbolic closed
geodesics have locally symplectically equivalent geodesic flows
near those orbits. (The author thanks M. Rouleux for corroborating
this point).

\subsection{Livsic cohomology}

The cohomology problem asks whether a function  (cocycle) $F \in
C^{\infty}(S^*M)$  satisfying $\int_{\gamma} F ds = 0$ for every
closed geodesic of the metric $g$ is necessarily a co-boundary, $F
= \Xi(f)$ where  $\Xi$ is the generator of the geodesic flow $G^t$
and $f$ is a function with some degree of regularity (see
\cite{L1} and other works of de la Llave and others for regularity
results).

It is relevant to the inverse length spectral problem for the
following reason, first observed by Guillemin-Kazhdan \cite{GK}:
Under a deformation $g_{\epsilon}$ of a metric $g = g_0$
preserving the extended $Lsp(M, g)$ (including multiplicities), one has
\begin{equation}\label{PER} \int_{\gamma} \dot{g} ds = 0,\;\;
\forall \gamma. \end{equation}  When the cohomology is trivial, one can therefore
write $\dot{g} = \Xi(f)$ for some $f$ with the given regularity.
Then one can study the harmonic analysis on $S^*M$ to see whether
the symmetric 2-tensor $\dot{g}$ can be expressed as $\Xi(f)$. The
answer is no for negatively curved surfaces, and that gave the
rigidity result of \cite{GK}. Since then their result has been
improved (see \cite{CS} for the most general result on rigidity
for negatively curved manifolds), but the general strategy retains
some  potential for other kinds of metrics.

It is known that the Livsic cohomology problem is always solvable
for an Anosov flow on a closed manifold. We refer to \cite{L1, Wa}
for recent results.  One might ask for the analogous result  for
hyperbolic billiard flows on bounded domains. It is difficult  to
formulate the analogous regularity result since the dynamics are
not smooth. For instance, the billiards are not even defined at
the corners of a domain with piecewise smooth boundary and with
concave boundary faces. One studies billiards on the domain by
puncturing out the measure zero set of orbits which ever run into
the corners. The resulting phase space is then not a closed
manifold and there would be complicated issues about regularity of
solutions as one approached the punctured set. Perhaps the
simplest setting would be that of Sinai billiards, i.e. the
exterior of a convex obstacle in a compact manifold (e.g. the
exterior of a disc in a two-dimensional torus). The billiards are
hyperbolic and there are no glancing orbits or corner orbits to
puncture out. But we are not aware of any work on the Livsic
equation in this setting.

The  Livsic equation can be studied for  for non-hyperbolic
 flows, but there do not seem to exist many studies of it
 other than manifolds all of whose geodesics are closed (in which case
 a simple Fourier analysis suffices).  One study in the near-integrable case
 is the article \cite{L2} of de la Llave.  One does not  expect the
cohomology to be trivial in general settings, but the results on
this equation might have interesting implications for the length
spectral deformation problem. As will be explained in \S
\ref{ISODEF}, in certain cases such as integrable systems, the
relevant homological equation is not the Livsic equation but an
analogue where one integrates over non-degenerate critical
manifolds of geodesics.

 \subsection{\label{DIP} Dynamical inverse problems}

 As mentioned above, inverse  spectral theory often
 proceeds by showing that certain dynamical invariants are
 invariants of the Laplace spectrum. This gives rise to inverse
 dynamical problems. In this section, we briefly survey some of the
 problems and results. A recent  survey of rigidity and conjugacy problems has been written
 by C. Croke \cite{C3}.  We are grateful to V. Baladi, G. Besson,  G. Courtois, C. Gordon, R. de La Llave,  and K. F.
 Siburg  for advice on this section.

 One of the important problems in dynamics is the following
 inverse spectral problem in dynamics.

\begin{prob}\label{MLS} For which $(M, g)$ does the marked length spectrum of $(M,
 g)$ determine $(M, g)$ up to isometry?
\end{prob}

In the positive direction,  V. Bangert \cite{B} proved that the
marked length spectrum
 of a flat two-torus determines the flat metric up to isometry.   J.
 P. Otal \cite{O}  and C. Croke \cite{C}  (see also  \cite{CFF, C3}) independently proved that the marked length spectrum
 determines surfaces of negative curvature; for more general results see \cite{C3}.   U.
 Hamenst\"adt \cite{H} proved that a locally symmetric manifold is determined by its  marked length spectrum.
An
 example due to F. Bonahon shows that
 metric structures more general than Riemannian metrics are not always determined by their marked length
 spectra.

 A
related problem concerns the conjugacy rigidity of Riemannian
manifolds. The geodesic flows  $G^t_j$ of Riemannian manifolds
$(M_j, g_j)$ are called $C^k$
 conjugate if there exists a time-preserving  $C^k$ conjugacy
 between them, i.e. a $C^k$ homeomorphism
 $\chi: S_{g_1}^*M_1 \to S_{g_2}^*M_2 $ satisfying
\begin{equation} \label{SYMP} \chi \circ G^t_1\circ \chi^{-1} = G_2^t.
\end{equation}
If $\chi$ extends to a homogeneous symplectic diffeomorphism of
the cotangent bundles, the flows are called symplectically
conjugate.  Two billiard maps are conjugate if there exists a
symplectic diffeomorphism $\chi: B^*\Omega_1 \to B^* \Omega_2$
such that
\begin{equation} \label{BILLSYMP} \chi \circ \beta_1 \circ \chi^{-1} = \beta_2.
\end{equation}

\begin{prob}\label{SYMPEQ} When does existence of a $C^k$ conjugacy between  two geodesic flows
(for a given $k$)  imply isometry of the metrics? Does symplectic
equivalence of billiard maps of convex domains imply their
isometry?
\end{prob}

In some cases, equality of marked length spectra implies $C^0$
conjugacy of geodesic flows. It was proved independently by  J. P.
Otal \cite{O}  and  by C. Croke  \cite{C} that negatively curved
surfaces with the same marked length spectrum  have $C^0$-
conjugate geodesic flows, and that such surfaces  must be
isometric. The corresponding statement in higher dimensions still
appears to be open in general. It is stated in \cite{C3} that
compact flat manifolds are $C^{\infty}$ conjugacy rigid.  In
\cite{H2}, Hammenst\"adt constructed a time-preserving $C^0$
conjugacy between negatively curved manifolds with the same marked
length spectrum. When this conjugacy is $C^1$ the work of Besson,
Courtois and Gallot \cite{BCG} proves that the manifolds must be
isometric. The method of  \cite{H} avoids this regularity issue.
The conjugacy problem for certain nilmanifolds was studied by C.
Gordon, D. Schueth and Y. Mao in \cite{GorM, GMS}. They proved for
special classes of $2$-step nilmanifolds that $C^0$ conjugacy of
the geodesic flow implies isometry. These manifolds occur in a
non-trivial isospectral deformation and have the same marked
length spectrum, so these examples also show that equality of
marked length spectra does not necessarily imply $C^0$ conjugacy
of the geodesic flows.

There is at least one known example where symplectic conjugacy
does not imply isometry, namely for surfaces all of whose
geodesics are closed. The geodesic flow of such a Zoll surface was
shown by Weinstein to be symplectically equivalent to that of the
standard $2$-sphere. Aside from related examples,  all of the
other results known to the author (see \cite{C3}) are rigidity
results showing that conjugacy implies isometry. The relations
between marked length spectral equivalence, conjugacy of geodesic
flows and isometry  do not seem to have been studied in other
settings.

One could pose local versions of the symplectic conjugacy problem,
where there should be  much less rigidity. A local version (on the
 level of germs around closed orbits) of symplectic
conjugacy is symplectic conjugacy between Poincare maps:
\begin{equation} \chi: S_{\gamma_1} \subset S^*_{g_1} M_1 \to
S_{\gamma_2} \subset S^*_{g_2} M_2, \;\;\;\; \chi {\mathcal
P}_{\gamma_1} \chi^{-1} = {\mathcal P}_{\gamma_2}. \end{equation}
One could ask:

\begin{prob}\label{LOCEQUIV}
When does local symplectic conjugacy of Poincare maps at a closed
geodesic (or local symplectic conjugacy of geodesic flows) imply
local isometry?

\end{prob}

As a special case of the local equivalence problem,
 suppose that $(M, g)$ is a hyperbolic manifold, and let
$\gamma$ be a closed geodesic. Let $g'$ be a second real analytic
metric on $M$. Suppose that there exists a closed geodesic for
$g'$ for which the Poincare maps ${\mathcal P}_{\gamma}$ are
symplectically conjugate. Must $g'$ be a hyperbolic metric?

Even less rigid  is the Birkhoff normal form, since in general
there only exists a formal (power series) canonical transformation
conjugating the germs of the metrics.

\begin{prob} \label{INVBNF} To what extent can the germ of a metric $g$ at $\gamma$ be
determined by its Birkhoff normal form at $\gamma$?

\end{prob}

 The sceptic might suspect  that when the formal conjugacy does not converge,
the Birkhoff normal form  gives little
 beyond
  a list of numerical invariants, similar to and  no simpler
than the wave invariants. When it does converge, there might be
too little rigidity to determine very much about the metric.

 For instance, it was
observed by  Y. Colin de Verdi\`ere \cite{CdV} the Birkhoff normal
form of the  billiard map around a bouncing ball orbit of a domain
with the symmetries of an ellipse determines all of the Taylor
coefficients of the domain at the endpoints of the orbit. Hence a
real analytic domain with two symmetries is determined by the
Birkhoff normal form of the Poincare map at a bouncing ball orbit.

If the domain has only one symmetry, then the calculation of
\cite{CdV} shows that  the Birkhoff normal form at a bouncing ball
orbit does not contain enough information to determine the Taylor
coefficients at both the top and bottom of the domain. This
simple, $2$-dimensional example shows how weak an invariant the
Birkhoff normal form is.

\subsubsection{Zeta functions of geodesic flows}

We add a final inverse problem which to our knowledge has not been
studied before. Let us consider two zeta functions associated to
the geodesic flow. The first is the Ruelle zeta function defined
by
\begin{equation} \label{RZETA} \log \; Z_{(M,g)}(s) = \sum_{\gamma \in {\mathcal P}} \sum_{k = 1}^{\infty}
\frac{ e^{-s k L(\gamma)}}{k | det(I - P_{\gamma}^k)|},
\end{equation}
Here, we are assuming that the metric is bumpy so that the fixed
point sets of the geodesic flow consist of isolated,
non-degenerate closed orbits. We  denote by ${\mathcal P} =
\{\gamma\}$ the set of primitive periodic orbits of the geodesic
flow $\Phi^t$ and by $P_{\gamma}$ the linear Poincare map of
$\gamma.$ As will be seen below, it is a spectral invariant as
long as the (extended) length spectrum is simple (multiplicity
free) or if the metric has no conjugate points.

The zeta function $\log \; Z_{(M,g)}(s)$  arises as the Laplace
transform of the so-called flat trace of the unitary Koopman
operator
$$W_t: L^2(S^*M) \to L^2(S^*M),\;\; W_t f(x, \xi) := f(G^t(x,
\xi)).$$ This operator does not possess a distribution trace in
the usual sense, but the integral of its diagonal part
$\int_{S^*M} W(t, \omega, \omega) dV(\omega)$ is well-defined in
the sense of pushing forward and pulling back distributions. In
this sense, it was observed by V. Guillemin \cite{G2} that
\begin{equation} Tr W(t) = C_n Vol(M, g) \delta^{(n)}(t)  +  \sum_{\gamma \in {\mathcal P}} \sum_{k = 1}^{\infty}
\frac{ \delta(t - L_{\gamma})}{ | det(I - P_{\gamma}^k)|}.
\end{equation} The Laplace transform may be viewed as the flat
trace of the resolvent of $\Xi$, so that
\begin{equation} \log \; Z_{(M,g)}(s) = Tr (\Xi - s)^{-1}.
\end{equation}

In a formal sense, the Ruelle zeta function is a  trace of a
spectral function of the geodesic flow and it is natural to wonder
how much of the spectrum of $W_t$ can be determined from the
traces. In particular, ergodicity, weak mixing  and mixing are
spectral properties of geodesic flows, i.e. can be read off from
the spectrum of $W_t$. For instance, ergodicity is equivalent to
the statement that the multiplicity of the eigenvalue $1$ for
$W_t$ equals one.

But the   trace is  not a distribution trace on $L^2$ in the sense
that if $\rho \in C_0^{\infty}$, $\int_{\R} \rho(t) W_t dt$  is
usually not a trace class operator on $L^2(S^*M)$. Hence, it is
far from clear that one can read off spectral properties of $W_t$
from these traces. At this time of writing, almost nothing seems
to be known about zeta functions except in the case of hyperbolic
flows, in which case one knows many properties of the spectrum of
the geodesic flow, e.g. that it is mixing. We refer to papers of
V. Baladi (see e.g. \cite{Ba} for her survey) for background on
spectral interpretations of zeta functions  and to C. Liverani
\cite{BKL} et. al. for some possible future extensions. C.
Deninger and others (see e.g. \cite{DS}) have recently studied the
flat trace in non-hyperbolic dynamical settings.

It might be interesting to develop the theory of zeta functions
for broader classes of flows. Indeed we are interested in the
inverse problem of determining the dynamical properties from the
trace.

A related zeta function is defined for $\Re s$ large  by
\begin{equation} \label{Rzeta} \log \; \zeta_{(M,g)}(s) = \sum_{\gamma \in {\mathcal P}} \sum_{k = 1}^{\infty}
\frac{ e^{-s k L(\gamma)}}{k | \det(I - P_{\gamma}^k)|^{1/2}},
\end{equation}
It is a Laplace spectral invariant of metrics without conjugate
points, even if (as with hyperbolic quotients) the extended length
spectrum has multiplicity.

\begin{prob} \label{ZETA} What are the analytical properties of these zeta functions
when the geodesic flow has some given dynamical signature.  Can
one determine apriori from $Z_{(M, g)}(s)$ whether the geodesic
flow is ergodic, weak mixing or mixing?
\end{prob}

The first step is to see how far the zeta function can be
analytically continued, and whether one can determine the
multiplicity of the pole at $s = 0$.  We then ask if the
multiplicity has a spectral interpretation as the dimension of the
invariant $L^2$-functions under the flow.

We observe that the zeta functions $\zeta_{(M, g)}(s)$ are Laplace
spectral invariant for metrics with multiplicity free length
spectra as well as  for those  without conjugate points. Hence
this dynamical inverse result would immediately imply that one can
determine spectral properties of the geodesic flow from the
Laplace spectrum.

\section{Wave invariants}

We now return to the Laplace spectrum and  introduce the wave
invariants, which are a more discriminating set of spectral
invariants than heat invariants. They arise from the trace of the
wave group $U(t) = e^{i t \sqrt{\Delta}}$ of $(M, g)$. One forms
the (distribution) trace
\begin{equation} Tr U(t) = \sum_{ \lambda_j \in Sp(\sqrt{\Delta})}
e^{i t \lambda_j}. \end{equation} It is a tempered distribution on
$\R$. We denote its singular support (the complement of the set
where it is a smooth function) by Sing Supp $Tr U(t)$.

The first result on the wave trace is the Poisson relation on a
manifold without boundary,
\begin{equation}\label{SS} \mbox{Sing Supp} Tr U(t) \subset \;\;
Lsp(M,g),
\end{equation}
proved by Y. Colin de Verdi\`ere \cite{CdV2, CdV3}, Chazarain
\cite{Ch2}, and Duistermaat-Guillemin \cite{DG} (following
non-rigorous work of Balian-Bloch \cite{BB2} and Gutzwiller
\cite{Gutz}). The generalization to manifolds with boundary was
proved by Anderson-Melrose \cite{AM} and Guillemin-Melrose
\cite{GM}. As above, we denote the length of a closed geodesic
$\gamma$ by $L_{\gamma}.$ For each $L = L_{\gamma} \in Lsp(M,g)$
there are at least two closed geodesics of that length, namely
$\gamma$ and $\gamma^{-1}$ (its time reversal). The singularities
due to these lengths are identical so one often considers the even
part of $Tr U(t)$ i.e. $Tr E(t)$ where $E(t)= \cos (t
\sqrt{\Delta}).$
\medskip

We emphasize that  (\ref{SS}) is only known to be a containment
relation. As will be seen below, cancellations could take place if
a length $L \in Lsp(M, g)$ is multiple, so that $Tr U(t)$ might be
smooth at $L \in Lsp(M, g)$.

\begin{prob} Is $Lsp(M, g)$ a spectral invariant? Are there any
examples where $Tr U(t)$ is $C^{\infty}$ at $t = L_{\gamma} \in
Lsp(M, g)$? \end{prob}

 No  example of this seems
to be known. Y. Colin de Verdi\`ere has  pointed out that it is
even unknown whether $Tr U(t)$ could be smooth on all of $\R
\backslash \{0\}$.  Moreover,  $Lsp(M, g)$  does not include
information about multiplicities of lengths.  Sunada-type
isospectral pairs always have multiple length spectra, and for
many (presumably, generic) examples, the length spectra have
different multiplicities. The reason is that the geodesics on the
base manifold of a Sunada quadruple can split in different ways in
the covers. Other examples of isospectral pairs (both in the sense
of eigenvalue and length spectra)  with different multiplicities
of lengths have been constructed by R. Gornet, C. Gordon and
others.

The indeterminacy of multiplicities raises  a natural problem. Let
us recall that the topological entropy of the geodesic flow of
$(M, g)$ is the exponential growth rate of the length spectrum:
\begin{equation} h_{top} : = \liminf_{L \to \infty} \log \#\{\gamma:
L_{\gamma} \leq L \}. \end{equation} Here, lengths are counted
with multipicity. In the case where geodesics come in families, we
count components of the fixed point sets.

\begin{prob} Is $h_{top}$ a spectral invariant?
\end{prob}

The question is only non-trivial when  multiplicities in the
length spectrum grow as fast as the length spectrum, as occurs for
compact hyperbolic manifolds in $\dim \geq 3$ and for arithmetic
hyperbolic quotients in dimension $2$ . As observed by
Besson-Courtois-Gallot \cite{BCG}, an affirmative answer would
show that that hyperbolic metrics are spectrally determined among
other negatively curved metrics, since they are the unique
minimizers of $h_{top}.$ A closely related problem is:

\begin{prob} (see \cite{BCG}) Suppose that $M$ with $\dim M \geq 3$
possesses a hyperbolic metric, and let ${\mathcal M}_-$ denote the
class of negatively curved metrics on $M$. Is the hyperbolic
metric $g$ determined by its spectrum among metrics in ${\mathcal
M}_-$? I.e. can there exist another non-isometric metric in this
class which is isospectral to the hyperbolic metric?
\end{prob}

\subsubsection{Singular support versus analytic singular support}
While discussing the Poisson relation, we pose the following
question concerning the analytic Poisson relation:

\begin{prob} \label{ASS}  Can one tell from  $Spec(\Delta)$  if a metric (or the
underlying manifold or domain) is real analytic? \end{prob}

The idea is to calculate the analytic wave front set of $Tr U(t)$,
i.e.
$$WF_a(\sum_{ \lambda_j \in Sp(\sqrt{\Delta})} e^{i t \lambda_j}).
$$ The analytic wave front set is the complement of the set where
the trace is real analytic.
 When $(M,g)$ is a $C^{\omega}$ (real analytic)
Riemannian manifold, the analytic wave front set $WF_a Tr U(t)$ is
also the set $Lsp(M,g)$, that is, $Tr U(t)$ is a real analytic
function outside of this discrete set. If  $(M,g)$ is a
$C^{\infty}$ but not a $C^{\omega}$ Riemannian manifold,  it is
plausible that $WF_a Tr U(t)$ could contain an interval or be all
of $\R$. Thus simply the discreteness of $WF_a Tr U(t)$ would say
that $(M, g)$ is real analytic.

\subsection{Singularity expansions}

Much more is true than the Poisson relation: $Tr U(t)$ has a
singularity expansion at each $L \in Lsp(M,g)$:
\begin{equation} \label{WTE} \begin{array}{l} Tr U(t) \equiv  e_0(t) + \sum_{L \in Lsp(M,g)} e_L(t)\;\; mod \;\; C^{\infty},
\end{array} \end{equation} where $e_0, e_L$ are Lagrangean distributions
with singularities at just one point, i.e. $sing supp e_0 = \{0\},
sing supp e_L = \{L\}$. When the length functional on the
loopspace of $M$ is a Bott-Morse functional, the terms have
complete asymptotic expansions.  In the Morse case (i.e. bumpy
metrics), the expansions take the form \begin{equation} e_0(t) =
a_{0,-n}(t+i0)^{-n} + a_{0, -n+1}(t+i0)^{-n+1}+\cdots
\end{equation}
\begin{equation}  \label{WEXP} \begin{array}{lll}
e_L(t) &=& a_{L,-1} (t-L+i0)^{-1} + a_{L,0}\log (t-(L+i0))\\[10pt]
&+&a_{L,1} (t-L+i0)\log (t-(L+i0)) +\cdots\;\;,\end{array}
\end{equation} where $\cdots$ refers to homogeneous terms of ever
higher integral degrees ([DG]). The wave coefficients $a_{0,k}$ at
$t=0$ are essentially the same as the singular heat coefficients,
hence are given by integrals over $M$ of $\int_M P_j(R,\nabla
R,...)\mbox{dvol}$ of homogeneous curvature polynomials. The wave
invariants for $t \not= 0$ have the form:
\begin{equation} \label{WINV} a_{L, j} = \sum_{\gamma: L_{\gamma}
= L} a_{\gamma, j}, \end{equation} where $a_{\gamma, j}$ involves
on the germ of the metric along $\gamma$. Here,  $\{\gamma\}$ runs
over the set of closed geodesics, and where $L_\gamma$,
$L_\gamma^{\#}$, $m_\gamma$, resp.\ $P_\gamma$ are the length,
primitive length, Maslov index and linear Poincar\'e map of
$\gamma$. (The primitive length of a closed orbit is the least
non-zero period of the orbit, i.e the length once around).  For
instance,
 the
principal wave invariant at $t = L$ in the case of a
non-degenerate closed geodesic is given by
 \begin{equation} \label{PRIN} a_{L,-1} = \sum_{\gamma:L_{\gamma}=L}
\frac{e^{\frac {i\pi }{4} m_{\gamma}
}L_\gamma^{\#}}{|\det(I-P_{\gamma})|^{\half}}. \end{equation} The
same formula for the leading singularity  is valid for periodic
reflecting rays of compact smooth Riemannian domains with boundary
and with Neumann boundary conditions, while in the Dirichlet case
the numerator must be multiplied by $(-1)^r$ where $r$ is the
number of reflection points (see \cite{GM, PS}).

The wave invariants for $t \not= 0$ are both less global and more
global than  the heat invariants.  First, they are more global in
that they are not integrals of local invariants,  but involve the
semi-global first return map ${\mathcal P}_{\gamma}$. One could
imagine different local geometries producing the same first return
map. Second, they are less global because they are determined by
the germ of the metric at $\gamma$ and are unchanged if the metric
is changed outside $\gamma$.

Thus, associated to any closed geodesic $\gamma$ of $(M, g)$ is
the sequence $\{a_{\gamma^r, j}\}$ of wave invariants of $\gamma$
and of its iterates $\gamma^r$. These invariants depend only on
the germ of the metric at $\gamma$. The principal question of this
survey may be stated as follows:

\begin{prob}\label{WI}  How much of the local geometry   of the metric $g$ at
$\gamma$ is contained in the wave invariants $\{a_{\gamma^r,
j}\}$? Can the germ of the metric $g$ at $\gamma$ be determined
from the wave invariants? At least, can the symplectic equivalence
class of its germ be determined?
\end{prob}

 As will be discussed in the next section \S \ref{BIRKHOFF},
one of the principal  results on this problem is the theorem (due
to V. Guillemin \cite{G, G2}, with some input and generalizations
by the author \cite{Z3, Z4}) that the classical Birkhoff normal
form of the metric (or the Poincar\'e map ${\mathcal P}_{\gamma}$)
at $\gamma$ is determined by the wave trace invariants. On a
global level:

\begin{prob} \label{GLOBALWI}  How much of the global geometry $(M, g)$ is
contained in the entire set of   wave invariants $\{a_{\gamma^r,
j}\}$?
\end{prob}

As these questions suggest, one may divide the potential use of
wave invariants into two classes: (i) those which use all of the
closed geodesics, and (ii) those which involve one or a few closed
geodesics. Obviously, (i) is more powerful if one can combine
information from all of the geodesics, since it adds precisely the
global feature which is lost by studying just one closed geodesic.
But it seems very difficult in general to combine information
coming from different geodesics.

To the author's knowledge, the global problem of combining  wave
invariants of all closed geodesics has only been led to successful
results on isospectral deformations. We will briefly survey the
methods and results in the next section \S \ref{ISODEF}.
Otherwise, the main results use only one or two closed geodesics,
and this  cannot possibly succeed unless the metrics or domains
are assumed real analytic. We will survey the results  in the
analytic case in \S \ref{ANALYTIC}.

Another dichotomy in the use of wave invariants is whether one
uses only `principal term' information at each geodesic, i.e. the
invariant $a_{\gamma, -1}$ (\ref{PRIN}), or whether one uses all
of the terms. Before the latter is possible, one needs to
calculate the lower order terms to some degree. We will discuss
the possible calculations in detail in later sections.

\subsection{Inverse spectral results using wave invariants}

We now consider  results which use wave invariants to specify a
class of metrics or domains which are spectrally determined in the
class or which admit no isospectral deformations. We will discuss
the proofs of some of the results in later sections to illustrate
the methods.

 The positive results based on wave invariant
analysis are as follows.

\begin{itemize}

\item Negatively curved compact manifolds are spectrally rigid  \cite{GK, C}.

\item Simple real analytic surfaces of revolution of `simple type' (with one critical
distance from the axis)  are spectrally determined
 within the class of such surfaces \cite{Z6}. Any other surface which is
 isospectral to a simple surface of revolution must be $C^0$-integrable \cite{Z7, S2}. Smooth surfaces
of revolution with a mirror symmetry thru the $x-y$ plane are
spectrally determined among metrics of this kind \cite{BH}.

 \item Simply connected analytic plane domains with two symmetry
 axes (i.e. with the symmetries of an ellipse) and with a bouncing
 ball orbit of fixed length $L$ are spectrally determined within
 this class  (\cite{Z1, ISZ, GM2}, see also \cite{CdV} for an earlier result proving spectral
rigidity of domains in this class).  A closely related result is
that  convex analytic domains with two
 symmetry axes are spectrally determined within this class
 \cite{Z1}.
 The shortest orbit is necessarily a bouncing ball orbit
and of course  its length is a spectral invariant \cite{Gh}.

 \item Simply connected analytic plane domains with one symmetry,
 and with a bouncing ball orbit whose orientation of a fixed
 length $L$ which is reversed by the symmetry, are spectrally
 determined within this class \cite{Z5}. (This implies the
 preceding result, but we state it separately since it is a new
 result based on different methods).

 \item There is a spectrally determined class of convex plane
 domains (ellipses?) for which each element is spectrally determined among all
 convex plane domains \cite{MM}. (Other isolated spectrally determined
 examples have recently been given in \cite{W}).

 \item The {\it mean minimal action} of a convex billiard table is
 invariant under isospectral deformations \cite{S1} (see \S
 \ref{MLSR} for further discussion).

 \item The exterior of a two-component analytic obstacle with two
symmetries around a bouncing ball orbit between the components is
determined by its resonances (poles of its scattering matrix)
among other such exterior domains \cite{Z10}. The  proof is
essentially the same  as in the interior case, once some known
facts on resonance poles (explained to the author by M. Zworski)
are added.

\end{itemize}

\subsubsection{Domains and metrics with the same wave invariants}

A Penrose mushroom type example  due to Michael Lifshits (see
\cite{Me},  \cite{Rau} for pictures and background) shows that
wave invariants are not sufficient to discriminate between all
pairs of smooth billiard tables. Indeed, Lifshitz constructs
(many) pairs of smooth domains $(\Omega_1, \Omega_2)$  which have
the same length spectra and the same wave invariants at
corresponding pairs of closed billiard orbits $\gamma_j$ of
$\Omega_j$ ($j = 1,2$). It follows that the Poincar\'e maps
${\mathcal P}_{\gamma_j}$ have the same Birkhoff normal forms as
well.

The idea is to exploit the complete integrability of the billiard
flow on an ellipse, i.e the fact that it is foliated by caustics.
Caustics are curves in the domain with the property that any
billiard trajectory which starts off tangent to the caustic will
remain tangent to it, as with the confocal ellipses and hyperbolae
of an ellipse.  The billiard trajectories fall into two families
separated by the bouncing ball orbit between the foci. One family
consists of trajectories which remain tangent to confocal ellipses
(which degenerate to the segment between the foci). The second
family consists of trajectories which are tangent to confocal
hyperbolae that pass through the segment between the foci. Now
divide the major axis into two parts, the segment between the foci
and the other two segments. Remove the segment between the foci
and replace it with any simple curve with the foci as endpoints,
e.g. a `tongue' below the segment. Next, remove the outer segments
of the axis between the foci and replace them with any `bumps'
below the segment. The trajectories which start in the outer bumps
never intersect the segment between the foci and therefore never
go into the tongue. Similarly trajectories which come into the
elliptical part from the tongue pass through the segment between
the foci and never go into the outer bumps.

It follows that the closed billiard orbits  fall into two
families: those which never intersect the segment between the foci
and therefore bounce back and forth between the bumps; or those
which do intersect this segment and never go into the outer bumps.

To obtain non-isometric domains with the same wave invariants it
suffices to reverse the relative orientations of the two families,
either by reflecting the segment between the foci at the center
(i.e. mirror reversing the tongue) or equivalently  by reversing
the outer bumps. Since the reversal involution is an isometry on
each `half' of the domain, it does  not change the wave invariants
for each half. But there is no `interaction' between the halves
(i.e. no closed geodesic intersecting both halves), so the
involution preserves all wave invariants. We refer to \cite{Me,
Rau} for pictures.

The question may occur whether such a domain and its reversal are
isospectral or not. This seems dubious, but we are not aware of a
proof of it. As mentioned above, equality of wave invariants
implies equality of Birkhoff normal form invariants at
corresponding pairs. Are the billiard maps of each pair
symplectically conjugate? Again dubious, but again we don't know a
proof.

\subsection{\label{ISODEF} Isospectral deformations and spectral rigidity}

One of the first uses of wave invariants and dynamics was the
proof by
 Guillemin-Kazhdan \cite{GK} that negatively curved surfaces are spectrally rigid.
In the negatively curved case, the Maslov indices are always zero
and no cancellation takes place in the wave trace formula as one
sums over closed geodesics of the same length. Hence, $Lsp(M, g)$
is a spectral invariant of negatively curved manifolds. An
isospectral deformation therefore preserves the length spectrum.

If $\gamma$ is an isolated, non-degenerate closed geodesic of $g$,
then for any deformation $g_t$ of  $g$,  $\gamma$ deforms smoothly
as a closed geodesic $\gamma_t$ of $g_t$ and one may define its
variation
\begin{equation} \dot{L_{\gamma}} = \frac{d}{dt}|_{t = 0} L_{\gamma_t}. \end{equation}
It is not hard to compute that \begin{equation} \label{DOTL}
\dot{L_{\gamma}} =  \int_{\gamma} \dot{g} ds,
\end{equation}
where $\dot{g}$ is viewed as a quadratic function on $TM$ and
$\gamma$ is viewed as the curve $(\gamma(s), \gamma'(s))$ in $TM$.

It follows that whenever the closed geodesics are non-degenerate
and of multiplicity one in the length spectrum, we have
\begin{equation}\label{GDOT}  \int_{\gamma} \dot{g} ds = 0, \;\;\; \forall
\gamma.
\end{equation}
Thus, the  Livsic cohomology problem enters. Since $Lsp(M, g)$ is
generically simple, this equation holds for generic isospectral
deformations.  We note that only the principal level of the wave
invariants was used here.

In the case of negatively curved surfaces, the geodesic flow is
Anosov and it is known that the cohomology is trivial, i.e. that
 $\dot{g} = \Xi (f)$ for some smooth $f$ (the regularity
problem is a separate issue, but for the sake of brevity we will
not consider it).  The next step is to study harmonic analysis on
the unit sphere bundle $S^*M$ to determine if there actually can
exist $f$ when $\dot{g}$ is a quadratic form. For  surfaces of
negative curvature, it was proved in \cite{GK} that there cannot
exist a smooth solution, and hence there exist no isospectral
deformations of negatively curved surfaces. After a series of
partial results by Guillemin-Kazhdan, Min-Oo and others, the
result was extended to higher dimensions by Croke-Sharafutdinov
\cite{CS}. In \cite{SU}, Sharafutdinov- Uhlmann  prove
infinitesimal spectral rigidity for the more general class of
closed two-dimensional manifold without focal points whose
geodesic flow is of Anosov type.

Isospectral deformations  also lead to the inverse marked length
spectral problem and  hence (in some cases) to the $C^0$ conjugacy
rigidity problem for the geodesic flows. An isospectral
deformation preserves lengths of closed geodesics in the generic
case where the multiplicities all equal one (or more generally
where no length is cancelled in the wave trace formula), and
therefore marks the length spectrum, i.e. gives a correspondence
between closed geodesics and lengths, by the length spectrum of
the initial metric.  More precisely, as one deforms the metric or
domain, each one-parameter family of closed geodesics
$\gamma_{\epsilon}$ stays in a fixed free homotopy class of the
fundamental group of $M$. Therefore, an isospectral deformation
gives rise to a one-parameter family of metrics or domains with
the same marked length spectrum. As discussed above, the marked
length spectrum determines the metric for surfaces of negative
curvature (Croke, Otal), which gave  a different proof of the
result of \cite{GK}.

There do not seem to exist many studies of the Livsic cohomology
equation for non-hyperbolic systems.
 It was proved by Kuwabara using heat invariants that flat tori
are spectrally rigid. It might be interesting to review the result
in terms of the Livsic equation, and possibly generalizing them to
a broader class of completely integrable systems. As will be
discussed in \S \ref{INT}, closed geodesics of integrable systems
(such as surfaces of revolution or Liouville tori) come in
families on invariant tori ${\cal T}$, and the analogue of
(\ref{GDOT}) would say that   \begin{equation} \label{FAMLIVSIC}
\int_{{\cal T}} \dot{g} dx = 0 \end{equation} for each invariant
torus of the geodesic flow. This is an analogue of the Livsic
cohomology equation when geodesics come in non-degenerate critical
manifolds. In the case of integrable systems, it  could be studied
 using Fourier analysis relative to the foliation by
invariant tori of a completely integrable system. However, it is
weaker than the Livsic equation. At least for some Riemannian
surfaces, isospectral deformations of metrics with integrable
geodesic flows must be through metrics whose geodesic flows are
$C^0$ integrable (see \S \ref{INT}), and the study of the Livsic
type equation (\ref{FAMLIVSIC}) might combine to form a more
powerful tool. The simplest case to consider is whether simple
surfaces of revolution on $S^2$ are spectrally rigid (see \S
\ref{INT} for background).

To give an extreme example of an open isospectral deformation
problem, consider the case of Zoll manifolds.   It is  simple to
see that isospectral deformations of Zoll manifolds must be Zoll
\cite{Z8}, but it is an open problem whether any non-trivial
isospectral deformations exist even for $M = S^2$. One of the main
problems is that all principal symbol level spectral invariants of
Zoll manifolds are the same, so one has to dig further into the
wave trace expansions to find obstructions to isospectral
deformability.

\subsubsection{Isospectral deformations of domains}

A few words on the analogues for Euclidean domains. To the
author's knowledge, there are no spectral rigidity results even
for Euclidean plane domains except in the case of a disc. It is
not even known whether ellipses are spectrally rigid. The length
spectrum of a bounded domain is a spectral invariant for both
Dirichlet and Neumann boundary conditions, so an isospectral
deformation must preserve the length spectrum and even the marked
length spectrum. We will discuss some results of K.-F. Siburg on
the marked length spectrum problem in \S \ref{MLSR}.

Since the length spectrum is generically an isospectral invariant,
it is natural to consider the analogue  of the Guillemin-Kazhdan
rigidity result  for domains with hyperbolic billiards. The
billiard flow is hyperbolic, so presumably the Livsic equation can
be solved (although the author is not aware of a specific
reference). It is not clear how the harmonic analysis would change
from the case of negatively curved surfaces.

There is another approach which reduces the wave group to the
boundary.
 Suppose that $\Omega_t$ is a smooth one parameter family of
smooth compact plane domains such that the Dirichlet (resp.
Neumann) spectrum $Spec(\Delta_D^{\Omega(0)})$ (resp.
$Spec(\Delta_N^{\Omega(0)})$ ) is constant in $t$. Since the area
and length of the boundary are spectral invariants, both must be
fixed under the deformation. With no loss of generality, we assume
that the variation is generated by a normal vector field $\rho
\nu$, where $\nu$ is the outward unit normal to $\Omega(0)$ and
where $\rho \in C^{\infty}(\partial \Omega(0)). $

The variations of the eigenvalues are given by  Hadamard's
variational formulae \cite{O}:
\begin{equation} \dot{\lambda}_j(0) =  \left\{ \begin{array}{ll}
\int_{\partial \Omega(0)}\; \rho \;|\partial_{\nu}
\phi_j(0)|_{\partial \Omega}|^2 d A, & Dirichlet \\
\\ \int_{\partial \Omega} \rho \; |\phi_j(0)|_{\partial \Omega(0)}|^2 \;dA, & Neumann.
\end{array}\right.
\end{equation}
Hence, the infinitesimal  deformation condition is that the right
hand sides are zero for all $j$. To normalize the problem, we
assume with no loss of generality that the deformation is volume
preserving, which implies that
\begin{equation} \int_{\partial \Omega} \rho dA = 0.
\end{equation}
Any such $\rho$ defines a volume preserving deformation of
$\Omega$.

Thus, the infinitesimal deformation is orthogonal to all boundary
traces of eigenfunctions:
\begin{equation} \dot{\lambda}_j(0) = 0 \forall j, \iff
 \left\{ \begin{array}{ll}
\int_{\partial \Omega(0)}\; \rho \;|\partial_{\nu}
\phi_j(0)|_{\partial \Omega}|^2 d A = 0, & Dirichlet \\
\\ \int_{\partial \Omega} \rho \; |\phi_j(0)|_{\partial \Omega(0)}|^2 \;dA = 0, & Neumann.
\end{array}\right.
\end{equation}
We may rewrite these conditions in terms of the boundary values of
the  wave  kernel:
\begin{equation}
E^b(t, x, x) :=   \left\{ \begin{array}{ll}
\partial_{\nu_x}\partial_{\nu_y} U(t, x, x)|_{x  \in \partial \Omega},
& Dirichlet \\  &  \\
 U(t, x, x)|_{x  \in \partial \Omega},
& Neumann\end{array}\right.
\end{equation}
as saying that
\begin{equation} \label{INTEB} \int_{\partial \Omega} E^b(t, x, x) \rho(x) dA(x)
= 0, \forall t.  \end{equation}

Using the calculations in  \cite{O}, one can obtain expansions for
(\ref{INTEB}) at  $t = 0$ similar to a heat kernel expansion in
terms of integrals of $\rho$ against  polynomials in the
(extrinsic) curvature invariants of the boundary. The
singularities at $t \not= 0$ in turn give integrals over closed
orbits of the billiard map. In the case of an ellipse, one can
directly study the boundary traces of the eigenfunctions, and it
appears  that $\int_{\Gamma} \rho d s = 0$ for every caustic
$\Gamma$ for the billiard map. It would be interesting to explore
this further.

\section{Formulae for wave invariants}

When studying the inverse spectral problem for individual pairs of
metrics  rather than isospectral deformations, it is difficult in
general   to bring in the global dynamics of the geodesic flow.
Nothing better is known than to concentrate on a single closed
geodesic. So to make progress on the inverse spectral problem, it
is crucial to have useful formulae for the wave invariants
$a_{\gamma^r, j}$  associated to a closed geodesic $\gamma$ and
its iterates $\gamma^r$. The purpose of this section is to survey
the known methods of calculation and the formulae which they
bring.

  There are several potential approaches:
\begin{itemize}

\item (i) Construct a microlocal parametrix for $e^{i t
\sqrt{\Delta}}$ at $\gamma$ and apply a stationary phase method to
calculate the wave invariants (cf. \cite{Z9}).

\item (ii) Construct a Birkhoff normal form for $e^{i t
\sqrt{\Delta}}$ at $\gamma$ and relate normal form invariants to
wave invariants (cf. \cite{G, Z1, Z3, Z4, Z9}.

\item (iii) Construct a Birkhoff normal form for the monodromy
operator ${\mathcal M}$ on the microlocal solution space
$\ker_{m_0} (\Delta - \lambda^2)$ at a fixed initial point $m_0$
of $\gamma$ and relate it to the wave invariants (cf. \cite{SjZ,
ISZ}).

\item (iv) For bounded domains, apply the Balian-Bloch or Calderon
projector methods (cf. \cite{Z5}\S \ref{BB}).

\end{itemize}

We now display the formulae and give a brief discussion of the
methods. In the following sections, we describe the methods in
much more detail.

\subsection{The parametrix method}

In boundaryless manifolds one can construct  a  Hadamard
parametrix for $\cos t \sqrt{\Delta} (x,x)$ for  small times. This
operator  solves the initial value problem
\begin{equation} \left\{ \begin{array}{ll} (\frac{\partial}{\partial t}^2 - \Delta ) u = 0& \\
u|_{t=0} = f & \frac{\partial}{\partial t} u |_{t=0} = 0
\end{array}\right.,\end{equation}
and has the form of the real part of the oscillatory integral
\begin{equation}  \int_{0}^{\infty} e^{i \theta (r^2-t^2)}
\sum_{j=0}^{\infty} U_j(x,y) \theta_{reg}^{\frac{n-1}{2} - j}
d\theta
 \;\;\;\mbox{mod}\;\;C^{\infty}  \end{equation}
where the Hadamard-Riesz coefficients $U_j$ are  determined
inductively by the transport equations
\begin{equation}\begin{array}{l}
 \frac{\Theta'}{2 \Theta} U_0 + \frac{\partial U_0}{\partial r} = 0\\ \\
4 i r(x,y) \{(\frac{k+1}{r(x,y)} +  \frac{\Theta'}{2 \Theta})
U_{k+1} + \frac{\partial U_{k + 1}}{\partial r}\} = \Delta_y WU_k.
\end{array}\end{equation} The solutions are given by:
\begin{equation}\label{HR} \begin{array}{l} U_0(x,y) = \Theta^{-1/2}(x,y) \\ \\
U_{j+1}(x,y) =  \Theta^{-1/2}(x,y) \int_0^1 s^j \Theta(x,
x_s)^{1/2} \Delta_2  U_j(x, x_s) ds
\end{array} \end{equation}
where $x_s$ is the geodesic from $x$ to $y$ parametrized
proportionately to arc-length, where  $\Theta(x, y)$  denotes the
volume density in normal coordinates centered at
 $x$, and where $\Delta_2$
operates in the second variable.

The simplest case is where the metric is without conjugate points,
in which case the parametrix is global in time  on the universal
cover.  It may then be  projected to $M$ in the usual way by
summing over the deck transformation group. One then takes the
trace and computes the coefficients by a stationary phase method.

This was first done by Donnelly \cite{D} for negatively curved
surfaces and then for all manifolds without conjugate points by
the author in \cite{Z9}. To get an impression of the complexity of
the result, here is the formula for the subprincipal wave
invariant (the coefficient of the logarithmic singularity):
$$\begin{array}{lll} a_{\gamma o} & = &   \frac{1}{|\det(I-P_\gamma)|^{\half}}\{C^o_{n,o,o,o} \int_{\gamma}
\Theta^{\half}_{\gamma}(\sigma)\{\Theta^{-\half}_{\gamma}
\mbox{Hess}(f_\gamma)^{-1}_\sigma (J_\gamma)((0,0), (\sigma,0)\\ &
& \\ & + &\mbox{Hess}(f_\gamma)^{-1}_\sigma
(\Theta^{-\half}_{\gamma}((0,0), (\sigma,0)
 d\sigma\; \\ & & \\
&+& C^o_{n,o,o,1}\int^{L_{\gamma}}_0
\int^{L_{\gamma}}_0\Theta^{\half}_{\gamma} (s\sigma_1)
\{\mbox{Hess}(f_\gamma)^{-1}_\sigma\}^2 (g_{\gamma}
\Theta^{-\half}_{\gamma} J_{[\gamma]})|_{\nu = 0}ds \\ & & \\
& + & C^o_{n,o,o,1}\int^{L_{\gamma}}_0
\int^{L_{\gamma}}_0\Theta^{\half}_{\gamma} (s\sigma_1)
\{\mbox{Hess}(f_\gamma)^{-1}_\sigma\}^3 (f_{\gamma})
\Theta^{-\half}_{\gamma}|_{\nu = 0}ds\\ & & \\ &+&  C^o_{n,o,1}
\int^{L_{\gamma}}_0 \int^{L_{\gamma}}_0\Theta^{\half}_{\gamma}
(\sigma_1) (\Delta_2\Theta^{-\half})((\sigma_1,0),
(\sigma_2,0))\cdot d\sigma_1d\sigma_2\}\;. \end{array} $$ Here,
 $\Theta_{\gamma}(x) = \Theta(x, \gamma x)$, $J_{\gamma}$ is
 a similar quantity using geodesic polar coordinates,
   $r = r(x,y)$  denote the distance function of the universal cover,
$\gamma\in \Gamma$, the fundamental group, $f_\gamma(x) = r(x,
\gamma x)^2$  denotes the displacement function. See \cite{Z9} for
details. The higher wave invariants are of course exponentially
more complicated.

In the case of bounded domains, the parametrix approach is even
messier. A microlocal parametrix for $\cos t \sqrt{\Delta}$ near a
transversal reflecting ray  was constructed by J. Chazarain
\cite{Ch} (see also \cite{GM, PS} for more details ). In principle
it could be used in the same way as a microlocal parametrix in the
boundaryless case, but in the boundary case there are almost
always conjugate points and additionally the transport equations
and phase are more complicated. Applying the method of stationary
phase to all orders to such a parametrix descends   into a jungle
of formulae. To our knowledge, no concrete inverse spectral
results have been proved using the parametrix method.

\subsection{ Wave invariants and quantum Birkhoff normal forms}

Evidently, some guiding principle is needed to civilize  the
wilderness of formulae.  Two such principles have emerged, at
least in the case of bounded domains. The first is the method of
quantum Birkhoff normal forms. As for classical Birkhoff normal
forms, one conjugates the Laplacian to a model normal form on a
model space. The wave invariants are then expressed in terms of
the coefficients of the normal form. The method is described in
detail in the next section. For the moment, we are interested in
the formula it gives for the wave invariants.

To state the result we will need some notation. The wave
invariants at a closed geodesic  $\gamma$ are  invariants of the
germ of the metric at $\gamma$. For simplicity, we assume that the
geodesic is isolated and non-degenerate $(\det (I - P_{\gamma})
\not= 0)$.  We assume $\dim M = n + 1$ and introduce   Fermi
normal coordinates $(s,y)$ along $\gamma$. We denote the
corresponding metric coefficients by $g_{ij}$, and  refer to the
coordinate vector fields $\frac{\partial}{\partial
s},\frac{\partial}{\partial y_j}$ ($j = 1, \dots, n)$ and their
real linear combinations  as Fermi normal vector fields along
$\gamma$. We denote the Riemannian connection, resp. curvature
tensor, by $\nabla$ resp. $R$ and refer   to contractions of
tensor products of the $\nabla^m R$'s with the Fermi normal vector
fields as {\it Fermi curvature polynomials}.
  Such polynomials will be called {\it
invariant} if they are invariant under the action of O$(n)$ in the
normal spaces.

We consider the Jacobi equation $Y'' + R(Y, \gamma') \gamma'= 0$
along $\gamma$. The space of complex Jacobi fields along $\gamma$
is denoted ${\cal J}_{\gamma}$. The linear Poincare map is the
monodromy map $Y(t) \to Y(t + L_{\gamma})$ on this space. We refer
to its eigenvectors $Y_j, \overline{Y}_j$ as Jacobi eigenvectors.
We denote by $y_{j k}$ their components relative to the Fermi
normal vector fields.

We then define {\it Fermi--Jacobi polynomials} to be   invariant
contractions of tensor products of $\nabla^m R$'s against
$\frac{\partial}{\partial s}$ and against the Jacobi eigenvectors,
with coefficients given by invariant polynomials in the components
$y_{jk}$. We will also use this term repeated  integrals over
$\gamma$ of such  polynomials. Finally, FJ polynomials whose
coefficients are given by polynomials in the Floquet invariants
$$\beta_j = (1- e^{i\alpha_j})^{-1}$$ are {\it
Fermi--Jacobi--Floquet} polynomials. Here, we assume for
simplicity that  the closed geodesic is elliptic and that the
eigenvalues of its Poincar\'e map are $\{e^{i \alpha_j}\}.$ (The
Floquet invariants in the general non-degenerate case are
analogous; see \cite{Z4} for details).

We define the  `weight' of the data going into  a
Fermi-Jacobi-Floquet polynomial  in terms of  its  scaling
behavior under $g \rightarrow \epsilon^2 g$. Thus, the variables
$g_{ij},\ D_{s,y}^{m} g_{ij}$ (with $m:=(m_1,\dots, m_{n+1})$),
$L:= L_{\gamma}, \alpha_j, y_{ij}, \dot{y}_{ij}$ have the
following weights: ${\rm wgt}( D_{s,y}^{m} g_{ij}) = -|m|,\ {\rm
wgt}(L) = 1, \ {\rm wgt}(\alpha_j) = 0, \ {\rm wgt}(y_{ij})=
\half, {\rm wgt}(\dot{y}_{ij}) = -\half.$    A polynomial in this
data is homogeneous of weight $s$ if all its monomials have weight
$s$ under this scaling.

The general result is \cite{Z4}:

\begin{theo}\label{BNFDATA}  Let ${\gamma}$ be a strongly
non-degenerate closed geodesic.  Then $a_{\gamma k} =
\int_{\gamma} I_{\gamma; k} (s; g)ds$ where:

\noindent(i)  $I_{\gamma; k}(s; g)$  is a homogeneous
Fermi--Jacobi--Floquet
 polynomial  of weight $-k-1$ in the data
$\{ y_{ij}, \dot{y}_{ij}, D^{m}_{s,y}g \}$ with $m=(m_1,\dots,
m_{n+1})$ satisfying $|m| \leq 2k+4$ ;

\noindent(ii) The degree of $I_{\gamma; k }$ in the Jacobi field
components is at most $6k+6$;

\noindent(iii) At most $2k+1$ indefinite integrations over
$\gamma$ occur in $I_{\gamma; k }$;

\noindent(iv) The degree of $I_{\gamma; k }$ in the Floquet
invariants $\beta_j$ is at most $k+2$.\end{theo}

The formula is simplest in dimension $2$, where there is only one
Floquet invariant $\beta$. We use the notation   $\tau$ for the
scalar curvature,  $\tau_{\nu}$ for its unit normal derivative
along $\gamma$, $\tau_{\nu \nu}$ for  the Hessian ${\rm
Hess}(\tau)(\nu,\nu)$. We denote by
 $Y$  the unique
normalized Jacobi eigenvector along $\gamma$ and by  $\dot{Y}$ its
time-derivative.

The subprincipal  wave invariant $a_{\gamma 0}$ is then given by:
\begin{equation} a_{\gamma 0} = \frac{a_{\gamma, -1}}{L^{\#}}[ B_{\gamma 0;4} (2\beta^2 - \beta
-\frac{3}{4}) + B_{\gamma 0;0} ]\end{equation} where:

\noindent (a) $a_{\gamma, -1}$ is the principal wave invariant
(\ref{PRIN});

\noindent(b) $L^{\#}$ is the primitive length of $\gamma$;
$\sigma$ is its Morse index; $P_{\gamma}$ is its Poincar\'e map;

\noindent(c) $B_{\gamma 0; j}$ has the form:
$$B_{\gamma 0;j} = \frac{1}{L^{\#}} \int_o^{L^{\#}} [a\; |\dot{Y}|^4 +
 b_1\; \tau |\dot{Y}\cdot Y|^2 + b_2\; \tau {\rm Re\,} (\bar{Y} \dot{Y})^2
+c \;\tau^2 |Y|^4 + d \;\tau_{\nu \nu}|Y|^4 + e\; \delta_{j0} \tau
]ds  $$
$$ +\frac{1}{L^{\#}}\sum_{0\leq m,n \leq 3; m+n=3}
C_{1;mn}\;\frac{\sin ((n-m)\alpha)}{|(1-e^{i(m-n)\alpha})|^2}\;
\left |\int_o^{L^{\#}}  \tau_{\nu}(s)\bar{Y}^m\cdot Y^n(s)ds
\right |^2 $$
$$+\frac{1}{L^{\#}}\sum_{0\leq m,n \leq 3; m+n=3} C_{2;mn}\;{\rm Im\,}
\; \left\{\int_o^{L^{\#}} \tau_{\nu}(s)\bar {Y}^m\cdot Y^n(s)
\left[\int_o^s\tau_{\nu}(t)\bar {Y}^n\cdot Y^m(t)dt\right] ds
\right\}$$ for various universal  coefficients.

The formula is very complicated, and it is only the second term in
the expansion! That is why one hopes to directly use the
information that these coefficients determine the Birkhoff normal
form.

In the case of surfaces of revolution, it is possible to simplify
the formulae to the point where one can determine the metric along
the invariant geodesics. The point is to recover the Taylor
expansion of the curvature and its normal derivatives along
invariant geodesics.  In  \cite{Z2}, a variant of this method was
used to show that `simple' surfaces analytic of revolution (those
with just one invariant geodesic) are determined by their wave
invariants, hence are spectrally determined among such surfaces.
The variant was to use the existence of a global Birkhoff normal
form for the metric and Laplacian, which greatly simplified the
calculations.

The wave invariants can also be calculated by the method of
Birkhoff normal forms in the boundary case, but the method becomes
appreciably more difficult. We therefore use a different approach
which is inspired by the Balian-Bloch approach to the Poisson
formula.

\subsection{Balian-Bloch approach in the  boundary case}

We now consider wave invariants associated to periodic reflecting
rays of  a bounded plane domain $\Omega$. For simplicity and
because they are often useful in applications, we consider wave
invariants at a a $2$-link periodic reflecting ray $\gamma$ of
length $ r L_{\gamma}$ (a `bouncing ball orbit'). We refer to  the
 \S \ref{POB} for the basic definitions.

 We  introduce  some further notation: we
write $y =f_{\sigma(j)}(x_j)$ where  $\sigma (2j + 1) = -,
\sigma(2j) = +.$ The length functional in Cartesian coordinates
for a given assignment $\sigma$ of signs is given by
\begin{equation} {\mathcal L}_{\sigma} (x_1, \dots, x_{2r}) =
 \sum_{j = 1}^{2r - 1} \sqrt{(x_{j + 1} - x_j)^2 + (f_{\sigma(j + 1) }(x_{j +1}) - f_{\sigma(j)}(x_j))^2}. \end{equation}
We will  need formulae for the entries of its  Hessian $H_{2r}$ in
Cartesian coordinates at the critical point corresponding to the
$r$th repetition of the bouncing ball orbit with a given
orientation. Thus, $x_j = 0$ for all $j$. We will assume, with no
essential loss of generality, that $q(x_{odd}) = A, q(x_{ev}) =
B.$ We put: $a = 1 - L f_+''(0), b = 1 - L f_-''(0)$.

The Hessian is given in either angular or Cartesian  coordinates
by: \begin{equation}  \label{KT}
 H_{2r} = \frac{1}{L} \left\{ \begin{array}{lllll} a & 1 & 0 & \dots & 1  \\ & & & & \\
  1 & b & 1 & \dots & 0 \\ & & & & \\ 0 & 1 &  a & 1 & 0  \\ & & & & \\ 0 & 0 & 1 & b & 1 \dots
 \\ & & & & \\\dots & \dots & \dots & \dots & \dots \\ & & & & \\
1 & 0 & 0 & \dots & b \end{array} \right\}, \end{equation} where
  $ a = 2(1- \frac{L}{R_A}),\;\;\; b =  2( 1- \frac{L}{R_B}).$
A well known formula relates the determinant of this Hessian to
that of the Poincar\'e map:
\begin{equation} \label{HESSPOIN} \det (I - P_{\gamma^r}) =
L^{2r} \det H_{2r};\;\;\; (\gamma \;\; 2-\mbox{link}.)
\end{equation}

We denote the matrix elements of the inverse Hessian $H_{2r}^{-1}$
at the bouncing ball orbit  by $h^{i j}.$ To be precise, the
bouncing ball orbit has two possible orientations, one (which we
denote by $+$ which starts at the top graph  $y = f_+$ and
proceeds to the bottom $y= f_-$ and the other $-$ which in the
reverse order. We then have two (closely related) length functions
and Hessians. We denote the matrix elements of their inverses by
$h_{\pm}^{i j}$. One has $h_-^{pq} = h_+^{p - 1, q - 1}. $

We  now state the formulae for the wave invariants.

\begin{theo} \label{BGAMMAJ} \cite{Z5}  Let $\gamma$ be a $2$-link reflecting ray of length $
L_{\gamma}$, and let $\gamma^{r}$, resp. $\gamma^{-r}$ be the
$r$th iterate of $\gamma, $ resp. $\gamma^{-1}$.
 Then there exist polynomials $p_{2, r, j}(\xi_1, \dots, \xi_{2j}; \eta_1, \dots, \eta_{2j})$ which are
 homogeneous of degree $-j$ under the dilation $f \to \lambda f,$ which are  invariant under the substitutions
$\xi_j \iff \eta_j$ and under $f(x) \to f(-x)$, and which have
degree $j + 1$ in the Floquet data $e^{i \alpha r},$ such that
$$a_{\gamma^r, j}  = p_{2, r, j}(f_-^{(2)}(0), f_-^{(3)}(0), \cdots,
 f_-^{(2j + 2)}(0); f_+^{(2)}(0), f_+^{(3)}(0), \cdots, f_+^{(2j + 2)}(0)).$$
The leading order term in derivatives of $f_{\pm}$ has the form
$$ \begin{array}{l}
 r \{(h^{22}_{+})^{j}  \sum_{q = 1; q \equiv 0 }^{2r} h^{2q}_{+} + (h^{22}_{+})^{j - 2} \sum_{q; q \equiv 0 }^{2r} (h^{2q}_{+})^3\} f_+^{(2j - 1)}(0) f_+^{(3)}(0) \\ \\+
r \{(h^{11}_{+})^{j }  \sum_{q; q \equiv 1 }^{2r}  h^{1q}_{-} +
(h^{11}_{+})^{j - 2}  \sum_{q; q \equiv 1 }^{2r} (h^{1q}_{+})^3\}
f_-^{(2j - 1)}(0) f_-^{(3)}(0)
\\ \\ -  r\{ ( h^{22}_+)^{j - 1} h^{11}_+
\sum_{q = 1; q \equiv 1}^{2r}  h^{2q}_{\pm} + (h^{22}_{+})^{j - 2}
\sum_{q= 1; q \equiv 1}^{2r} (h^{2q}_{+})^3 \}f_+^{(2j - 1)}(0)
f_-^{(3)}(0)\\ \\ - r\{(h^{11}_{+})^{j - 1} h^{22}_{+}  \sum_{q =
1; q \equiv 0  }^{2r} h^{1q}_{\pm} + (h^{11}_{+})^{j - 2} \sum_{q
= 1; q \equiv 0 }^{2r}(h^{1q}_{+})^3 \}f_-^{(2j - 1)}(0)
f_+^{(3)}(0)
\\ \\ + R_{2r} (j^{2j - 2} f_+(0), j^{2j - 2} f_-(0)),\end{array} $$
where the remainder $ R_{2r} (j^{2j - 2} f_+(0), j^{2j - 2}
f_-(0))$ is a polynomial in the designated jet of $f_{\pm}.$

\end{theo}

In the case where the domain has a symmetry interchanging the top
and bottom, so that $f_{\pm}$ are mirror images of the graph of $y
= f(x)$. The length functions  for $\gamma^r$ and $\gamma^{-r}$
(and hence their Hessians) are equal, so   we may drop the
subscripts $\pm$. The formula simplifies as follows:

\begin{cor} \label{BGAMMAJSYM} Suppose that $\gamma$ (as above) is invariant under an isometric involution $\sigma.$ Then, modulo
the error term  $R_{2r} (j^{2j - 2} f(0))$, we have:
$$ \begin{array}{l}  a_{\gamma^r, j - 1} =  r \{ 2 (h^{11})^j f^{(2j)}(0) + \{2  (h^{11})^j \frac{1}{2 - 2 \cos \alpha/2} +  (h^{11})^{j - 2} \sum_{q = 1}^{2r} (h^{1 q})^3\} f^{(3)}(0) f^{(2j - 1)}(0)\}\}.\end{array} $$
Here,   we sum over repeated indices.
\end{cor}

The inverse spectral problem is then reduced to analyzing Hessian
coefficients. It is not evident to the author how this reduction
of the inverse spectral problem would be visible using Birkhoff
normal forms.

\subsubsection{Melrose-Marvizi invariants}

In \cite{MM}, Melrose-Marvizi introduce further spectral
invariants of a convex domain which could  be interpreted as
quantum normal form invariants, namely the normal form of $\Delta$
around the closed geodesic $\partial \Omega$. The calculation in
\cite{MM} involves Melrose's normal form for glancing
hypersurfaces and may be viewed as a normal form construction, but
one might also follow Lazutkin's construction of whispering
gallery quasi-modes to put $\Delta$ into normal form as a function
of a so-called Airy operator. This would presumably give a new way
to calculate the invariants. J. Toth and the author have partial
results in this direction (unpublished), and from discussions with
G. Popov it appears that they are known to others.

 The first two
Melrose-Marvizi integrals have the form \cite{MM}:
$$I_1 = - 2 \int_{\partial \Omega} \kappa^{2/3} ds,\;\; 1080  I_2
=\int_{\partial \Omega} (9 \kappa^{4/3} + 8 \kappa^{-8/3}
\dot{\kappa}^2) ds.$$

\section{\label{BIRKHOFF}Calculation of wave invariants I: Birkhoff normal forms}

We now outline the construction of the Birkhoff normal form of
$\Delta$ in both the boundaryless and boundary cases, and indicate
how it is used to calculate the wave invariants. Our goals are to
explain how the calculations in Theorem \ref{BNFDATA} are done.

Let $\gamma$ be a non-degenerate  closed geodesic on an
$n$-dimensional Riemannian manifold, and at first let us  assume
it to be elliptic. During the 70's, various authors (Babic,
Lazutkin, Ralston, Guillemin-Weinstein) constructed a series of
quasi-modes and quasi-eigenvalues associated to $\gamma$. Roughly
speaking, the results showed that for each transversal quantum
number $q \in \Z^{n-1}$, there was an approximate eigenvalue of
the form
\begin{equation} \label{APPROXEIG} \lambda_{kq} \equiv r_{kq} +
\frac{p_1(q)}{r_{kq}} + \frac{p_2(q)}{r_{kq}^2} +
...\end{equation} where
$$r_{kq} = \frac{1}{L} (2 \pi k + \sum_{j=1}^n (q_j + \frac{1}{2}) \alpha_j)$$
where the coefficients are polynomials of specified degrees and
parities. We refer to \cite{BB} for a clear exposition and for
details.

A natural question is whether the wave invariants $a_{\gamma j}$
of (\ref{WEXP}) can be determined from the quasi-eigenvalues
(\ref{APPROXEIG}), i.e. from the polynomials $p_j(q)$, which could
be computed explicitly from the quasimode constructions described
in \cite{BB}.  This question was difficult to answer or even to
formulate precisely until the article of V. Guillemin \cite{G} on
quantum Birkhoff normal forms for the Laplacian around elliptic
closed orbits. This article does not refer to quasi-modes or
quasi-eigenvalues but it effectively proves that the wave
invariants may indeed be expressed in terms of the polynomials
$p_j(q)$ and conversely that these polynomials are spectral
invariants. A somewhat different proof of this result, which
constructed quantum Birkhoff normal forms by adapting the
algorithm in \cite{BB} for constructing quasi-modes, was given by
the author for elliptic orbits in \cite{Z2}, and for general
orbits in  \cite{Z3}.

The results of \cite{G} are stated in the terms of quantum
Birkhoff normal forms.  A number of expositions now exist which
explain this notion in detail (see for instance \cite{Z9} in
addition to the original articles), so we will only briefly review
the notion. To put $\Delta$ into normal form, is first  to
conjugate it
  into a distinguished
maximal abelian algebra ${\mathcal A}$ of pseudodifferential
operators on a model space, the cylinder
 $S^1_L \times \R^n$, where $S^1_L$ is the circle of length $L$.
  The algebra is  generated by the tangential operator
$D_s:=\frac{\partial}{i \partial s}$ on $S^1_L$ together with the
transverse action operators. The nature of these action operators
depends on $\gamma$. When $\gamma$ is elliptic, the action
operators are harmonic oscillators
$$I_j=I_j(y,D_y) := \frac{1}{2} (D_{y_j}^2 + y_j^2),$$
while in the real  hyperbolic case they have the form
$$I_j = y_j D_{y_j} + D_{y_j} y_j.$$
When $\gamma$ is non-degenerate,  they involve some mixture of
these operators (and also complex hyperbolic actions)  according
to the spectral decomposition of $P_{\gamma}$.
 For notational simplicity we restrict to the elliptic case and
 put
$$H_{\alpha}:= \frac{1}{2}\sum_{k=1}^n \alpha_k I_k $$
where $e^{\pm i \alpha_k}$ are the eigenvalues of the Poincare map
$P_{\gamma}$.

To put $\Delta$ into normal form is to conjugate it to the model
space and algebra as a function of $D_s$ and the action operators.
The conjugation is only defined in a neighborhood of $\gamma$ in
$T^*M-0$, i.e. one constructs  a microlocally elliptic
  Fourier Integral operator $W$
from the conic neighborhood  of $\R ^+\gamma$ in $T^*N_{\gamma}-0$
to a conic neighborhood of $T^*_{+}S^1_L$ in $T^*(S^1_L \times
R^n)$ such that:

\begin{equation} \label{QBNF} W \sqrt{\Delta_{\psi}}W^{-1}
 \equiv {\mathcal D}  +
\frac{\tilde{p}_1(\hat{I}_1,\dots,\hat{I}_n)}{L {\mathcal D}} +
 \frac{\tilde{p}_2(\hat{I}_1, \dots,\hat{I}_n)}{(L {\mathcal D})^2}
+\dots+\frac{\tilde{p}_{k+1}(\hat{I}_1,\dots,\hat{I}_n)}{(L
{\mathcal D})^{k+1}}+ \dots \end{equation} where the numerators
$p_j(\hat{I}_1,...,\hat{I}_n),
\tilde{p}_j(\hat{I}_1,...,\hat{I}_n)$ are polynomials of degree
j+1 in the variables $\hat{I}_1,...,\hat{I}_n$,  where $W^{-1}$
denotes a microlocal inverse to $W$. Here, ${\mathcal D} =  D_s +
\frac{1}{L}H_{\alpha}.$ The kth remainder term lies in the space
$\oplus_{j=o}^{k+2} O_{2(k+2-j)}\Psi^{1-j}$, where $\Psi^s$
denotes the pseudo-differential operators on the model space of
order $s$ and where  $O_{2(k+2-j)}$ denotes the operators whose
symbols vanish to the order $2(k + 2 - j)$ along $\gamma$. We
observe that (\ref{QBNF}) is an operator version of
(\ref{APPROXEIG}).
\bigskip

The inverse result of \cite{G}, see also \cite{Z3, Z4}] is:
\bigskip

\noindent{\bf Theorem  } {\it Let $\gamma$ be a non-degenerate
 closed geodesic. Then the quantum Birkhoff normal form
around $\gamma$ is a spectral invariant; in particular the
classical Birhoff normal form is a spectral invariant.}
\bigskip

In other words, one can determine the polynomials
$p_j(\hat{I}_1,...,\hat{I}_n)$ from the wave trace invariants of
$\Delta$ at $\gamma.$ When $\gamma$ is a non-elliptic  (e.g.
hyperbolic) geodesics, the action operators have continuous
spectra and there are no approximate eigenvalue expansions as in
(\ref{APPROXEIG}). Thus the quantum Birkhoff normal form approach
is  conceptually clearer.

We give a brief review of how wave invariants and normal form
invariants for the Laplacian are calculated by the algorithm in
\cite{Z2, Z3, Z4,Z5} in the case of an elliptic closed geodesic.
The algorithm is similar in the general non-degenerate case.  We
then consider the algorithm for calculating the normal form of the
monodromy operator in \cite{ISZ}.

\subsection{Manifolds without boundary}

Our algorithm for conjugating $\Delta$ to normal form along a
closed geodesic $\gamma$ is inspired by the constructions due to
Lazutkin \cite{L} and Babich-Buldyrev \cite{BB} for constructing
quasimodes associated to elliptic closed orbits. The same
conjugation method works in the general non-degenerate case.

Since we have recently written an exposition of the method
\cite{Z9}, we only include only a few formal aspects of the
calculations, in the hope  they provide sufficient explanation of
how the   wave invariants are calculated. In particular, we wish
to emphasize:
\begin{itemize}

\item The {\it semiclassical normal form}, and how it arises from a
  a semi-classical scaling of the
Laplacian along $\gamma$. Our method differs from the others
\cite{G, ISZ} by working entirely on the quantum level and
inductively on the Taylor expansion of the metric around $\gamma$.
The homogeneous normal form (\ref{QBNF}) is obtained from the
semiclassical normal form by (roughly speaking) replacing the
large parameter by $|D_s|$ along $\gamma$.

\item The homological equations and the obstructions to their
solvability along a closed orbit. The obstructions determine the
normal form.

\end{itemize}

\subsection{Quasi-mode heuristics}

We wish to conjugate the Laplacian to a normal form on a model
space,  namely   the normal bundle $N_{\gamma}$ of $\gamma$, or
equivalently,   the cylinder $S^1_L \times \R^n$, where as above
$S^1_L = \R / L\Z$. On the phase space level, the model is
$T^*(S^1_L \times \R^n)$ or  more precisely the cone $|\eta| \leq
\epsilon \sigma, |y| \leq \epsilon$ in the natural symplectic
coordinates $(s,\sigma,y,\eta)$ corresponding to Fermi normal
coordinates along $\gamma$.

 In the construction of  quasi-modes associated to $\gamma$, one uses the WKB ansatz
$$\Phi_{kq}(s,\sqrt{r_{kq}}y) = e^{ir_{kq}s} U_{kq}(s,
\sqrt{r_{kq}}y,r_{kq}^{-1}),$$ where $r_{kq} = \frac{2\pi}{L}(k +
(q + \frac{1}{2}) \alpha)$ and where
$$ U_{kq}(s, \sqrt{r_{kq}}y,r_{kq}^{-1})
 \sim  \sum_{j=0}^{\infty} r_{kq}^{-\frac{j}{2}}
U_q^{\frac{j}{2}}(s, \sqrt{r_{kq}}y,r_{kq}^{-1}). $$ One then
solves asymptotically the eigenvalue
$$\Delta_y e^{ir_{kq}s} U_{kq}(s, \sqrt{r_{kq}}y,r_{kq}^{-1}) \sim
\lambda_{kq}e^{ir_{kq}s}
 U_{kq}(s, \sqrt{r_{kq}}y,r_{kq}^{-1}),$$
 satisfying the periodicity condition of being well defined on the
 cylinder.
 The  intertwining operator $W_{\gamma}$ to  normal form may be thought of as the
 operator taking the model eigenfunctions to the quasimodes,
$$W_{\gamma} \phi_{kq}(s,y) = \Phi_{kq}(s, \sqrt{r_{kq}}y).$$
We now change our point of view and concentrate on the
construction of the intertwining operator.

\subsection{Semiclassically scaled Laplacian}

In view of the form of the quasi-eigenvalue problem,
$$\Delta_u e^{\frac{i}{hL}s}U(s, h^{-\half} u,h) = \lambda(h) e^{\frac{i}{hL}s}U(s,
h^{-\half}u,h),$$ it is natural to rescale the Laplacian before
conjugating it to normal form. In fact, this rescaling is all we
retain from the quasi-mode construction. We therefore introduce
the  unitary operators $T_h$ and $M_h$ on ${\cal H}_T$ or
equivalently on the 1/2-density version $L^2_T(\R^1\times \R^n,
\Omega_{1/2})$ given by
$$T_h (f(s,u)|ds|^{1/2}|du|^{1/2}):= h^{-n/2} f(s, h^{-\half}u) |ds|^{1/2}|du|^{1/2} \leqno (2.5a)$$
$$M_h(f(s,u)|ds|^{1/2}|du|^{1/2}) := e^{\frac{i}{hL}s} f(s,y)|ds|^{1/2}|du|^{1/2} \leqno (2.5b).$$

\noindent{\bf Definition} {\it The rescaling of an operator $A_u =
a(s,D_s,u,D_u)$ of the adapted model  is the operator
$$A_h := T_h^*M_h^*AT_hM_h.$$}

We now rescale the Laplacian
 in  Fermi normal coordinates. It is convenient to first conjugate
 to the unitarily equivalent  1/2-density Laplacian
 $$\Delta_{1/2} := J^{1/2} \Delta J^{-1/2},$$
which can be written  in the form:
$$-\Delta_{1/2} = J^{-1/2}\partial_s g^{oo}J \partial_s J^{-1/2}
+\sum_{ij =1}^{n} J^{-1/2}\partial_{u_i} g^{ij} J \partial_{ u_j}
J^{-1/2}\leqno(2.9)$$
$$\equiv g^{oo}\partial_s^2 + \Gamma^o \partial_s +
 \sum_{ij=1}^n g^{ij} \partial_{u_i}\partial_{u_j} + \sum_{i=1}^{n} \Gamma^{i}
\partial_{u_i} + \sigma_o.$$
Here, $\partial_{x} := \frac{\partial}{\partial x},$ and  $
J=J(s,u)= \sqrt{g}$ is the volume density in these coordinates.

We then have:
$$-M_h^* \Delta M_h = -(hL)^{-2}g^{oo} + 2i(hL)^{-1}g^{oo} \partial_s + i(hL)^{-1}\Gamma^o +
\Delta \leqno (2.10)$$ Conjugation with $T_h$ then gives
$$-\Delta_{h} =  -(hL)^{-2} g^{oo}_{[h]} + 2i(hL)^{-1}g^{oo}_{[h]}\partial_s + i(hL)^{-1}\Gamma^o_{[h]}+
 h^{-1}( \sum_{ij=1}^n g^{ij}_{[h]}\partial_{u_i}\partial_{u_j}) + h^{-\half}(\sum_{i=1}^{n} \Gamma^{i}_{[h]}
\partial_{u_i}) + (\sigma)_{[h]},\leqno(2.11)$$
 the subscript $[h]$ indicating to dilate the coefficients of the operator in the form,
$f_h(s, u):=f(s, h^{\half} u).$

Expanding the coefficients in Taylor series at $h=0$, we obtain
the asymptotic expansion
$$\Delta_h \sim \sum_{m=0}^{\infty} h^{(-2 +m/2)}{\mathcal L}_{2-m/2} \leqno (2.12)$$
where ${\mathcal L}_2 = L^{-2},$ ${\cal L}_{3/2}=0$ and where
$${\mathcal  L}_1 = 2 L^{-1}[i  \frac{\partial}{\partial s} + \half \{\sum_{j=1}^{n}
\partial_{u_j}^2 -
\sum_{ij=1}^{n} K_{ij}(s) u_i u_j\}] .$$

\subsection{\label{SCNF} Conjugation of scaled Laplacian to semi-classical
normal form}

We now conjugate the semi-classically scaled Laplacian to normal
form. The conjugating operators are $s$-dependent operators acting
on the transverse space to $\gamma$ (they are sometimes called
tangential semi-classical Fourier integral operators).

 The first step is to conjugate the
principal term into quadratic normal form on the transverse space.
Since it is quadratic,  there exists a metaplectic conjugation
$\mu(s)$ depending on $s$ on the transverse space which conjugates
it to the operator  ${\cal D}$ of (\ref{QBNF}). That is, for each
$s \in [0, L]$ there exists an operator $\mu(s)$ in the
metaplectic representation of the metaplectic group $ML(n, \R) \to
Sp(n, \R)$ which acts on the transverse $\R^n$ of the model space.
For background on metaplectic operators and details on the
conjugation we refer to \cite{Z3, Z4}.  This conjugates the other
operators ${\mathcal L}_j$ to new operators ${\mathcal D}_j$.

Once  the principal term is in normal form, we continue by
conjugating  the lower order terms to functions of the actions by
using perturbation theory, i.e. formal series in $h$. Thus we wish
to
  put the formal series $\Delta_h$
into a {\it semi-classical normal form} by an infinite sequence of
conjugations, each one putting one new order in $h$ into normal
form. We carry out the procedure to two orders to illustrate the
main points.

The first step is to
  construct  $Q_{\frac{1}{2}}(s, x, D_x)$
 such that
$$e^{-ih^{\frac{1}{2}}Q_{\frac{1}{2}}} {\mathcal D}_h e^{ih^{\frac{1}{2}}Q_{\frac{1}{2}}}|_o
=[- h^{-2}+ 2 h^{-1}{\mathcal D} + {\mathcal D}^{\half}_o +
\dots]|_o$$ where \begin{equation} \label{SLASH} |_o \;\;
\mbox{denotes the restriction of the operator to functions of}\;\;
y,
\end{equation} and where the dots $\dots$ indicate higher powers in
$h$. Introduction of the  $|_o$ operator is motivated by the
construction of quasi-modes, since the normal form is only being
applied to the amplitude (a function of $y$). A more conceptual
explanation (suggested by the work of
Iantchenko-Sj\"ostrand-Zworski \cite{ISZ}) is that the correct
Hilbert space on which to define the operators is the microlocal
solution space of $\Delta - \lambda^2$ along $\gamma$. These are
essentially the space of quasi-modes along the open arc $(0, L)$
of $\gamma$. Thus,  $|_o$ restricts the operator to the microlocal
solution space.

 Expanding the exponential, we find that the
 operator $Q_{\frac{1}{2}}$ then must satisfy
the {\it homological equation}
$$\{ [L^{-1}{\mathcal D},Q_{\frac{1}{2}}]+ {\mathcal D}_{\frac{1}{2}} \}|_o= 0.$$
  One may solve this equation explicitly by further conjugating ${\mathcal D}$ to $D_s$.
  When  ${\mathcal D} = D_s$, the  homological equation  becomes
$$\{ [L^{-1} D_s, Q_{\frac{1}{2}}]+
{\mathcal D}_{\frac{1}{2}} \}|_o= 0,$$ that is,
$$L^{-1}\partial_s \; Q_{\frac{1}{2}}|_o =
 - i\{{\mathcal D}_{\frac{1}{2}}\}|_o.$$
Here,  $\partial_s A$ is the Weyl operator whose complete symbol
is the $s$-derivative of that of $A$. To solve this  equation  we
rewrite it in terms of  complete Weyl symbols. We will use the
notation $A(s,x,\xi)$ for the complete Weyl symbol of the operator
$A(s,x,D_x)$.  We then arrive at the homological equation
$$L^{-1}\partial_s \tilde{Q}_{\half}(s,x,\xi)= -i {\mathcal D}_{\half}|_o(s,x,\xi)
$$
We solve with the Weyl symbol
$$\tilde{Q}_{\half}(s,x,\xi) = \tilde{Q}_{\half}(0,x,\xi) + L \int_0^s
-i {\cal D}_{\half}|_o(u,x,\xi)du$$ where
$\tilde{Q}_{\half}(0,x,\xi)$ is determined by the consistency
condition
$$\tilde{Q}_{\half}(L,x,\xi) - \tilde{Q}_{\half}(0,x,\xi) =
L \int_0^L -i {\cal D}_{\half}|_o(u,x,\xi)du.)$$

To solve the equation, we invoke the fact (which is not obvious)
that ${\mathcal D}_{\half}|_o(u,x,\xi)$ is a polynomial of degree
3 in $(x,\xi)$. We also switch to complex coordinates $z_j = x_j +
i\xi_j$ and $\bar z_j = x_j - i \xi_j$ in which the action of
$r_{\alpha}(L)$ is diagonal. The homological equation becomes
$$\tilde{Q}_{\half}(0,e^{i\alpha}z, e^{-i\alpha}\bar z) -
\tilde{Q}_{\half}(0,z,\bar z) = L\int_0^L -i {\cal
D}_{\half}|_o(u,z,\bar z ) du.$$ We put:
$$\tilde{Q}_{\half}(s,z,\bar z) = \sum_{|m|+|n|\leq 3} q_{\half;mn}(s) z^m \bar z^n$$
and
$$ {\mathcal D}_{\half}|_o(s,z,\bar z ) du = \sum_{|m|+|n|\leq 3} d_{\half;mn}(s)
z^m \bar z^n$$ then the homological equation becomes
$$\sum_{|m|+|n|\leq 3} (1 - e^{(m - n) \alpha}) q_{\half; mn}(0) z^m \bar z^n
= -i L^2 \sum_{|m|+|n|\leq 3} \bar d_{\half;mn} z^m \bar z^n.$$
Since there are no terms with $m=n$ in this (odd-index) equation,
and since the $\alpha_j$'s are independent of $\pi$ over $\Z$,
there is no obstruction to the solution. Thus we can simply
eliminate the term  ${\mathcal D}_{\half}$.

\subsubsection{The normal form coefficients}

We now consider the second step, where the first obstruction
occurs.  We thus seek a pseudodifferential operator
$\tilde{Q}_1(s,x,D_x)$ and a quadratic polynomial
$f_o(I_1,...,I_n)$ in the action operators so that
$$e^{-i h \tilde{Q}_1}{\mathcal D}^{\half} e^{i h \tilde{Q}_1} = h^{-2}L^{-2} + h^{-1}L^{-1}D_s +
h^{-\half}{\mathcal D}^{\half}_{\half} + {\mathcal D}^1_o(s,D_s,
x,D_x) + \dots$$ with
$${\mathcal D}_o^1(s,D_s,x, D_x)|_o = f_o(I_1,...,I_n).$$
As usual, the dots signify terms of higher order in $h$. The
homological equation is then
$$\{[D_s,\tilde{Q}_1] + {\mathcal D}_o^{\half}\}|_o = f_o(I_1,...,I_n),$$
 or equivalently
$$\partial_s\tilde{Q}_1|_o =\{-{\mathcal D}_o^{\half} + f_o(I_1,...,I_n)\}|_o.$$

We rewrite the equation in terms of the complete Weyl symbols, to
obtain
$$L^{-1} \partial_s \tilde{Q}_1(s,z, \bar z) = -i \{{\mathcal D}^{\half}_o|_o (s,z,\bar z) -
f_o(|z_1|^2,\dots, |z_n|^2)\} $$ or equivalently
$$\tilde{Q}_1(s,z,\bar z) = \tilde{Q}_1(0,z,\bar z) -i L\int_0^s
[{\mathcal D}^{\half}_o|_o (u,z,\bar z) - f_o(|z_1|^2,\dots,
|z_n|^2)] du,  $$ or again,
$$\tilde{Q}_1(0,e^{i\alpha}z,e^{-i\alpha}\bar z) - \tilde{Q}_1(0,z,\bar z) = -i L \{\int_0^L
{\mathcal D}^{\half}_o|_o (u,z,\bar z)du  - L f_o(|z_1|^2,\dots,
|z_n|^2) \}.$$ By construction,  ${\mathcal D}^{\half}_o|_o
(u,z,\bar z)$ is a polynomial of degree 4, so if we put
$$ \tilde{Q}_1(s,z,\bar z) = \sum_{|m|+|n|\leq 4} q_{1;mn}(s) z^m \bar z^n, \;\;\;\;\;\;
 f_o(|z_1|^2,\dots, |z_n|^2) = \sum_{|k|\leq 2} c_{o k} |z|^{2k} $$
and
$${\mathcal D}^{\half}_o|_o (s,z,\bar z)du :=\sum_{|m|+|n|\leq 4} d_{o; mn}^{\half}(s) z^m \bar
z^n,\;\;\;\;\;\;\;\;\;\bar d^{\half}_{o;mn}:=
\frac{1}{L}\int_o^Ld^{\half}_{o;mn}(s)ds, $$  we can solve for the
off-diagonal coefficients,
$$q_{1; mn}(0) = -i L^2 (1 - e^{i (m-n)\alpha})^{-1} \bar d_{o; mn}^{\half}.$$
We also must set the diagonal coefficients equal to zero, and this
determines the  $c_{o k}$ coefficients:
$$c_{o k} =   \bar d_{1; kk}^{\half}. $$
These coefficients are the quantum normal form coefficients. Since
we used only algebraic operations on the rescaled Laplacian, it is
clear that the coefficients can be calculated in terms of the
 data described in Theorem \ref{BNFDATA}.
 We refer to \cite{Z4, Z9} for further details.

\subsection{Bounded domains} There are two serious differences in the
normal form for the Laplacian around a bouncing ball orbit of a
bounded domain.  First, we need to straighten out the domain to
define a model space. Secondly, we need the conjugation to normal
form to incorporate the boundary conditions.

\subsubsection{The model domain $\Omega_o$}

The  configuration space of the model is the infinite strip $
\Omega_o$. We denote the coordinate on $[0,L]$ by  $s$  and that
on $\R$ by  $y$ , with dual cotangent coordinates $\sigma, \eta$
on $T^* [0,L] \times T^* \R.$
  As in the boundaryless case,  we would like to view the model
space as  the normal bundle of the orbit.  In the case of a
bouncing ball orbit in the  boundary case, the normal bundle and
exponential map are ill-defined at the reflection points but
Lazutkin has  constructed a nice replacement for them.  Namely, he
constructs a formal power series map  from $\Omega_{\epsilon}$ to
$\Omega_0$,
$$ \Phi : \Omega_{\epsilon}
\rightarrow U_{\epsilon},\;\;\;\;\; \Phi(s, y) = (\s, \y) $$ of
the form
$$ \begin{array}{l} \s = s + \sum_{m = 2}^{\infty}
\kappa_{m m} (s)   y^m \\ \\
\y = \sum_{p=1}^{\infty}  \psi_{p p}(s)   y^p \end{array}$$ with
real valued analytic  coefficients in $s$ extending analytically
to a neighborhood of $[0, L]$ in $\C$  and satisfying the boundary
conditions:
$$ \{\s = 0\} \cup
\{\s = L\} = \Phi (\partial \Omega_{\epsilon}). $$

\begin{lem} Suppose that $\overline{AB}$ is a bouncing ball orbit. Then there exists a transversal power series  map $\Phi : (\Omega_{\epsilon}, \partial
\Omega_{\epsilon}) \rightarrow (\Omega_o, \partial \Omega_o)$
which straightens the domain and puts  $\bDelta$ in the form:
$$\begin{array}{ll} \bar{\Delta} \sim D_{\s}^2 + B(\s, \y)(\y^2 D_{\s}^2 + D_{\y}^2) +
\Gamma'_s D_s + \Gamma'_y D_y & \rm{elliptic \; case} \\ & \\
\bar{\Delta} \sim D_{\s}^2 + B(\s, \y)( D_{\y}^2 - \y^2 D_{\s}^2)
+ \Gamma'_s D_s + \Gamma'_y D_y & \rm{hyperbolic \; case}
\end{array}
$$
in the sense that the left and right sides agree to infinite order
at $y = 0.$ Here, $B(s,y)$ is a transversal power series.
\end{lem}

By choosing $\Phi$ carefully, one can arrange that the Laplacian $
(\Phi^{* -1}\Delta \Phi^*)_{\half}$ in the transformed coordinates
has the  principal term  $D_{\s}^2 + (a_{22}(s) \y^2 D_{\s}^2 +
e_{11}(s) \y D_{\y} + b_{00}(s) D_{\y}^2.$ We  can further put the
principal part into the normal form
\begin{equation}\begin{array}{ll}   D_{\s}^2  +      \half \dot{b}_{00}(s)[\y^2 D_{\s}^2 +  D_{\y}^2] & \rm{elliptic \;
case} \\ & \\D_{\s}^2  +      \half \dot{b}_{00}(s)[  D_{\y}^2 -
\y^2 D_{\s}^2 ] & \rm{hyperbolic \; case}\end{array} .
\end{equation}

\subsection{Semi-classical scaling of the Laplacian}

Having straightened the domain, and hence having transferred the
information about the boundary into the metric, we now  follow the
approach in the boundaryless case by scaling the Laplacian.
Following \cite{L, Z2}, we use the notation $N = h^{-1}$.

 We define
operators $T_N, M_N$ on the model space $L^2(\Omega_o)$ by
\begin{equation} \begin{array}{ll} \bullet &  T_N f(\s,\y):= N f(\s, N \y) \\
\\ \bullet &  M_N f(\s,\y) := e^{i N^2 \s} f(\s,\y) \end{array} \end{equation}

\begin{defn} The semiclassically scaled Laplacian is the operator
on $\Omega_o$ defined by
$$\bDelta_N = M_N^* T_N^* \bar{\Delta} T_N M_N .$$ \end{defn}

In the straightened form, the rescaled  Weyl symbol of
$\bar{\Delta}$ has the  form:
\begin{equation}\sigma_{\bar{\Delta}_N}^w
 \sim   (\sigma + N)^2 + B(s, \frac{1}{N} y) (\pm y^2 N^{-2}(\sigma + N^2)^2 + N^2 \eta^2)
+ K(s, \frac{1}{N} y).  \end{equation}

\subsection{\label{SCNFB} Conjugation to semiclassical normal form}

As in the boundaryless case, we begin by conjugating the principal
term (which again is quadratic)
\begin{equation} \sigma_{\bar{\Delta}_N}^w  =  N^4 + 2 N^2 [  D_{\s}
+    \dot{b}_{00}(s) \hat{I}] \;\;\rm{mod}\;\;N
\end{equation}
into quadratic normal form by an $s$-dependent metaplectic
conjugation on the transverse space to the bouncing ball orbit.
Here, $\hat{I}$ denotes the quantum action operator:  $\hat{I}^e =
\half (D_{\y}^2 +  y^2)$ in the elliptic case and  $\hat{I}^h =
\half (D_{\y}^2 -  y^2)$ in the hyperbolic case.

 In \cite{Z2},  we prove that there
exists an $SL(2, \R)$-valued function $a_{\alpha}(\s)$ so that
$$\begin{array}{l} \mu(a_{\alpha})^* [ D_{\s} + b_{00}(\s) \hat{I}] \mu(a_{\alpha}) =
 D_s + \frac{\alpha}{L} \hat{I},\;\;\; \mu(a_{\alpha})(0) = \mu(a_{\alpha}) (L) =
 Id,
\end{array}$$
 where $\alpha = \int_0^L b_{00}(s) ds.$
We then conjugate  the resulting operator
\begin{equation} {\cal R}_N : = \mu(a_{\alpha})^* \bDelta_N \mu(a_{\alpha}) \end{equation}
 to the semiclassical
normal form
\begin{equation} \label{SCNFR} F_N(\hat{I})^2 \sim N^4 + N^2 \frac{\alpha \hat{I}}{L} + p_1(\hat{I}) +
N^{-2} p_2(\hat{I}) + \cdots \end{equation}
 by conjugations as in the boundaryless case, but now additionally
 preserving the boundary conditions.

We assume for simplicity that the bouncing ball orbit is elliptic,
but the same  argument and result hold in the hyperbolic case. We
again use the notation of (\ref{SLASH}): we denote by $A|_o$ the
restriction of $A$ to functions of $y$ only. In the boundary case,
we need to construct complex valued  Weyl symbols $P_{j/2},
Q_{j/2}$ so that iterated composition with $e^{ N^{-j}(P + i
Q)_{j/2}^w}$ will successively remove the lower order terms in
${\cal R}_N$ after restriction by $|_o$ and so that the boundary
condition is satisfied. It turns out that the boundary condition
on the conjugating operator involves only the odd part of the real
part and the even part of the imaginary part:
\begin{equation} \label{BC} \begin{array}{l}  P_{j/2}^o(0, \y, \bareta) = P_{j/2}^o(L, \y, \bareta) = 0,\;\;\;\;\;\;
 Q_{j/2}^e(0, \y, \bareta) = Q_{j/2}^e(L, \y, \bareta) = 0. \end{array} \end{equation}
We emphasize that there is no condition on  the odd part
$Q_{j/2}^o$ of the imaginary part or the even part  $P^e_{j/2}$ of
the real part.

\subsubsection{A Sturm-Liouville homological equation}

As before, we see what happens in the first and second steps. We
first find  $P_{\frac{1}{2}}(\s, \y, D_{\y}),\;
Q_{\frac{1}{2}}(\s, \y, D_{\y}) $ so that the boundary conditions
are satisfied and so that
\begin{equation} e^{-  N^{-1}(P + i  Q)_{\frac{1}{2}}} {\cal R}_N e^{ N^{-1}
(P+ i Q)_{\frac{1}{2}}}|_o = \{ N^4 + N^2 [D_s + \frac{\alpha}{L}
\hat{I}]+ \cdots \}|_o.
\end{equation} Expanding the exponential, we get
the {\it homological equation}:
\begin{equation}  \{ [D_s + \frac{\alpha}{L} \hat{I} , (P+ i Q)_{\frac{1}{2}}]+ {\cal R}_{\frac{1}{2}} \}|_o= 0. \end{equation}
Taking the complete symbol of both sides we get the symbolic
homological equation:
\begin{equation} i \{\sigma + \frac{\alpha}{L} I, P_{\half} +  i Q_{\frac{1}{2}}\}+ {\cal R}_{\frac{1}{2}}|_o = 0.\end{equation}
The equation may be rewritten in the form:
\begin{equation} \label{SOL} \partial_{\bar{s}} (P + i Q)_{\frac{1}{2}}(\s, r_{\alpha}(\s)(\y, \bareta))=
 -i  {\cal R}_{\frac{1}{2}}|_o, \end{equation}
whose solution is given by
\begin{equation} (P + i Q)_{\half}(\s,r_{\alpha}(\s)(\y, \bareta))  =   (P + i Q)_{\half}(0) -
 i \int_0^{\s}  {\cal R}_{\half}(u, \y, \bareta) du.\end{equation}
We thus  need to determine $ (P + i Q)_{\half}(0)$ so that the
boundary conditions at $\s = 0$ and $\s = L$ are satisfied.  Thus,
we must solve the system
\begin{equation} \label{EO} \left\{\begin{array}{ll} \partial_{\s} P_{\half}^e + \frac{\alpha}{L}  \{I, P^o_{\half}\} =  \Im {\cal R}_{\frac{1}{2}}^e|_o &
\partial_{\s} P_{\half}^o + \frac{\alpha}{L}  \{I, P^e_{\half}\} =  \Im {\cal R}_{\frac{1}{2}}^o|_o
\\ & \\
\partial_{\s} Q_{\half}^e + \frac{\alpha}{L}  \{I, Q^o_{\half}\} = -\Re  {\cal R}_{\frac{1}{2}}^e|_o &
\partial_{\s} Q_{\half}^o + \frac{\alpha}{L}  \{I, Q^e_{\half}\} = - \Re {\cal R}_{\frac{1}{2}}^o|_o \end{array}\right. \end{equation}
 with the boundary conditions
\begin{equation} P_{\half}^o(0) = P^o_{\half}(L) = 0,\;\;\;\;\;\;Q_{\half}^e(0) = Q_{\half}^e(L) = 0.\end{equation}

As in the boundaryless case, ${\mathcal R}_{\half}^{\pm}(u, \y,
\bareta ) =
  \y \circ (a r_{\alpha}^{-1} (u)) I$ is a   polynomial of degree 3 in $(\y, \eta)$ in
  which every term is of odd degree in $(\y, \eta)$.   It follows that $P_{\half}, Q_{\half}$
   are  also odd polynomials
of degree 3. Assuming we are in the elliptic case, we  change
coordinates to
 the complex cotangent variables $z = \y + i \bareta,
\bz = \y - i \bareta$  and write
\begin{equation}\left\{\begin{array}{l}P_{\half}^e(\s, z, \bz) = \sum_{m,n: m + n \leq 3} p_{\half mn}^e(\s) (z^m \bz^n +
\bz^m z^n), \\ \\
P_{\half}^o(\s, z, \bz) = \sum_{m,n: m + n \leq 3} p_{\half
mn}^o(\s) (z^m \bz^n - \bz^m z^n).
\\ \\
Q_{\half}^e(\s, z, \bz) = \sum_{m,n: m + n \leq 3} q_{\half
mn}^e(\s) (z^m \bz^n +
\bz^m z^n), \\ \\
Q_{\half}^o(\s, z, \bz) = \sum_{m,n: m + n \leq 3} q_{\half
mn}^o(\s) (z^m \bz^n - \bz^m z^n).
\end{array} \right.\end{equation}

Then (\ref{EO}) may be rewritten in terms of these coordinates. In
the elliptic case, we have:
\begin{equation} \label{EOC} \left\{\begin{array}{l}
\frac{d}{ds} p^e_{\half mn}(s) +  i \frac{\alpha}{L}  (m - n) p^o_{\half mn}(s)  =   \Im {\cal R}_{\frac{1}{2}}^e|_o \\ \\
 \frac{d}{ds} p^o_{\half mn}(s) + i \frac{\alpha}{L} (m - n)  p^e_{\half mn}(s)  =  \Im  {\cal R}_{\frac{1}{2}}^o|_o \\ \\
 \frac{d}{ds} q^e_{\half mn}(s) +  i \frac{\alpha}{L}  (m - n) q^o_{\half mn}(s)  = - \Re {\cal R}_{\frac{1}{2}}^e|_o \\ \\
 \frac{d}{ds} q^o_{\half mn}(s) + i \frac{\alpha}{L} (m - n)   q^e_{\half mn}(s)  = - \Re {\cal R}_{\frac{1}{2}}^o|_o
 \end{array}\right. \end{equation}
 Since $m \not= n$, we
can eliminate  $p^e_{\half mn},  q^o_{\half mn}$ and reduce the
homological equation to (uncoupled) second order Sturm-Liouville
boundary value problems
\begin{equation} \label{2O} \left\{ \begin{array}{l} - \frac{d^2}{d s^2}p^{o}_{\half mn}(s) -  [\frac{\alpha}{L} (m - n)]^2 p^{o}_{\half mn}(s) =
 \frac{L}{\alpha (m-n)}  \frac{d}{ds} \Im  {\cal R}_{\frac{1}{2}}^e|_o +   \Im {\cal R}_{\frac{1}{2}}^e|_o
\\ \\
- \frac{d^2}{d s^2}q^{e}_{\half mn}(s) -  [\frac{\alpha}{L} (m -
n)]^2 q^{e}_{\half mn}(s) = - \frac{L}{\alpha (m-n)}  \frac{d}{ds}
\Re {\cal R}_{\frac{1}{2}}^e|_o +  \Re  {\cal
R}_{\frac{1}{2}}^e|_o \\ \\  p^o_{\half mn}(0)  =  0, \;
p^o_{\half mn}(L)  =  0; \;\; q^e_{\half mn}(0)  =  0, q^e_{\half
mn}(L)  =  0 \end{array} \right.
\end{equation}
for the independent variable $p^o_{\half mn}, q^e_{\half mn}(s)$.
 The  boundary value problem (\ref{2O}) is always
solvable  unless $0$ is an eigenvalue of the operator $ D_{\s}^2 -
[\frac{\alpha}{L}  (m - n)]^2$ (elliptic case).   The
eigenfunction would have to have the form $\sin (\frac{\alpha}{L}
(m - n) \s)$, hence a sufficient condition for solvability is that
$\alpha/\pi \notin  \Q$. The analogous boundary problem in the
hyperbolic case is always solvable.

Thus, the  conjugation eliminated the sub-principal term
${\mathcal R}_{\frac{1}{2}}$ and put  our operator in  normal form
up to order $N$.  The exponents $P_{\half}, Q_{\half}$ are odd
polynomial differential operators of degree 3 with smooth
coefficients defined in a neighborhood of $[0, L]$.

\subsubsection{The normal form coefficients emerge}

 We  carry the process forward one
more step because, as in the boundaryless case, the even steps
produce the non-trivial normal form coefficients. In outline, the
conjugated operator in the second step, ${\cal R}^{\half}_N$, is
further conjugated with an exponential of the form $e^{N^{-2}  (P
+ i Q)_1}.$ The homological equation is again a Sturm-Liouville
boundary problem, with the new feature that diagonal terms with $m
= n$ do occur in (\ref{EOC}). In the elliptic case, a diagonal
term is a function of $|z|^2$, hence is even under the involution
$z \to \bz$.
 We write the diagonal terms in the form:
\begin{equation}\left\{ \begin{array}{l} P_{j}^d(\s, |z|^2) = \sum_{m: m \leq j} p_{j m}^d(\s) |z|^{2m} \\ \\
Q_{j}^d(\s, |z|^2) = \sum_{m: m \leq j} q_{j m}^d(\s) |z|^{2m}
\end{array}\right. \end{equation}

It is impossible to solve the inhomogeneous  boundary problem
(\ref{EOC})  for the diagonal terms as in the odd case with zero
boundary conditions. To obtain a solvable system, we must add new
terms $f_j(|z|^2)$ to the right side:
\begin{equation} \left\{ \begin{array}{l}
\frac{d}{ds} p_{j m}^d(s)   =   \Im {\cal R}_{2 - j }^d|_o - \Im f_j(|z|^2) \\ \\
 \frac{d}{ds} q^d_{jm}(s)  = - \Re {\cal R}_{2 - j}^d|_o + \Re f_j(|z|^2)  \\ \\
 p_{j m}^d(0) = p_{j m}^d(L) = 0\\ \\  q^d_{jm}(0) = q^d_{jm}(L) = 0. \end{array}
 \right. \end{equation}
 The right hand sides are determined by the condition that
equations with
\begin{equation}\left\{ \begin{array}{l} p^{d}_{j m}(s)   = \int_0^{s} \{ -\Im  {\cal R}_{2 -j}^{j, d}|_o(u, |z|^2) + \Re f_j (|z|^2)\} \\ \\
q^{d}_{j m}(s)   = \int_0^{s}
 \{  \Re {\cal R}_{2-j}^{j, d}|_o(u, |z|^2) + \Im  f_j (|z|^2)\}  du, \end{array}\right. \end{equation}
satisfies the boundary conditions. This forces
\begin{equation} f_j(|z|^2) = -  \frac{1}{L} \int_0^{L}  {\cal R}_{2-j}^{j, d}|_o (u, |z|^2) du. \end{equation}

In the case $j = 1$, we have then conjugated ${\cal R}_N^+$ to
$N^4 + N^2 {\cal R} + Op^w(f_1(I)) + O(N^{-1}).$ The coefficients
of the  polynomial $f_1(I)$ are the first Birkhoff invariants.
We note that $Op^w(f_1(I))$ is a function of $\hat{I}$, so we have
conjugated to normal form to fourth order. We now repeat the
process to complete the normal form.

\subsection{Monodromy operator}

We now outline a different approach to normal forms due to
Iantchenko-Sj\"ostrand-Zworski \cite{SjZ, ISZ}. (See also the
appendix \cite{SjZ2}).  Their goal is to construct the normal form
of a so-called {\it quantum monodromy operator} rather than of the
Laplacian or wave group. The normal form of the monodromy operator
is apparently the same as the semi-classical normal form of the
wave group that we constructed in \S \ref{SCNF} and \S
\ref{SCNFB}. However, it is presented in a conceptually clearer
way which works equally well in the boundary or boundaryless case.
On the other hand, the method presented above also gives an
algorithm for calculating the normal form, which is at present
lacking in the approach of \cite{ISZ}.

The new concepts introduced in \cite{SjZ, ISZ}
 {following  earlier one-dimensional definitions in works of
Helffer-Sj\"ostrand and Colin de Verdiere- Parisse) are the
following:

\begin{itemize}

\item The monodromy operator ${\mathcal M}(\lambda): \ker_{m_0} (\sqrt{\Delta} - \lambda) \to
\ker_{m_0}(\sqrt{\Delta} - \lambda)$  acting on the space
$\ker_{m_0}(\sqrt{\Delta} - \lambda) $ of microlocal solutions of
the equation $ (\sqrt{\Delta} - \lambda)  u(h) =
O(\lambda^{-\infty}) \;\;\; \mbox{near}\;\; m_0$, where  $m_0$ is
 is an arbitrarily chosen base point of the  given geodesic $ \gamma$.

\item The {\it flux norm} on the microlocal solution space $\ker_{m_0}(\sqrt{\Delta} -
\lambda)$, with respect to which ${\mathcal M}(\lambda)$ is
unitary;

\item The Grushin reduction of the wave operator to the monodromy
operator and the expression of the wave trace as a trace of the
monodromy operator.

\end{itemize}

The microlocal solution space is essentially the space of
quasi-modes along arcs of $\gamma$,  so  the monodromy approach is
not in essence so different from the approach we outline above. We
give a brief exposition of the key ideas, which are  discussed in
more detail in the  appendix by Sj\"ostrand-Zworski \cite{SjZ2}.

Assuming for simplicity that $\gamma$ projects to an embedded
curve in $M$, we again work on the model space $S^1 \times
\R^{n-1}$ given by the normal bundle $N_{\gamma}$ to $M$. We then
consider the universal cover $\R \times \R^{n-1}$,  fix the base
point $m = (0, 0),$ and   consider microlocal solutions along the
geodesic $\R \times \{0\}$.  In the simply connected space  $\R
\times \R^{n-1}$, there is no obstruction to constructing a global
quasimode along $\R \times \{0\}$. So we may define
\begin{equation} {\mathcal M}(\lambda): \ker_0 (\sqrt{\Delta} - \lambda) \to
\ker_L (\sqrt{\Delta} - \lambda) \end{equation} where $L =
L_{\gamma}$ by
\begin{equation}  {\mathcal M}(\lambda) [u_{\lambda}]_0 = [u_{\lambda}]_{(L,0)} \end{equation}
where $[u_{\lambda}]_x$ is the germ of the quasimode at the point
$x$. We obtain a quasi-mode when $(I - {\mathcal M}(\lambda))
[u_{\lambda}]_0 = 0$.

To make contact with the wave group, we quote the following
observation of Sj\"ostrand-Zworski:

\begin{prop} \label{MW} \cite{SjZ2} We have:
$${\mathcal M}(\lambda) = \exp(- i L (\sqrt{\Delta} - \lambda)): \ker_0 (\sqrt{\Delta} - \lambda) \to
\ker_L (\sqrt{\Delta} - \lambda). $$
\end{prop}

\noindent{\bf Sketch of Proof} (See \cite{SjZ2} for more details)
It is obvious that $\exp(- i L (\sqrt{\Delta} - \lambda))$ takes
microlocal solutions into themselves since it commutes with
$\sqrt{\Delta}$. The rest of the proof is to verify
 that the complete symbol of the quasimode $\exp(- i L
(\sqrt{\Delta} - \lambda)) u_{\lambda}(s, y)$ at $(0, 0)$ is the
same as that of $u_{\lambda}(s, y)$ at $s = L$. QED

\subsubsection{Flux norm} We now describe  the {\it
flux norm} $\langle \cdot, \cdot \rangle$ on $\ker_{m}
(\sqrt{\Delta} - \lambda).$ It is motivated by the properties of
the probability current density ${\bf j}$ of a solution $\psi(x,
t)$ of the Schr\"odinger equation in quantum mechanics (see e.g.
\cite{LL}, \S 19). It is defined by
$${\bf j} = \frac{\hbar}{i}  (\psi \nabla \overline{\psi} -
\overline{\psi} \nabla \psi). $$ The integral $\int_S {\bf j}
\cdot dA$ over a surface measures the probability that the
particle described by $\psi$ will cross the surface $S$ in a unit
of time. As is pointed out in \cite{LL}, if $\psi = A
e^{\frac{i}{\hbar} S}$ then ${\bf j} = |A|^2 \nabla S. $

 The flux
norm on microlocal solutions is an invariantly defined version of
this. We consider the microlocal solution space $\ker_{m_0}$ for
$P = \hbar^2 \Delta - 1$ along an initial arc of  a closed
geodesic $\gamma$. We let $Y$ denote a transversal to $\gamma$ in
$M$.  If $\nu_Y$ is the unit normal to $Y$, then
$$||u||_{QF}^2 = \int_Y {\bf j} \cdot \nu_Y dvol_Y,\;\;\;u \in
\ker_{m_0}(P). $$ More precisely and generally,
$$||u||_{QF}^2 = \frac{i}{h} \;\langle  [ P, \chi] u, u \rangle,$$
where we use the semi-classical notation   $P = -h^2 \Delta - 1$
with $h = \lambda^{-1}$, and where $\chi$ is a microlocal cutoff
defined in the universal cover of a small tube around $\gamma$,
whose complete symbol  equals  $0$ before $m_0$ and $1$ after
$m_0$ in a somewhat smaller tube. We use the semi-classical
notation to conform to the notation of \cite{ISZ, SjZ} and also
because their results apply to much more general semi-classical
operators $P(h)$. For further details we refer to \cite{ISZ, SjZ}.

\begin{prop} \cite{ISZ} The monodromy operator ${\cal M}$ is a unitary
operator on $\ker_{m_0}(\sqrt{\Delta} - \lambda)$ with respect to
the flux norm $||\cdot||_{QF}. $ \end{prop}

\subsubsection{Grushin problem}

One of the basic steps of   \cite{SjZ} is the proof of a trace
formula which relates the trace of the wave group at $\gamma$ to
the traces involving the monodromy operator.  Formally, the  trace
formula says roughly that
\begin{equation} \label{GRUSHINTRACE} tr \rho (P/h)\chi(P) A =
\sum_{\pm} tr \int_{\R} \rho (z/h) \frac{d}{dz} \log det (I +
{\mathcal M}(z \pm i 0)) \chi(z) dz,
\end{equation}
where $\hat{\rho} \in C_0^{\infty}(L_{\gamma} - \epsilon,
L_{\gamma} + \epsilon)$ (compare to (\ref{RRHO})), and $\chi$ is a
cutoff to the sphere bundle.  The precise statement is given in
Theorem 1 of \cite{SjZ}:
\begin{equation} tr \rho (P/h)\chi(P) A = \frac{1}{2 \pi} \sum_{-N - 1}^{N + 1}
tr \int_{\R} \rho (z/h) {\cal M}(z, h)^k \frac{d}{dz} {\cal M} (z,
h) \chi(z) dz + O(h^{\infty}). \end{equation} This formula is
analogous to one in the case of domains with boundary, where the
boundary replaces the transversal,  see Proposition \ref{MAINCOR}.

To prove the formula, one sets up a microlocal  Grushin problem
near the closed trajectory $\gamma$. One begins by forming a
microlocally invertible system
$${\mathcal P}(z) = \left( \begin{array}{ll} \frac{i}{h} P(z) & R_-(z) \\ & \\
R_+(z) & 0 \end{array} \right): {\mathcal D}'(M) \times {\mathcal
D}'(\R^{n-1}) \to  {\mathcal D}'(M) \times {\mathcal
D}'(\R^{n-1}),  $$ where $R_{\pm}(z)$ are defined microlocally
near $\gamma$. The microlocal inverse of ${\mathcal P}(z)$ is
given by:
\begin{equation} \label{GRUSHIN} {\mathcal E}(z) = \left( \begin{array}{ll} E(z) & E_+(z) \\ & \\
E_-(z) & E_{-+}(z) \end{array} \right) \end{equation} We thus
obtain the key formula
\begin{equation} \label{GRUPOT} (\sqrt{\Delta} - \lambda)^{-1} \sim E(\lambda) - E_+(\lambda) ( I - {\mathcal M})^{-1}
E_-(\lambda),
\end{equation}  which reduces the eigenvalue problem to the
microlocal invertibility of $( I - {\mathcal M})$. Using this
formula one proves (\ref{GRUSHINTRACE}).  This  equation is
analogous to one  which occurs in the Fredholm-Neumann reduction
of the Dirichlet problem to the boundary of this type (see
\ref{POT}).

The authors of \cite{ISZ, SjZ, ISj} then put ${\mathcal
M}(\lambda)$ into quantum Birkhoff normal form. In view of
Proposition \ref{MW}, the normal form of the monodromy operator is
just what we called the semi-classical normal form of the wave
group. The operation $|_o$ thus has been interpreted as
restriction to the microlocal solution space, or as  part of a
Grushin reduction to the transversal. In \cite{Z4, Z5}, we proved
that the semi-classical normal form was a spectral invariant by
first turning into the homogeneous normal form, by proving that
the homogeneous normal form was a spectral invariant and then by
relating the semi-classical and homogeneous normal forms. The
formula (\ref{GRUSHINTRACE}) eliminates this latter step and
directly shows that the semi-classical normal form is a spectral
invariant by using the Grushin reduction.

At  the present time, the only application of the monodromy method
to inverse spectral results for plane domains is for the case of
bouncing ball orbits of analytic domains with  two symmetries as
in \cite{Z2}. This does not involve any new calculations since the
principal symbol of the normal form of the monodromy operator is
the classical Birkhoff normal form of ${\mathcal P}_{\gamma}$ and,
as mentioned above, this normal form determines $\Omega$ when it
is analytic with two symmetries (as proved in \cite{CdV}).

\section{\label{BB} Calculation of Wave Invariants II: Balian-Bloch approach}

In this section, we describe a method for calculating wave
invariants of bounded domains which does not use normal forms or
parametrices, but rather is based on an exact formula for the
Dirichlet (or Neumann) resolvent in terms of the `free resolvent'.
Our goal is to explain how the calculations of Theorem
\ref{BGAMMAJ} and Corollary \ref{BGAMMAJSYM} are done.

 For simplicity, we will confine ourselves
to calculating wave invariants at a bouncing ball orbit of a
bounded simply connected domain $\Omega \subset \R^2$. The method
extends with no difficulty to higher dimensions, more general
metrics and domains (in particular, non-convex domains) and more
general periodic orbits, but the special case is already
sufficiently rich and difficult.

The method is based on the use of a classical  formula due to C.
Neumann and I. Fredholm  for the resolvent $R_B^{\Omega}(z)$ for
the boundary problem in $\Omega$ with boundary conditions $B$ in
terms of layer potentials and the boundary integral operators they
induce.
 This approach was first used by the physicists
Balian-Bloch \cite{BB1, BB2} in their well-known  work on the
Poisson relation on three dimensional Euclidean domains, and we
will refer to it as the Balian-Bloch approach.  The crucial
advantage of this approach over the use of normal forms or
monodromy operators is its computability. For the first time, one
can calculate wave invariants of all orders explcitly and find out
what information they contain.

\subsection{Wave invariants as semi-classical  resolvent trace invariants}

First, we explain  how wave trace expansions are equivalent to
esolvent expansions. We fix
 a non-degenerate bouncing ball orbit  $\gamma$  of length
 $L_{\gamma}$, and
  let  $\hat{\rho} \in C_0^{\infty}(L_{\gamma} - \epsilon, L_{\gamma} + \epsilon)$ be a cutoff,  equal to one on an interval  $(L_{\gamma} - \epsilon/2, L_{\gamma} + \epsilon/2)$ which contains
no other lengths in Lsp$(\Omega)$ occur in its support. We then
define the  regularized resolvent by
\begin{equation}\label{RRHO} R_{\rho}^{\Omega}(k + i \tau):= \int_{\R} \rho
 (k - \mu) (\mu + i \tau)  R_B^{\Omega}(\mu + i \tau) d \mu.  \end{equation}
From the resolvent identity (e.g.)
$$ R_{B}^{\Omega} (\mu + i \tau) = \frac{1}{\mu + i \tau} \int_0^{\infty} e^{i (\mu + i \tau) t} E_{B}^{\Omega^c}(t) dt,\;\;\;\; $$
it follows that
\begin{equation} \label{RW} R_{\rho B}^{\Omega}(k + i \tau)
= \int_0^{\infty} \hat{\rho}(t) e^{i (k + i \tau) t}
 E_{B}^{\Omega}(t) dt. \end{equation}
When $\gamma, \gamma^{-1}$ are the unique closed orbits of length
$ L_{\gamma}$,  it follows from the Poisson relation for manifolds
with boundary (see \S 3 and \cite{GM, PS}) that the trace $Tr
1_{\Omega} R_{\rho}((k + i \tau))$ of the regularized resolvent on
$L^2(\Omega)$  admits a complete asymptotic expansion  of the
form:
\begin{equation} \label{PR} Tr 1_{\Omega} R_{\rho}(k + i \tau) \sim  e^{ (i k - \tau) L_{\gamma}}  \sum_{j = 1}^{\infty} (B_{\gamma, j}
+  B_{\gamma^{-1}, j}) k^{-j},\;\;\; k \to \infty \end{equation}
with coefficients $B_{\gamma, j},  B_{\gamma^{-1}, j}$ determined
by the jet of $\Omega$ at the reflection points of $\gamma$. The
coefficients $B_{\gamma, j},  B_{\gamma^{-1}, j}$ are thus
essentially the same as the wave trace
 coefficients at the singularity  $t = L_{\gamma}$.

 \subsection{Reduction to the boundary}

To calculate $Tr 1_{\Omega} R_{\rho}(k + i \tau)$ asymptotically,
we use  the Fredholm-Neumann reduction of the Dirichlet or Neumann
problems in a bounded domain to the boundary. The key formula is
similar to (\ref{GRUPOT}):
\begin{equation} \label{POT}  R_{\Omega}(k + i \tau) =
 R_0(k + i \tau) - {\mathcal D} \ell(k + i \tau) (I +  N(k + i \tau))^{-1}
 \gamma
 {\mathcal S} \ell^{tr}(k + i \tau), \end{equation}
where $R_0(k + i \tau)$ is the free resolvent  $- (\Delta_0 + (k +
i \tau)^2)^{-1}$  on $\R^2$. Here, $\gamma: H^{s}(\Omega) \to
H^{s-1/2}(\partial \Omega)$ is the restriction to the boundary,
and
 ${\mathcal D} \ell (k + i \tau)$ (resp. ${\mathcal S} \ell (k + i \tau)$) is the double
 (resp. single) layer potential is the operator from $H^s(\partial \Omega) \to H^{s+ 1/2}_{loc}
(\Omega)$ defined by
\begin{equation} \label{layers}\left\{ \begin{array}{l} {\mathcal S} \ell(k + i \tau)  f(x) = \int_{\partial \Omega} G_0(k + i \tau, x, q) f(q)
d s (q), \\ \\   {\mathcal D} \ell(k + i \tau)f(x) =
\int_{\partial \Omega}\frac{\partial}{\partial \nu_y} G_0(k + i
\tau, x, q) f(q) ds(q),
 \end{array} \right. \end{equation}
where $ds(q)$ is the arc-length measure on $\partial \Omega$,
where $\nu$ is the interior unit normal to $\Omega$, and where
$\partial_{\nu} = \nu \cdot \nabla$.  Also, $\SL^{tr}$ is its
transpose (from the interior to the boundary), and $G_0(z, x, y)$
is the free Green's function, i.e. the  kernel of $R_0(z)$.
Further,
\begin{equation} \label{blayers} N(k + i \tau)f(q) =  2 \int_{\partial
\Omega}\frac{\partial}{\partial \nu_y} G_0(k + i \tau, q, q')
f(q') ds(q')
\end{equation} is the boundary integral operator induced by ${\mathcal D} \ell$. It is classical that this
operator satisfies
\begin{equation}  \label{ISO} \begin{array}{ll} (i) &
\; N(k + i \tau) \in \Psi^{-1}(\partial \Omega), \\
&  \\ (ii) &
 (I + N(k + i \tau)): H^s(\partial \Omega) \to H^s(\partial \Omega)\;\; \mbox{is\; an\;
isomorphism}. \end{array} \end{equation}

The reader might now compare  the formulae (\ref{POT}) and
(\ref{GRUPOT}) to see how $N(\lambda)$ is analogous to ${\mathcal
M}(\lambda)$. Moreover, it is used in an analogous way to reduce
the calculation of  the (distribution) trace of $R_{\Omega}(k + i
\tau)$, formula (\ref{POT})  to the boundary. In doing so, we
found it \cite{Z10} more convenient to combine the interior and
exterior problems as follows:  We write $L^2(\R^2) = L^2(\Omega)
\oplus L^2(\Omega^c)$
 and  let $R_{N }^{\Omega},$ resp. $ (k + i \tau), R_{ D}^{\Omega^c}(k + i \tau)$ denote the Neumann resolvent
 on the exterior domain, resp. the Dirichlet resolvent on the interior domain.  We then regard
 $R_{ D}^{\Omega^c}(k + i \tau)
 \oplus  R_{N }^{\Omega} (k + i \tau) $ as an operator on this
 space. The reason for combining the interior/exterior problems is
 that we can cycle around the layer potentials in (\ref{POT}) when
 taking the trace and simplify the formula to:

\begin{equation} \label{IO}
\begin{array}{l}
 Tr_{\R^2} [R_{ D}^{\Omega^c}(k + i \tau)
 \oplus  R_{N }^{\Omega} (k + i \tau)  - R_{0 }(k + i \tau)] \;= \;  \frac{d}{d k } \log \det
 \bigg(I + N(k + i \tau)\bigg)
,\end{array}
\end{equation}
where the determinant is the usual Fredholm determinant. We refer
to \cite{Z10} for a proof of this apparently well-known formula.
This gives:

\begin{prop} \label{MAINCOR}  Suppose that $L_{\gamma}$ is the only length
in the support of $\hat{\rho}$. Then,
$$\int_{\R} \rho(k -
\lambda) \frac{d}{d \lambda} \log \det (I + N(\lambda + i \tau)) d
\lambda  \sim \sum_{j = 0}^{\infty} B_{\gamma; j} k^{-j}, $$ where
$B_{\gamma; j}$ are the wave invariants of $\gamma$.
\end{prop}

\subsection{Semi-classical analysis of $N(k + i \tau)$}

The next step is to analyse the operator $N(k + i \tau)$ and the
geometric series expansion of $(I + N(k+ i \tau)^{-1}$.  The
operator $N(k + i \tau)$
 has the singularity of a homogeneous
pseudodifferential operator of order $-1$ on the diagonal (in
fact, it is of order $-2$ in dimension $2$) and that  is the way
it is normally described in potential theory.
 However, away from the
diagonal, it has a WKB approximation which exhibits it as a
semi-classical Fourier integral operator with phase $d_{\partial
\Omega}(q,q') = |q - q'|$ on $\partial \Omega \times
\partial \Omega$, the boundary distance function of $\Omega$.
Indeed, the free Green's function in dimension two  is given by:
$$\begin{array}{l} G_0(k + i \tau, x, y) =
  H^{(1)}_0((k + i \tau) |x - y|)  = \int_{\R^2} e^{i \langle x - y, \xi \rangle} (|\xi|^2 - (k + i \tau)^2)^{-1} d \xi. \end{array}$$
Here, $H^{(1)}_0(z)$ is the Hankel function of index $0$. It has
the near diagonal and off-diagonal asymptotics,
\begin{equation}\label{LK} (i) \;\; H^{(1)}_{0} (z) \sim \left\{ \begin{array}{l}
 - \frac{1}{2\pi} \ln |z| \;\rm{as}\; |z| \to 0,\\ \\(ii) \;\;
 e^{i (z |x - y| -  \pi/4)} \frac{1}{|z|^{1/2}}\;\rm{as}\;
 |z| \to \infty.  \end{array} \right. \end{equation}

By the explicit formula we have:
$$ \begin{array}{lll}
\frac{1}{2} N(k + i \tau, q(\phi), q(\phi')) & = &
\partial_{\nu_y} G_0(\mu,  q(\phi), q(\phi')) \\ [6pt] & = &
-\ (k + i \tau ) H^{(1)}_1 (k + i \tau|q(\phi) - q(\phi')|)\\
[6pt] & &\times\  \cos \angle(q(\phi) - q(\phi'), \nu_{q(\phi)}).
\end{array}$$
Combining with  (\ref{LK}(ii)), we see that when $ |q(\phi) -
q(\phi')|  \geq |k|^{1 - \epsilon}$ for some $ \epsilon < 1$, , $
N(k + i \tau, q(\phi), q(\phi'))$ is a semi-classical Fourier
integral kernel whose phase is the boundary distance function. For
a convex domain, the boundary distance function generates the
billiard map of $\partial \Omega$ and hence we view $N(k + i
\tau)$ as a global  quantization of the billiard map. For
non-convex domains, the boundary distance function additionally
generates `ghost orbits' which in part exit the domain, but these
only present a mild complication. We refer to \cite{HZel} for
discussion of these orbits and to \cite{Z2, HZel} for further
discussion of $N(k)$.

 We now explain how to use  Proposition (\ref{MAINCOR}) to  calculate
the wave trace coefficients $B_{\gamma, j}$ at a closed geodesic.
When calculating these coefficients we first   compose with a
special kind of semiclassical  cutoff operator $\chi (x, k^{-1}
D_x)$ to a neighborhood of $\gamma$  on both sides of (\ref{POT}).
Since
 we are not dealing with conventional Fourier integral
operators,  it must be proved that composition with such a cutoff
does in fact microlocalize to $\gamma$.   We will suppress this
issue until the end.

At least formally, we
 expand  $(I +  N(k + i \tau))^{-1}$ in a finite geometric
series plus remainder:
\begin{equation} \label{GS} (I \!+\! N(\lambda \!+\! i \tau))^{-1} = \sum_{M = 0}^{M_0} (-1)^M \; N(\lambda)^M
+
 (-1)^{M_0 + 1} \; N(\lambda)^{M_0 + 1}  (I
\!+\! N(\lambda \!+\! i \tau))^{-1}. \end{equation} To calculate a
given wave invariant, we need to show that, for each  order
$k^{-J}$ in the trace expansion of Corollary (\ref{MAINCOR}),
there exists $M_0(J)$ such that
\begin{equation} \label{ASY} \begin{array}{ll} (i) &
\sum_{M = 0}^{M_0} (-1)^M  Tr  \int_{\R} \rho(k - \lambda)\;
N(\lambda)^M N'(k + i \tau) d \lambda \\ & \\ &  = \sum_{j =
0}^{J} B_{\gamma; j} k^{-j} + O(k^{-J - 1}) , \\&  \\(ii) & Tr
\int_{\R} \rho(k - \lambda)  N(\lambda)^{M_0 + 1}  (I \!+\!
N(\lambda \!+\! i \tau))^{-1} N'(k + i \tau) d \lambda = O(k^{-J -
1}).
\end{array} \end{equation}

We outline the method for obtaining the wave trace asymtotics at a
closed geodesic $\gamma$  from the first term. Since $N$ has
singularities on the diagonal, one cannot just apply the
stationary phase method to the trace. Rather one has to regularize
the operator.  We do this by  separating out the tangential and
transversal parts of $N$ by introducing  a cutoff of the form
$\chi(k^{1 - \delta} |q - q'| )$ to the diagonal, where $\delta >
1/2$ and where $\chi \in C_0^{\infty}(\R)$ is a cutoff to a
neighborhood of $0$. We then put
\begin{equation} \label{N01}  N(k + i \tau) = N_0(k + i \tau) + N_1(k + i
\tau), \;\; \mbox{with}  \end{equation}
\begin{equation} \label{N01DEF} \left\{\begin{array}{l}  N_0(k + i
\tau, q, q') = \chi(k^{1 - \delta} |q - q'| ) \;  N(k + i \tau, q,
q'), \\ \\ N_1(k + i \tau, q, q') = (1 - \chi(k^{1 - \delta} |q -
q'| ))\; N(k + i \tau, q, q'). \end{array} \right.
\end{equation}
The term $N_1$ is a semiclassical Fourier integral kernel
quantizing the billiard map, while $N_0$ behaves like an Airy
operator close to the diagonal with the singularity of a
homogeneous pseudodifferential operator on the diagonal.

Now consider the powers $N(k + i \tau)^M$ in the first term (i)
of (\ref{ASY}).  We write
\begin{equation} \label{BINOMIAL} (N_0 + N_1)^M = \sum_{\sigma: \{ 1, \dots, M\} \to \{0, 1\}}
N_{\sigma(1)} \circ N_{\sigma(2)} \circ \cdots \circ
N_{\sigma(M)}. \end{equation} We  regularize $N^M$ by eliminating
the factors of $N_0$ from each of these terms. This is obviously
not possible for the term $N_0^M$ but  it is possible for the
other terms. By explicitly writing out the composition in terms of
Hankel functions and using the basic identities for these special
functions, we prove that $N_0 \circ N_1 \circ \chi_0 (k + i \tau,
\phi_{1}, \phi_{2}) )$ is a semiclassical Fourier integral
operator on $\partial \Omega$ of order $-1$ associated to the
billiard map. Thus,  composition with $N_0$ lowers the order. The
remaining terms $N_0^M$, when composed with a cut-off to $\gamma$,
do not contribute asymptotically to the trace.

The successive removal of  the factors of $N_0$ thus gives a
semi-classical quantization of the billiard map near $\gamma$. We
then calculate the traces of each term by the stationary phase
method and obtain the result stated in (\ref{BGAMMAJ}). The terms
displayed there with the maximum number of derivatives of the
defining function of $\partial \Omega$ come only from the  $N_1^M$
terms. Thus, although
 the removal of the $N_0$
factors does change  the amplitude of the $N_1$ factors and does
contribute to the wave trace, it turns out to be negligeable in
the inverse spectral problem. This is one of the principal virtues
of the Balian-Bloch approach.

We will explain how the explicit formulae for the wave invariants
are derived in the next section.

\subsubsection{Remainder estimate} We now address some sketchy remarks to the remainder
estimate, in the hopes of convincing the reader that it is
plausible to expect the remainder trace to be small.  The ability
to insert a microlocal cutoff to $\gamma$ is crucial here. The
main obstacle to the remainder is that the norm of  $N(k + i
\tau)$ fails to decrease with increasing $\tau$ due to the Airy
part associated to creeping rays. This again is a difference to
the monodromy operator.

 To obtain a small
remainder we set the spectral parameter in $N$ equal to $k + i
\tau \log k$. The presence of the $\log k$ in the imaginary part
changes the wave trace expansions by $k^{- C r L_{\gamma}}$, but
this does not hurt the expansions such the remainder estimates
will be of lower order.  We then estimate the remainder
$$Tr \int_{\R} \rho(k - \lambda)  N(\lambda)^{M_0 + 1}  (I \!+\!
N(\lambda \!+\! i \tau))^{-1} \; \chi_{\gamma} \; N'(k + i \tau) d
\lambda$$ by applying the Schwarz inequality for the
Hilbert-Schmidt inner product, and using the relation $(I \!+\!
N(\lambda \!+\! i \tau))^{-1}$ to the Poisson kernel to estimate
this factor. These estimates leave a trace of $N^M N^{*M}$ for
fixed $M$,  microlocalized to $\gamma$. We regularize these traces
as above and obtain osillatory integrals whose phases  have
critical points corresponding to $M$-link closed circuits which
being  at some point $q$, bounce along the boundary until $q'$ and
then return to $q$ by traversing the links in reverse order.
However, the cutoffs to $\gamma$    force the links in critical
paths to point in the direction of $\gamma$ and hence to be of
length roughly $M L_{\gamma}$. The imaginary part $i \tau \log k$
of the semiclassical parameter then contributes a damping factor
of $e^{- \tau M L_{\gamma} \log k}$ for each link. The links
correspond to the $N_1$ factors. Thus, for each string,  we have
one $k^{-1}$ for each $N_0$ factor and one  $e^{- \tau M
L_{\gamma} \log k}$ for each $N_1$ factor.  For sufficiently large
$\tau$ these combine to give a factor of $k^{-R}$ for any
prescribed $R$.

\subsubsection{$N(\lambda)$ versus ${\mathcal M}(\lambda)$}

We digress momentarily to compare $N(\lambda)$ and the monodromy
operator ${\mathcal M}(\lambda)$. Both operators arise as a
reduction of the wave group to the boundary (or a transversal) and
are quantizations of the billiard map, $N(\lambda)$ globally and
${\mathcal M}(\lambda)$ microlocally at a closed orbit $\gamma$.
(An earlier reduction also occurs in \cite{MM, P, P3} } but in a
rather different way. ) In applications to boundary problems, they
are quite similar, as evidence by the comparisons of
(\ref{GRUSHINTRACE}) with Proposition \ref{MAINCOR} and of
(\ref{GRUPOT}) with (\ref{POT}).

But there do exist  rather important differences. The operator
$N(\lambda)$ is not a standard Fourier integral operator on
$\partial \Omega$, while ${\mathcal M}(\lambda)$ is one.  The
diagonal singularities of $N(\lambda)$ require a complicated
regularization procedure. Moreover  it is a global invariant of
$\partial \Omega$, not a microlocal one at $\gamma$; On the other
hand, $N(\lambda)$ is just the  restriction of a canonically
defined free Green's kernel to $\partial \Omega \times
\partial \Omega$, and thus is an elementary and computable object.
By comparison, ${\mathcal M}(\lambda)$ must be constructed by some
kind of parametrix method. It was precisely the complexity of
micrlocal  parametrices for the wave group, even at periodic
reflecting rays, which motivated our turning to  the Balian-Bloch
approach.

\subsection{Stationary phase expansion}

After regularizing the traces, we end up with oscillatory
integrals in the standard sense and obtain expansions by applying
stationary phase. The amplitudes and phases are canonical, since
we began with canonical amplitudes and phases and since the
regularization procedure is essentially the same for all domains.
So most of the complexity of the expansion is due to the
stationary phase method.

The key point in inverse spectral theory is to identify data in
the stationary phase term of order $k^{-j}$ which represents `new
data' not contained in the previous terms. Terms of the
coefficient of  $k^{-j}$ which contain the maximum number of
derivatives of the phase are the most important ones. Thus we face
the combinatorial problem of locating such terms in the stationary
phase expansion. The
 Feynman diagram method of assigning labelled
graphs to each term in the expansion proves to be very effective
for this purpose.

Consider a general  oscillatory integral $Z_k = \int_{\R^n} a(x)
e^{ik S(x)} dx$ where $a \in C_0^{\infty}(\R^n)$ and where $S$ has
a unique critical point in supp$a$  at $0$. Let us write $H$ for
the Hessian of $S$ at $0$. The stationary phase expansion takes
the form:
$$\begin{array}{l}
 Z_k = \bigg(\frac{2\pi}{k}\bigg)^{n/2} \frac{e^{i \pi sgn
(H)/4}}{\sqrt{|det H|}} e^{i k S(0)} Z_k^{h \ell}, \\ [6pt]
\hspace{70pt}\mbox{where}\;\;  Z_k^{h \ell} =  \sum_{j =
0}^{\infty} k^{-j} \bigg\{\sum_{(\Gamma, \ell): \chi_{\Gamma'} =
j} \frac{I_{\ell} (\Gamma)}{S(\Gamma)}\bigg\}.
\end{array} $$

Here, the sum runs over the  set ${\mathcal G}_{V, I}$  of
labelled graphs $(\Gamma, \ell)$ with $V$ closed  vertices of
valency $\geq 3$ (each corresponding to the phase), with one open
vertex (corresponding to the amplitude), and with $I$ edges.
Further, the  graph $\Gamma'$ is defined to be $\Gamma$ minus the
open vertex, and $\chi_{\Gamma'} = V - I$ equals its Euler
characteristic. We note that there are only finitely many graphs
for each $\chi$ because the valency condition forces $I \geq 3/2
V.$ Thus, $V \leq 2 j, I \leq 3 j.$

The function $\ell$ `labels' each end of each edge of $\Gamma$
with an index $j \in \{1, \dots, n\}.$ Also, $S(\Gamma)$ denotes
the order of the automorphism group of $\Gamma$, and  $I_{\ell}
(\Gamma)$ denotes the `Feynman amplitude' associated to $(\Gamma,
\ell)$. By definition, $I_{\ell}(\Gamma)$  is obtained by the
following rule: To each edge with end labels $j,k$ one assigns a
factor of $\frac{-1}{ik} h^{jk}$ where $H^{-1} = (h^{jk}).$ To
each closed vertex one assigns a factor of $i k
\frac{\partial^{\nu} S (0)}{\partial x^{i_1} \cdots \partial
x^{i_{\nu}}}$ where $\nu$ is the valency of the vertex and $i_1
\dots, i_{\nu}$ at the index labels of the edge ends incident on
the vertex. To the open vertex, one assigns  the factor
$\frac{\partial^{\nu} a(0)}{\partial x^{i_1} \dots \partial
x^{i_{\nu}}}$, where $\nu$ is its valence.   Then
$I_{\ell}(\Gamma)$ is the product of all these factors.  To the
empty graph one assigns the amplitude $1$.  In summing over
$(\Gamma, \ell)$ with a fixed graph $\Gamma$, one sums the product
of all the factors as the indices run over $\{1, \dots, n\}$.

\subsubsection{ The data
\boldmath{$f_{\pm}^{2j}(0)$}}

An analysis of the diagrams and amplitudes shows that the $j$th
even Taylor coefficients  $f_{\pm}^{(2j)}(0)$ of the boundary
defining functions  appear first in the $k^{-j + 1}$ term.

When the domain has one symmetry axis, which we visualize as an
 up/down symmetry, the terms with this data have the form  $$ 2 r L (h^{11})^j f^{(2j)}(0) +
\cdots,$$ where  $\cdots$ refers to terms with $\leq 2j - 1$
derivatives.

When the domain has two symmetries, a left/right symmetry in
addition to an up/down symmetry, the odd Taylor coefficients
vanish and we see immediately that the even Taylor coefficients
can be determined inductively from the wave trace invariants. This
gives a new proof that analytic domains with two symmetries can be
determined from the wave trace invariants at a bouncing ball orbit
which is one of the symmetry axes.

This does not quite prove that such domains are spectrally
determined among other analytic domains with two symmetries, since
the length of the bouncing ball orbit must be known in order to
obtain the wave invariants. This length is a spectral invariant if
the domains are additionally convex \cite{Gh}. For non-convex
domains we needed to add as an assumption that the bouncing ball
symmetry axis had a fixed length $L$ in \cite{Z1}.

\subsubsection{ The data
\boldmath{$f_{\pm}^{(2j - 1)}(0)$}}

When we do not assume a left/right symmetry, the odd Taylor
coefficients are non-zero in general, and the problem arises
whether there is sufficient information in the wave invariants to
determine all of the even and odd Taylor coefficients of the
boundary defining function (or curvature function) of a domain
with one symmetry.  We assume the axis of symmetry is a bouncing
ball orbit whose orientation is reversed by the symmetry.

An analysis of the diagrams and amplitudes for the odd
 Taylor coefficients $f_{\pm}^{(2j -
1)}(0)$ show that they appear first in the term of order $k^{-j +
1}. $  It turns out that five diagrams contain this data, but the
amplitudes of three  automatically vanish.  They two amplitudes
have the following  forms:

\begin{itemize}

\item (i)
$
(h^{pp}_{\pm})^{j - 1} h^{qq}_{\pm} h^{pq}_{\pm} f^{(2j - 1)}(0)
f^{(3)}(0).  $

\item (ii)  $(h^{pp}_{\pm})^{j - 2} (h^{pq}_{\pm})^3
f^{(2j - 1)}(0) f^{(3)}(0). $

\end{itemize}

To decouple them, we need to analyze the behavior of power sums of
columns
 the Hessian matrix elements. In \cite{Z5} we proved that cubic column sums
 are
 linearly independent from linear ones as  $r \to \infty$.

\subsection{\label{ANALYTIC} Positive results for analytic domains
and metrics }

We now review the proof in \cite{Z5} that  the Taylor coefficients
$f^{2j-1}(0), f^{2j}(0)$ can be determined from the wave
invariants $B_{\gamma^r,j}$ as $r$ varies over $r = 1, 2, 3,
\dots.$

It suffices to  separately determine the two terms
\begin{equation} \label{TWOTERMS}\begin{array}{l}
 2 (h^{11}_{2r})^2 \bigg\{f^{(2j)}(0) + \frac{1}{2 - 2
\cos \alpha/2}\ f^{(3)}(0)
f^{(2j - 1)}(0) \bigg\}, \\[9pt]
 \mbox{and} \quad \bigg\{  \sum_{q = 1}^{2r} (h^{1 q}_{2r})^3
\bigg\} f^{(3)}(0) f^{(2j - 1)}(0).\end{array}
\end{equation}
It is easy to see that the terms decouple   as $r$ varies if and
only if the cubic sums
 $  \sum_{q = 1}^{2r} (h^{1 q}_{2r})^3$ are non-constant in $r = 1, 2, 3,
 \dots$.
By the explicit calculation in \cite{Z2}, we  have:
$$
\begin{array}{clll} \sum_{q =
1}^{2 r } ( h^{pq})^3 =\\[9pt]
 2r  \!\!\!\sum_{k_1, k_2 = 0}^{2r} \frac{1 }{(\cosh \alpha/2
\!+\! \cos \frac{k_1 \pi}{r}) (\cosh \alpha/2 \!+ \!\cos \frac{k_2
\pi}{r})(\cos \alpha/2 \!+\! \cosh \frac{(k_1 \!+\!k_2) \pi}{r})}
. \end{array}$$ It is obvious that the sum is strictly increasing
as $r$ varies over even integers.

This independence of the linear and cubic sums permits us to use
an inductive argument when the domain has the given symmetry.
From the $j = 0$ term we determine $f''(0).$ Indeed,
 $(1 - L f^{(2)}(0) = \cos(h) \alpha/2$ and $\alpha$ is a wave trace invariant.
  From the $j = 2$ term we recover $f^{3} (0), f^{4}(0).$  The induction hypothesis is then that the Taylor polynomial
  of $f$ of degree $2j - 2$ has been recovered by the $j - 1$st stage. By the decoupling argument we can determine
  $f^{2j}(0), f^{2j - 1}(0)$. (Strictly speaking, there is the
  minor annoyance of the factor of $f^{(3)}(0)$, which is resolved
  in
 \cite{Z5}).

\subsubsection{Final comments on the Balian-Bloch invariants}

 It would of course be desirable to remove all the symmetry
 assumptions, if possible. Thus one would need to recover the
 Taylor coefficients $f_{\pm}^{(2j)}(0), f_{\pm}^{(2j - 1)}(0)$
 when $f_+ \not= f_-$ from the wave invariants. The problem is
 that the induction does not work because the term
$R_{2r} (j^{2j - 2} f(0))$, which we do not need to know in the
symmetric case, might have a different dependence on the Taylor
coefficients of $f_{\pm}$ than the `principal term', i.e. the one
with the highest derivatives. The latter is basically a sum of the
coefficients of $f_{\pm}$ while $R_{2r} (j^{2j - 2} f(0))$ could
involve any symmetric polynomial in these coefficients.  At the
present time, we do not even know whether the set
$\{f_{+}^{(2)}(0), f_-^{(2)}(0)\}$ of second derivatives can be
determined from the second  wave coefficient.

 Some hope has been provided by computer calculations of
 C. Hillar \cite{Hil}. The  results suggest that Hessian column power sums with different powers are
 linearly independent as $r \to \infty$. This would give  quite
 a large supply of Taylor coefficients.

\section{Surfaces and domains with integrable dynamics}

As mentioned above, Birkhoff normal forms at a periodic orbit (on
both the classical and quantum level) are approximations to the
Hamiltonian by an integrable one. When the Hamiltonian is
completely integrable, its dynamics and Birkhoff normal forms are
very special  and the metrics might be rather spectrally rigid, at
least in low dimensions. This is the case with flat metrics
\cite{Ku3}.  In this section, we consider some model inverse
spectral problems on two-dimensional surfaces or domains with
integrable dynamics.

Let us recall that  the geodesic flow of an $n$-dimensional
Riemmanian manifold $(M, g)$ is {\it completely integrable} if it
commutes with a Hamiltonian action of $\R^n$ on $T^*M-0$ ($n =
\dim M$). That is, the metric Hamiltonian $|\xi|_g$ satisfies
$$\{|\xi|_g, p_j\} = 0 = \{p_i, p_j\},\;\; i,j = 1, \dots, n$$
where $p_j: T^*M-0 \to \R$ are homogeneous of degree one and are
independent in the sense that
$$dp_1 \wedge dp_2 \wedge \cdots \wedge dp_n \not= 0\;\; \mbox{on
a dense open subset}\; U \subset T^*M. $$

The orbits of an $\R^n$ action give a (usually singular) foliation
of $S_g^*M$ (the unit sphere bundle for the metric), called the
Liouville foliation,  by affine manifolds of the form $\R^n \cdot
(x, \xi) \equiv \R^n/\Gamma$ where $\Gamma$ is the isotropy group
at $(x, \xi)$. If it is a lattice of full rank, the orbit is a
torus $T^n$ of dimension $n$. If it has less than full rank, the
orbit type is $T^k \times \R^{n - k}$. The isotropy group might
have the form $\R^k \times \Z^{n - k}$ in which case the orbit
becomes a singular torus of dimension $n - k.$ The orbits which
contain periodic geodesics are sometimes called `periodic tori',
and they furnish the components ${\cal T}$ in the trace formula.

In the case of bounded plane domains, complete integrability of
the billiard map $\beta$  means that $B^* \partial \Omega$ is
foliated by invariant curves.  The well-known conjecture of
Birkhoff is that ellipses are the only example of compact smooth
Euclidean plane domains with integrable billiard map. There are
more examples if curved metrics are allowed, see Popov-Topalov
\cite{PT}.

\subsection{\label{INT}Trace formulae for integrable systems}

So far, we have mainly considered the trace of the wave group
around non-degenerate periodic  orbits. For integrable systems,
the  periodic orbits usually come in families filling out
invariant tori. We now consider the appropriate notions of
non-degeneracy in this context.

We will always assume that the closed geodesics come in clean
families in the following sense:
\begin{defn}  A metric $g$ on a compact manifold $M$ will be
said to have a {\em simple clean length spectrum} if the length
function $L_g$ on the loop space $Map (S^1, M)$ is a Bott-Morse
function which takes distinct values at distinct components of its
critical set, $Crit(L_g)$. \end{defn}

The term Bott-Morse means that each component of $Crit(L_g)$ is a
manifold, whose tangent space is the kernel of $d L_g$.
Equivalently, each component is a clean fixed point set for
$G^t_g$.  One needs the clean (Bott-Morse) condition to get a nice
wave trace expansion and one needs the simple length spectrum
condition to determine geometric information from the expansion.

Under this assumption,  the trace of its wave group has the form:
\begin{equation} \label{SCTR} Tr e^{i t \sqrt{\Delta_g}} = e_0(t)  +
\sum_{{\cal T}} e_{\cal T}(t) \end{equation}  where the singular
term  $e_0(t)= C_n \;Vol (M,g) \;(t + i0)^{-n} + \dots$ at $t = 0$
is the same as in the non-degenerate case, where $\{{\cal T}\}$
runs over the critical point components of $L_g$, and where
$$e_{\cal T} = c_{{\cal T}, d_T} (t - L_{{\cal T}} + i0)^{-d_T/2} +
c_{{\cal T}, d_T - 1} (t - L_{{\cal T}} + i0)^{(-d_T /2 + 1)} +
\dots, \;\;\; d_{{\cal T}} = \dim {\cal T}. $$
 Here, $d_T$ is the
dimension of the symplectic cone formed by the  family of closed
geodesics within $T^*_g M$ and $L_{{\cal T}}$ is the common length
of the closed geodesic in ${\cal T}$. For instance, in the
non-degenerate case, ${\cal T} = \R_+ \gamma$ is the symplectic
cone generated by $\gamma \subset S^*_g M$ and $d_{{\cal T}} = 2$.

\subsection{Spectral determination of simple surfaces of revolution}

The simplest example of what can be done with integrable systems
is given by `simple' analytic surfaces of revolution. We first
review the proof in \cite{Z2} that they are spectrally determined
among other simple analytic surfaces of revolution and then sketch
a result observed independently by the author and G. Forni and by
K.F. Siburg regarding their spectral determination among all
metrics on $S^2$.

The precise class of  metrics we consider are those  metrics $g$
on $S^2$ which belong to the
 class ${\cal R}^*$ of real analytic, rotationally
invariant metrics on $S^2$ with simple length spectrum in the
above sense and satisfying the following `simplicity' condition
$$\begin{array}{ll} \bullet & g = dr^2 + a(r)^2 d\theta^2 \\ & \\
\bullet & \exists ! r_0 : a'(r_0) = 0; \\ & \\
\bullet & \mbox{The Poincare map}\;\; {\cal P}_0\;\; \mbox{of}\;\;
r = r_0 \;\;\mbox{is elliptic of twist type} \end{array} $$

Convex analytic surfaces of revolution are examples, but of course
there are others. The unique isolated closed geodesic (at distance
$r_0$) is an elliptic orbit.  Peanut surfaces are obviously
non-simple, since they possess three isolated closed geodesics of
which one is hyperbolic. The standard round metric and other Zoll
metrics of revolution are  non-simple since the Poincare maps are
not twist maps.

The following results shows that one can solve the inverse
spectral problem in the class of analytic simple surfaces of
revolution.

\begin{theo}\label{ISPSR} (\cite{Z2}) Suppose that $g_1, g_2$ are two real analytic metrics on
$S^2$ such that $(S^2, g_i)$ are simple surfaces of revolution
with simple length spectra.  Then Sp($\Delta_{g_1}$) =
Sp($\Delta_{g_2}$) implies $g_1 = g_2.$ \end{theo}

The proof is based on quantum Birkhoff normal forms for the
Laplacian $\Delta$. But the special feature of simple surfaces of
revolution is that there exists a global Birkhoff normal form as
well as local ones around the critical closed orbit or the
invariant tori. This is because $\Delta$ is a toric integrable
Laplacian in the following sense: there exist commuting first
order pseudo-differential operators $\hat{I}_1, \hat{i}_2$ such
that:
\begin{itemize}

\item the joint spectrum is integral, i.e.
 $Sp({\cal I}) \subset \Z^2 \cap \Gamma + \{\mu\}$ where $\Gamma$ is the cone $I_2 \geq |I_1|$ in $\R^2$.

\item The square root of $\Delta$ is a first order polyhomogenous function  $\sqrt{\Delta} =
\hat{H}(\hat{I}_1, \hat{I_2})$ of the action operators.
\end{itemize}

By polyhomogeneous, we mean that $\hat{H}$ has an asymptotic
expansion in homogeneous functions of the form:
$$\hat{H} \sim H_1 + H_o + H_{-1} + \dots, \;\;\;\;\;H_j( r I) = r^j
H_j(I).$$ The principal symbols $I_j$ of the $\hat{I}_j$'s
generate a classical Hamiltonian torus action on $T^*S^2 - 0$.
Analysis of the normal form shows that $H_o = 0$.

It follows that
$$Sp(\sqrt{\Delta_g}) = \{ \hat{H} ( N + \mu): N \in \Z^2 \cap \Gamma_o\},$$
where the eigenvalues have expansions
$$\lambda_N \sim H_1(N + \mu) + H_{-1}(N + \mu) + \dots.$$

\begin{theo}(\cite{Z2}) Let $(S^2, g)$ be an analytic simple surface of revolution with simple
length spectrum.  Then the normal form $\hat{H}(\xi_1, \xi_2)$  is
a spectral invariant.\end{theo}

The normal form and the proof are very different from the
non-degenerate case in \cite{G, Z3, Z4}, although the philosophy
of the proof is similar.

To complete the proof of  Theorem \ref{ISPSR}, we need to show
that $\hat{H}$ determines a metric in ${\cal R}.$ As in the
bounded domain case outline above, the crucial point  is to
calculate the normal form invariants. It turns out to be
sufficient to calculate $H_1 = H$ and $H_{-1}$ in terms of the
metric (i.e. in terms of $a(r)$) and then to invert the
expressions to determine $a(r)$.

The method given in \cite{Z2} for calculating  $H$ and $H_{-1}$
was to  study the spectral asymptotics of $\sqrt{\Delta} =
\hat{H}(\hat{I}_1, \hat{I_2})$ along `rays of representations'  of
the quantum torus action, i.e. along multiples of a given lattice
point $(n_o, k_o).$    The  lattice points $(n_o, k_o + \half)$
parametrize tori $T_{n_o, k_o}$ satisfying so-called
Bohr-Sommerfeld quantization conditions, which imply that one can
construct associated   joint eigenfunctions $\phi_{n_o, k_o}$ of
$(\hat{I}_1, \hat{I}_2)$ by the WKB method. This is reminiscent of
the quasi-mode method mentioned above around non-degenerate closed
geodesics, but there we studied the conjugation to normal form
rather than the asymptotics of quasi-modes. Here, the existence of
a global torus action makes the quasi-modes easier to study, and
indeed the  $\{\phi_{n,k}\}$ are actually modes  (eigenfunctions)
of $\Delta$ with complete asymptotic expansions along rays. By
studying the eigenvalue problem as $|(n,k)|\rightarrow \infty$ we
 determine the $H_{-j}$'s.

Once it is known that the global Birkhoff normal for $\hat{H}$ is
a spectral invariant,  it follows that for each $n_0 \in \Z$, the
function  $\hat{H}(n_o, \hat{I}_2)$ is a known function of the
variable $I:=\hat{I}_2$. Its principal symbol $H_{n_o}(I):= H_1
(n_o, I)$ is then a known function and its inverse function
$$I_{n_o} (E) = \int_{r_{-}(E)}^{r_+ (E)} \sqrt{E - \frac{1}{a(r)^2}} dr$$
is also known. Here, $r_{\pm}(E)$ are the upper and lower values
of the radius of the projection of the torus of energy $E$ and
angular momentum $n_0$ to $S^2$.  We may write the integral in the
form
$$ \int_{\R} ( E - x)_+^{\half} d\mu(x)$$
where $\mu$ is the distribution function $\mu (x):= |\{r:
\frac{1}{a(r)^2} \leq x\}|$ of   $\frac{1}{a^2}$, with $|\cdot|$
the Lebesgue measure. This  Abel transform is invertible and hence
$$d \mu(x) = \sum_{r: \frac{1}{a(r)^2} = x} |\frac{d}{dr} \frac{1}{a(r)^2} |^{-1}
dx$$ and therefore
$$ J(x) := \sum_{r: a(r) = x} \frac{1}{|a'(r)|}$$
are spectral invariants. By the simplicity assumption on $a$,
there are just two solutions of $a(r) = x$; the smaller will be
written $r_{-}(x)$ and the larger, $r_{+}(x).$ Thus,  the function
$$J(x) = \frac{1}{|a'(r_{-}(x))|} + \frac{1}{|a'(r_{+}(x))|} \leqno(6.9)$$
is a spectral invariant.

By studying  $H_{-1}$, we find in a somewhat similar way that
$$K(x) = |a'(r_{-}(x))| + |a'(r_{+}(x))|$$
is a spectral invariant. It follows that we can determine
$a'(r_+(x))$ and $a'(r_-(x))$.

Since both metrics $g_1$ and $g_2$ are assumed to belong to ${\cal
R}^*$, they are determined by their respective functions $a_j(r)$.
We conclude that $a_1 = a_2$ and hence that $g_1 = g_2$.

\subsection{Unconditional spectral determination}

We now ask  whether we can  remove the assumption that $g_2 \in
{\cal R}^*$ in Theorem \ref{ISPSR}?  The question is, if
Spec$\Delta_g =$ Spec$\Delta_h$ and $g \in {\cal R}^*$, then is $h
\in {\cal R}^*$? An affirmative answer would give a large
 class of metrics which are spectrally determined.  To the
 author's knowledge, the
only  metric on $S^2$  known to be  spectrally determined is the
canonical round one.

We  cannot answer this question, but we will give a partial result
suggesting that it is true.  The following theorem was stated in
\cite{Z7} and was worked out in a conversation with G. Forni in
1997. A similar but somewhat stronger conclusion  was drawn by K.
F. Siburg (Theorem 5.2 of \cite{S2}), under the stronger
hypothesis of an isospectral (or length-spectral) deformation.
After sketching our proof, we will also sketch his, which has
other applications.

\begin{theo} \label{FZ} Let $g \in {\cal R}^*$ and suppose that $h$ is any
metric on $S^2$ with simple clean length spectrum for which
Spec$\Delta_g =$ Spec$\Delta_h$. Then $h$ has the following
properties:
\begin{itemize}

\item (i) It has just one  isolated
non-degenerate closed geodesic $\gamma_h$ (up to orientation); all
other closed geodesics come in one-parameter families lying on
invariant tori in $S^*_h S^2$;

\item (ii) The Birkhoff normal form of $G_h^t$ at $\gamma_h$ is
identical to that of $G^t_g$ at its unique non-degenerate closed
orbit. Hence it is convergent.

 \item (iii) the geodesic flow $G^t_h$  of $h$ is $C^0$- integrable. That is, $S^*_h S^2$
 has
a $C^0$-foliation by $2$-tori invariant under $G^t_h$.
\end{itemize}
\end{theo}

 If we knew in (ii) that the Birkhoff transformation conjugating $G_h^t$
 to its Birkhoff normal form was convergent, then it would follow
 that $G_h^t$ is  completely integrable with global action-angle
 variables, and that it  would commute with a Hamiltonian torus
 action. We conjecture that this is the case. Statement (iii) shows that it is at least integrable in the $C^0$
 sense. At the present time, metrics on $S^2$ whose geodesic flows
 commute with Hamiltonian torus actions have not been classified.
 For a result which classifies such metrics on the torus (they must be flat), see
 \cite{LS}.

\subsubsection{Wave invariants for $g \in {\cal R}^*$}

The proof is based on a study of the wave trace formula in this
setting.

\begin{prop}  Suppose that $g \in {\cal R}^*$.  Then the trace of its wave group
has the form:
$$Tr e^{i t \sqrt{\Delta_g}} = e_0(t) + e_{\gamma_g}(t) + \sum_{{\cal T}} e_{\cal
T}(t)$$ where $$e_0(t)= C_n area(M,g) (t + i0)^{-n} + \dots$$ is
singular only at $t = 0$, where
$$e_{\gamma_g}(t) = c_{\gamma} (t - L_{\gamma} + i0)^{-1} + a_{\gamma 0} \log
(t - L_{\gamma} + i0) + \dots$$ and where
$$e_{{\cal T}} = c_{{\cal T}} (t - L_{{\cal T}} + i0)^{-3/2} + \dots.$$\end{prop}

Suppose now that $h$ is any other metric with Spec $\Delta_h$ =
Spec$\Delta_g$ and with simple length spectrum.  Then the wave
trace of $h$ has precisely the same singularities as the wave
trace of $g$. Since there is only one singularity of the order $(t
+ i 0)^{-1}$, there can exist only one non-degenerate closed
geodesic, proving (i). By Guillemin's inverse result, the Birkhoff
normal form of the metric and Poincar\'e map for $\gamma_h$ is the
same as for $\gamma_g$. In particular, the Poincare map ${\cal
P}_h$ of $\gamma_h$ is elliptic of twist type.
 From the fact that all other critical components have the singularity
 of a three-dimensional cone, it follows that in $S^*_hS^2$ the
 other closed geodesics come in one-parameter families. They weep
 out a surface foliated by circles, which can only be a two
 dimensional torus, proving (ii)

It is the third statement (iii) which requires a new idea. So far
we only know that periodic orbits lie on invariant tori, but we do
not know what lies between these tori. Aubry-Mather theory will
now close the gaps.

We recall that  Aubry-Mather theory is concerned  with an
area-preserving diffeomorphism  $\phi$ of an annulus $A= S^1
\times (a, b)$ \cite{KH, MF, Ka,  M, S1, S2}. Let $\tilde{\phi}$
denote a lift to $\R \times (a, b)$ with $\tilde{\phi}(x + 1, y) =
\tilde{\phi}(x, y) + (1, 0)$. The map $\phi$ is called a  {\it
monotone twist mapping} if it preserves the orientation of $A$, if
it preserves the boundary components and if the lift
$\tilde{\phi}(x_0, y_0) = (x_1, y_1)$ satisfies:
\begin{itemize}

\item The twist condition: $\frac{\partial x_1}{\partial y_0} >0;$

\item The exactness condition: $y_1 dx_1 - y_0 dx_0 = dh(x_0, x_1).$

\end{itemize}

If $a, b$ are finite, $\tilde{\phi}$ extends continuously to the
boundary as a pair of `rotations':
$$\tilde{\phi}(x, a) = (x + \omega_-, a),;\;\;\;\tilde{\phi}(x, b) = (x + \omega_+,
b).$$ The interval $(\omega_-, \omega_+)$ is called the {\it twist
interval} of $\phi$.

Let $\{(x_i, y_i)\}$ be an orbit of $\tilde{\phi}$. Its {\it
rotation number} is defined to be
$$\lim_{|i| \to \infty} \frac{x_i - x_0}{i}. $$
A curve $C \subset A$ is called an {\it invariant circle} if it is
an invariant set which is homeomorphic to the circle and which
separates boundary components. According to Birkhoff's invariant
circle theorem, an invariant circle is a Lipschitz graph over the
factor $S^1$ of $A$. Any invariant circle has a well-defined
rotation number (the common rotation number of orbits in the
circle) and the rotation number belongs to the twist interval.

 An orbit $\{(x_i, y_i)\}$ is determined by the sequence
$\{x_i\}$ of its $x$-coordinates. It is called {\it minimal} if
every finite segment is action-minimizing with fixed endpoints:
$$\sum_{i = k}^{n - 1} h(x_i, x_{i + 1}) \leq \sum_{i = k}^{n-1}
h(\xi_i, \xi_{i + 1}),\;\;\; \forall (\xi_k, \dots, \xi_n) \;
\mbox{with}\;\; \xi_k = x_k, \xi_n = x_n.
$$
The corresponding orbit orbit $(x_i, y_i)$ is called a {\it
minimal orbit}. The Aubrey-Mather theorem states (cf. \cite{Ka},
Theorem 1;) :
\medskip

{\it A monotone twist map possesses minimal orbits for each
rotation number $\omega \in (\omega_1, \omega_+)$ in its twist
interval. Every minimal orbit lies on a Lipschitz graph over the
$x$-axis.  For each rational rotation number $\omega =
\frac{p}{q}$ , there exists a periodic minimal orbit of rotation
number $\frac{p}{q}$. When $\omega$ is irrational, there exists
either an invariant circle with rotation number $\omega,$ or an
invariant Cantor set $E$.}
\medskip

The theorem also describes three possible orbit types in both the
rational or irrational case.

We now return to our problem on simple surfaces of revolution. We
fix  local transverse discs (Poincar\'e sections) $S_g,$ resp.
$S_h$ to the geodesic flows $G^t_g,$ resp. $G^t_h$ at the orbits
$\gamma_g,$ resp. $\gamma_h$. Concretely, the transversals can be
taken to be small variations of $\gamma'_g$, resp. $\gamma'_h$
moved up and down a small orthogonal geodesic arc to $\gamma'_g$,
resp. $\gamma_h'.$ Since $\gamma_g$ is non-degenerate elliptic,
its  Poincar\'e map ${\cal P}_g$ defines  area-preserving map of
the symplectic disc $S_g$ with a non-degenerate  elliptic fixed
point corresponding to $\gamma_g$. Since the Birkhoff normal forms
of ${\cal P}_g$ and ${\cal P}_h$ are the same, ${\cal P}_h$
defines the same kind of map of $S_h$. To obtain a twist map of an
annulus, we puncture out the fixed point of ${\cal P}_g, $ resp.
${\cal P}_h.$ We define the rotation angle $\omega_0$
corresponding to this orbit by continuity from nearby orbits,
which can be read off from the Birkhoff normal form. Indeed, the
Poincar\'e maps have the form
$$\left( \begin{array}{l} x \\ \\ y \end{array} \right) \to
\left( \begin{array}{ll} \cos 2 \pi \theta  & - \sin 2 \pi \theta
\\ &  \\ \sin 2\pi \theta & \cos 2 \pi \theta \end{array} \right)
\cdot
 \left(
\begin{array}{l} x \\ \\ y \end{array} \right) \; + \; O(|x|^2 + |y|^2),$$
with $\theta \sim  \omega_0 + \beta  (x^2 + y^2) + \cdots, $ as
$x^2 + y^2 \to 0.$ For further discussion of area-preserving maps
around elliptic fixed points and twist maps, see \cite{S2}.

The  foliation of $S^*S_g \backslash \gamma_h$ by $2$-tori
intersects $S_g$ in a foliation by invariant circles converging to
the fixed point. Circles with rational rotation numbers contain
only periodic orbits, and conversely all periodic orbits belong to
invariant circles with rational rotation numbers.

In the case of $h$, we know that once $\gamma_g$ is punctured out,
all  periodic orbits come in two-tori which project to invariant
curves in $S_h$ which are diffeomorphic to circles. The Birkhoff
normal form shows that they are invariant circles in the sense
that they are also homotopically non-trivial in the punctured
$S_h$ (indeed, the action-angle variables are essentially polar
coordinates on this disc).

  Let $C \subset S_h$ be an
invariant circle for ${\cal P}_h.$ Set:

\begin{itemize}
\item  $I(C)$ = area enclosed by $C$.

\item  $C_I$ = invariant circle enclosing an area $I$. (It is
clearly unique).

\item $\omega_I$ = rotation number of ${\cal P}_H |_{C_I}.$

\end{itemize}

Since ${\cal P}_h$ is a twist map,  the rotation number $\omega_I$
is a monotone increasing function of $I$.

Let $I_{+}$ denote the area of $S_h$ with respect to the
symplectic form. By shrinking $S_h$  we may assume that the
boundary of $S_h$ is a periodic circle $C_{I_+}$.

Since the origin is a fixed point of ${\cal P}_h$ and since the
rotation number is increasing, it is clear that the set of
rotation numbers lie in the interval $[0, \omega_{I_+}].$  It
follows by the Aubry-Mather theorem that every rational number
$p/q \in [0, \omega_{I_+}]$ is the rotation number of a periodic
circle $C_I \subset S_h.$

\begin{lem} Let $\alpha \in [0, \omega_{I_+}]$.  Then there exists an invariant
circle for ${\cal P}_h$ of rotation number $\alpha.$ \end{lem}

\noindent{\bf Proof} This follows from Corollary 6.1 of \cite{Ka},
which shows that as soon as Birkhoff periodic points of all
rational rotation numbers are constructed, then there exist orbits
of each irrational rotation number which are dense in an invariant
circle. Alternatively one could let $p_n/q_n \to \alpha$ and let
$C_{p_n/q_n}$ be the corresponding periodic circles.  Each
$C_{p_n/q_n}$ is a Lipschitz circle and the sequence of these
circles tends monotonically to a limit circle. It is Lipschitz and
its rotation number is $\alpha.$

To complete the proof, we make the observation:

\begin{lem} $S_h$ is foliated in the $C^0$ sense by invariant circles for
${\cal P}_h.$ \end{lem}

\noindent{\bf Proof}:  If not there exists an annulus $A \subset
S_h$ with boundary consisting of two invariant circles and
containing no invariant circles in its interior.  But by the
Aubrey-Mather theorem, there must exist a periodic point in $A$.
Since the periodic points come in circles, there must exist a
periodic circle, contradicting the non-existence of invariant
circles in $A$. QED

This completes the proof of Theorem \ref{FZ}.

\subsubsection{Isospectral class of an ellipse}

Instead of surfaces in ${\cal R}^*$, one can apply this reasoning
to the Dirichlet (or Neumann) problem for an ellipse $E_{a,b} =
\{(x, y): \frac{x^2}{a^2} + \frac{y^2}{b^2} = 1 \}.$ It is
well-known that ellipses have integrable billiards.

The ellipse has three distinguished periodic billiard orbits:
\begin{itemize}

\item The bouncing ball orbit along the  minor axis, which is a non-degenerate elliptic orbit;

\item The bouncing ball orbit along the  major axis, which is a non-degenerate hyperbolic orbit;

\item Its boundary.

\end{itemize}

All other periodic orbits come in one-parameter families. The
existence of a hyperbolic orbit means that the billiard flow is
very different from geodesic flows of metrics in ${\cal R}^*$.

The wave trace formula shows that any domain $\Omega$ with
$Spec(\Omega) = Spec(E_{a, b})$ has precisely one isolated
elliptic orbit, one isolated hyperbolic orbit. The accumulation
points in the length spectrum must be multiples of the perimeter
of the domain (a spectral invariant), so the boundary must be a
closed geodesic as well.

One can apply the  twist map theory either to the boundary orbit
or to the unique non-degenerate  elliptic orbit with isolated
length in the length spectrum. The argument above shows that there
exists a $C^0$ foliation by invariant circles at least near these
two orbits.

\subsection{\label{MLSR} Marked length spectral rigidity of
domains}

We now review a result due to  K.F. Siburg \cite{S1, S2} on
isospectral deformations of integrable systems  which has
interesting applications to metrics in ${\cal R}^*$ and also to
bounded plane domains.

The key invariant is the
 {\it mean minimal action}
\begin{equation} \alpha: [ \omega_-, \omega_+] \to \R,
\end{equation}  of a twist map $\phi$, which associates to a rotation number $\omega$ in
the `twist interval' the mean action
\begin{equation} \alpha(\omega) = - \lim_{N \to \infty} \frac{1}{2N} \sum_{i = -N}^N h(x_i, x_{i +
1})  \end{equation} of a minimal orbit $(q_i, \eta_i)$ of $\phi$
of rotation number $\omega$. It is a strictly convex function
which is differentiable at all irrational numbers. If $\omega =
p/q$, then $\alpha$ is differentiable at $\omega$ if and only if
there exists an invariant circle of rotation number $p/q$
consisting entirely of periodic minimal orbits. If a monotone
twist map possesses an invariant circle of rotation number
$\omega$, then every orbit on the circle is minimal (\cite{MF},
Theorem 17.4). In the case of a bounded plane domain, $h(q, q') =
- |q - q'|.$

It is observed by  K. F. Siburg \cite{S1} (Theorem 4.1)   that the
marked length spectrum  is essentially the same invariant as the
mean minimial action. The mean minimal action is therefore an
isospectral deformation invariant. He used this to give a proof
that there cannot exist isospectral deformations within ${\cal
R}^*$ (or within more general classes of deformations, see
\cite{S2}, Theorem 4.5). Indeed, the mean minimal actions $\alpha$
would all be the same. Hence if that  $\alpha_0$ of the original
surface is differentiable at all rationals, so are they all. But
this implies that they all have invariant circles of rational
rotation numbers. By taking limits, one obtains invariant circles
of all rotation numbers. (We remark that the trace formula already
shows that there existed invariant circles of all rational
rotation numbers).

Siburg  further connects the mean minimal  to the Melrose-Marvizi
invariants: Let  $\alpha^* : [-1, 1] \to \R$ denote its convex
conjugate of $\alpha$.  Then Siburg shows (loc. cit. p. 300) that
the Melrose-Marvizi invariants are algebraically equivalent to the
Taylor coefficients of $(\alpha^*)^{2/3}$ at $-1$.

It appears that the  only explicitly known mean minimal action is
that   of the disc $D$, where it is given by $\alpha (\omega) =
\frac{-1}{\pi} \sin \pi \omega$ (\cite{S1}); $\alpha$ is only
smooth when $\Omega = D$ (\cite{S1}, Theorem 4.6) It is likely
that it is computable in the case of an ellipse (perhaps in terms
of elliptic functions). Perhaps it can be proved that ellipses are
the unique domains with these particular minimal action functions
(which could be simpler than the long outstanding problem of
proving that they are the unique integrable billiard systems). If
so, this would prove that ellipses are spectrally rigid.

\begin{prob} \label{MEANMIN} Is the map from curvature functions $\kappa$ of convex plane
domains to the   mean minimal action $\alpha$ of the associated
convex domain  injective or finitely many to one? at least near
ellipses or under some additional analyticity or  discrete
symmetry condition?  \end{prob}


\begin{thebibliography}{HHHH}

\bibitem[AM]{AM} K. G. Andersson and R. B.  Melrose,
The propagation of singularities along gliding rays. Invent. Math.
41 (1977), no. 3, 197--232.

\bibitem[AKN]{AKN} V. I. Arnol'd, V.V.  Kozlov, and A. I.  Neushtadt,
{\it  Mathematical aspects of classical and celestial mechanics.
Dynamical systems, III}, pp. vii--xiv and 1--291, Encyclopaedia
Math. Sci., 3, Springer, Berlin, 1993.

\bibitem[BB]{BB} V.M.Babic, V.S. Buldyrev: {\it Short-Wavelength Diffraction Theory},
Springer Series on Wave Phenomena 4, Springer-Verlag, New York
(1991).

\bibitem[Ba]{Ba} V. Baladi, Periodic orbits and dynamical spectra. Ergodic Theory Dynam. Systems 18 (1998), no. 2, 255--292.

\bibitem[BB1]{BB1} R. Balian and C. Bloch, Distribution of eigenfrequencies for
the wave equation in a finite domain I: three-dimensional problem
with smooth boundary surface, Ann. Phys. 60 (1970), 401-447.

\bibitem[BB2]{BB2} R. Balian and C. Bloch,  Distribution of eigenfrequencies for the wave equation in a finite domain. III. Eigenfrequency density oscillations. Ann. Physics 69
(1972), 76--160.

\bibitem[B]{B} V. Bangert,
Geodesic rays, Busemann functions and monotone twist maps.
 Calc. Var. Partial Differential
Equations 2 (1994), no. 1, 49--63.

\bibitem[BLW]{BLW} A. Banyaga, R.  de la Llave, and C. E. Wayne,
Cohomology equations near hyperbolic points and geometric versions
of Sternberg linearization theorem.  J. Geom. Anal. 6 (1996), no.
4, 613--649 (1997).

\bibitem[BGR]{BGR} C. Bardos, J. C. Guillot, and J. Ralston,
La relation de Poisson pour l'\'equation des ondes dans un ouvert
non-borné, Comm. Partial Differential Equations 7 (1982), 905--958

\bibitem[BK]{BK} M. I. Belishev and Y. V.  Kurylev,  To the reconstruction of a Riemannian manifold via its spectral data (BC-method). Comm. Partial Differential Equations 17 (1992), no. 5-6, 767--804


\bibitem[Be]{Be} P. B\'erard, Transplantation et isospectralit\'e. Math. Ann. 292, 547--560 (1992)

\bibitem[Be2]{Be2} P. B\`erard,
Transplantation et isospectralité. II.  J. London Math. Soc. (2)
48 (1993), no. 3, 565--576.

\bibitem[Be3]{Be3} P. B\`erard, {\it Spectral geometry: direct and inverse problems}.
 With appendixes by Gérard Besson, and by B\`erard and Marcel Berger.
 Lecture Notes in Mathematics, 1207. Springer-Verlag, Berlin, 1986.


 \bibitem[Ber]{Ber} M. Berger, P.  Gauduchon, and E.  Mazet, {\it Le spectre d'une vari\`et\`e riemannienne.} Lecture Notes in Mathematics, Vol. 194 Springer-Verlag, Berlin-New York 1971


 \bibitem[BCG]{BCG} G. Besson, G. Courtois, and S.  Gallot,
Entropies et rigidit\'es des espaces localement symétriques de
courbure strictement négative.  Geom. Funct. Anal. 5 (1995), no.
5, 731--799.

 \bibitem[Bia]{Bia} M. Bialy,Convex billiards and a theorem by E. Hopf. Math. Z. 214 (1993), no. 1, 147--154


\bibitem[BKL]{BKL}  M.  Blank, G.Keller, and C.  Liverani, Ruelle-Perron-Frobenius spectrum for Anosov maps Authors:
(arXiv preprint nlin.CD/0104031).


\bibitem[BJP]{BJP} D. Borthwick, C. Judge, and P. A. Perry,
 Determinants of Laplacians and isopolar metrics on surfaces of infinite area. Duke Math. J. 118 (2003), no. 1, 61--102.


\bibitem[BG]{BG} T.P. Branson and P.  Gilkey,
The asymptotics of the Laplacian on a manifold with boundary.
Comm. Partial Differential Equations 15 (1990), no. 2, 245--272.

\bibitem[BGKV]{BGKV} T. Branson, P. B. Gilkey, K. Kirsten and D. Vassilevich,
 Heat kernel asymptotics with mixed boundary conditions, Nuclear Phys. B563 (1999), 603--626.

 \bibitem[BP]{BP} R. Brooks and P. A.  Perry,
Isophasal scattering manifolds in two dimensions.  Comm. Math.
Phys. 223 (2001), no. 3, 465--474.


\bibitem [BH]{BH}  J.Bruning and E.Heintze, Spektrale starrheit gewisser
Drehflachen, Math.Ann. 269 (1984), 95-101.

\bibitem[Ch]{Ch} J. Chazarain,
Construction de la paramétrix du problème mixte hyperbolique pour
l'équation des ondes.  C. R. Acad. Sci. Paris Sér. A-B 276 (1973),
A1213--A1215.

\bibitem[Ch2]{Ch2} J. Chazarain,
Formule de Poisson pour les variétés riemanniennes. Invent. Math.
24 (1974), 65--82


\bibitem[CdV]{CdV} Y. Colin de Verdi\`ere, Sur les longueurs des trajectoires périodiques d'un billard, in:
Dazord, Desolneux (eds.): {\it G\'eom\'etrie symplectique et de
contact}, Sem. Sud-Rhod. Géom. (1984), 122 139.

\bibitem[CdV2]{CdV2} Y. Colin de Verdi\`ere,
Spectre du laplacien et longueurs des géodésiques périodiques.
 C. R. Acad. Sci. Paris Sér. A-B 275 (1972), A805--A808.

 \bibitem[CdV3]{CdV3} Y.  Colin de Verdi\`ere,
Spectre du laplacien et longueurs des géodésiques périodiques. I,
II.  Compositio Math. 27 (1973), 83--106; ibid. 27 (1973),
159--184.

\bibitem[CdV4]{CdV4} Y. Colin de Verdi\`ere,  Spectre conjoint d'op\'erateurs pseudo-différentiels qui commutent.
II. Le cas int\'egrable.  Math. Z. 171 (1980), no. 1, 51--73.

\bibitem[C]{C} C. Croke,
Rigidity for surfaces of nonpositive curvature. Comment. Math.
Helv. 65 (1990), no. 1, 150--169.



\bibitem[C2]{C2} C. Croke,
Rigidity and the distance between boundary points. J. Differential
Geom. 33 (1991), no. 2, 445--464.

\bibitem[C3]{C3} C. Croke, Rigidity theorems in Riemannian geometry,
to appear in IMA Volume 137: {\it
Geometric Methods in Inverse Problems and PDE Control}. C.B.
Croke, I. Lasiecka, G. Uhlmann, and M. S.Vogelius eds.

\bibitem[CFF]{CFF} C.  Croke, A.  Fathi and J. Feldman,
The marked length-spectrum of a surface of nonpositive curvature.
Topology 31 (1992), no. 4, 847--855.

\bibitem[CS]{CS} C. B. Croke and V.  Sharafutdinov,
Spectral rigidity of a compact negatively curved manifold. Topology 37 (1998), no. 6, 1265--1273.

\bibitem[DS]{DS} C. Deninger and W.  Singhof,
A note on dynamical trace formulas. Dynamical, spectral, and
arithmetic zeta functions (San Antonio, TX, 1999), 41--55,
Contemp. Math., 290, Amer. Math. Soc., Providence, RI, 2001.

\bibitem[D]{D} H. Donnelly,
On the wave equation asymptotics of a compact negatively curved surface. Invent. Math. 45 (1978), no. 2, 115--137.

\bibitem[DG]{DG} J.J.Duistermaat and V.Guillemin, The spectrum of positive
elliptic operators and periodic bicharacteristics, Inv.Math. 24
(1975), 39-80.

\bibitem[Gh]{Gh} M. Ghomi, Shortest periodic billiard trajectories in convex bodies, to appear in Geom. Funct. Anal.

\bibitem[Gor]{Gor} C. Gordon,Survey of isospectral manifolds. {\it Handbook of differential geometry}, Vol. I,
 747--778, North-Holland, Amsterdam, 2000.

 \bibitem[Gor2]{Gor2} C. Gordon, CBMS Lecture Notes (to appear).

\bibitem[GorM]{GorM} C. Gordon and Y.  Mao,
Geodesic conjugacies of two-step nilmanifolds. Michigan Math. J.
45 (1998), no. 3, 451--481.

\bibitem[GMS]{GMS} C. Gordon, Y. Mao, and D. Schueth,
Symplectic rigidity of geodesic flows on two-step nilmanifolds.
 Ann. Sci. École Norm. Sup. (4)
30 (1997), no. 4, 417--427.




\bibitem[GWW]{GWW} C. Gordon, D.  Webb, and S. Wolpert,
 Isospectral plane domains and surfaces via Riemannian orbifolds. Invent. Math. 110 (1992), no. 1, 1--22.

\bibitem[G]{G} V. Guillemin,  Wave-trace invariants. Duke Math. J. 83 (1996), no. 2, 287--352.

\bibitem[G2]{G2} V. Guillemin,
Lectures on spectral theory of elliptic operators. Duke Math. J.
44 (1977), no. 3, 485--517.

\bibitem[G3]{G3} V. Guillemin,
Wave-trace invariants and a theorem of Zelditch. Internat. Math.
Res. Notices 1993, no. 12, 303--308

\bibitem[GK]{GK} V. Guillemin and D. Kazhdan,
 Some inverse spectral results for negatively curved $2$-manifolds. Topology 19 (1980), no. 3, 301--312.

 \bibitem[GM]{GM} V. Guillemin and  R. B. Melrose, The Poisson summation formula for manifolds with boundary,
  Adv. Math. 32 (1979), 204--232.

  \bibitem[GM2]{GM2} V. Guillemin and R. B. Melrose, The wave equation on manifolds with
  boundary, Colloque Boutet de Monvel (2003).

  \bibitem[GM3]{GM3} V. Guillemin and R. B. Melrose,
An inverse spectral result for elliptical regions in $ R\sp{2}$.
Adv. in Math. 32 (1979), no. 2, 128--148.

\bibitem[GM4]{GM4} V. Guillemin and R. B.  Melrose,
A cohomological invariant of discrete dynamical systems. E. B.
Christoffel (Aachen/Monschau, 1979), pp. 672--679, Birkhäuser,
Basel-Boston, Mass., 1981


  \bibitem[Gutz]{Gutz} M. C. Gutzwiller, J. Math. Phys. 12, 343 (1971).

  \bibitem[H]{H} U. Hamenst\"adt,
Cocycles, symplectic structures and intersection. Geom. Funct.
Anal. 9 (1999), no. 1, 90--140.

\bibitem[H2]{H2} U.  Hamenst\"adt,
Time-preserving conjugacies of geodesic flows. Ergodic Theory
Dynam. Systems 12 (1992), no. 1, 67--74.

\bibitem[HZel]{HZel} A. Hassell and S. Zelditch, Quantum ergodicity of boundary values of
eigenfunctions (2003), arxiv preprint math.SP/0211140.

\bibitem[HZel2]{HZel2} A. Hassell and S.  Zelditch,
Determinants of Laplacians in exterior domains. Internat. Math.
Res. Notices 1999, no. 18, 971--1004.



\bibitem[HZ]{HZ} A. Hassell and M.  Zworski,
Resonant rigidity of $S\sp 2$, J. Funct. Anal. 169 (1999), no. 2,
604--609.

\bibitem[Hil]{Hil} C. Hillar (private communication).

\bibitem[ISj]{ISj} A. Iantchenko and J. Sj\"ostrand,  Birkhoff normal forms for Fourier integral operators. II.
 Amer. J. Math. 124 (2002), no. 4, 817--850.

 \bibitem[ISZ]{ISZ} A. Iantchenko, J.  Sj\"ostrand and M. Zworski, Birkhoff normal forms in semi-classical
  inverse problems. Math. Res. Lett. 9 (2002), no. 2-3, 337--362.

\bibitem[I1]{I1} H.Ito, Convergence to Birkhoff normal forms for integrable
systems, Comment.Math.Helv. 64 (1989), 412-461.

\bibitem[I2]{I2} H.Ito, Integrability of Hamiltonian systems and Birkhoff normal
forms in the simple resonance case, Math.Ann.292 (1992), 411-444.

\bibitem[K]{K} W. Klingenberg, {\it Lectures on Closed
Geodesics}, Grundlehren der.\ math.\ W. {\bf 230}, Springer-Verlag
(1978).



  \bibitem[KKL]{KKL} A.  Katchalov, Y.  Kurylev, and Lassas, {\it Inverse boundary spectral problems}.
  Chapman and Hall/CRC Monographs and Surveys in Pure and Applied Mathematics, 123.  Boca Raton, FL, 2001
Ruishi On isospectral deformations of Riemannian metrics.
Compositio Math. 40 (1980), no. 3, 319--324.

\bibitem[Ka]{Ka} A. Katok,
Periodic and quasiperiodic orbits for twist maps. {\it Dynamical
systems and chaos} (Sitges/Barcelona, 1982), 47--65, Lecture Notes
in Phys., 179, Springer, Berlin, 1983.

\bibitem[KH]{KH} A. Katok and B. Hasselblatt, {\it Introduction to the modern theory of dynamical systems.}
 Encyclopedia of Mathematics and its Applications, 54. Cambridge University Press, Cambridge, 1995.

 \bibitem[KT]{KT} V.V. Kozlov and D. V. Treshch\"e v, {\it
 Billiards}, Trans. Math. Mono. 89, AMS, Providence R.I. (1991).

  \bibitem[Ku2]{Ku2} R. Kuwabara,
On isospectral deformations of Riemannian metrics. II. Compositio
Math. 47 (1982), no. 2, 195--205.

\bibitem[Ku3]{Ku3} R. Kuwabara,
On the characterization of flat metrics by the spectrum. Comment.
Math. Helv. 55 (1980), no. 3, 427--444.

\bibitem[LL]{LL} L. D. Landau and E.M. Lifshitz, {\it Quantum
Mechanics}, Course of Theoretical Physics Vol. 3, Third Edition,
Pergmanon Press (1977).

  \bibitem[LU]{LU} M. Lassas and G. Uhlmann,  On determining a Riemannian manifold from the
  Dirichlet-to-Neumann map. Ann. Sci. École Norm. Sup. (4) 34 (2001), no. 5, 771--787

  \bibitem[L]{L} V. F. Lazutkin,
Construction of an asymptotic series for eigenfunctions of the
"bouncing ball" type. In {\it  Asymptotic methods and stochastic
models in wave propagation problems. } (Russian) Trudy Mat. Inst.
Steklov. 95 (1968), 106--118, 215

\bibitem[La2]{La2} V. F. Lazutkin,
Existence of caustics for the billiard problem in a convex domain.
Math. USSR-Izv. 7 (1973), 185--214.




\bibitem[Ll]{L1} R. de la Llave, Analytic regularity of solutions of Livsic's cohomology
equation and some applications to analytic conjugacy of hyperbolic dynamical systems
, Ergodic Theory Dynam. Systems 17 (1997), 649--662.

\bibitem[L2]{L2} R. de la Llave,
On necessary and sufficient conditions for uniform integrability
of families of Hamiltonian systems. International Conference on
Dynamical Systems (Montevideo, 1995), 76--109, Pitman Res. Notes
Math. Ser., 362, Longman, Harlow, 1996.

\bibitem[LS]{LS} E. Lerman and N.  Shirokova,  Completely integrable
 torus actions on symplectic cones. Math. Res. Lett. 9 (2002), no. 1, 105--115.


\bibitem[Ma]{Ma} A. Majda,  High frequency asymptotics for the scattering matrix and
 the inverse problem of acoustical scattering. Comm. Pure Appl. Math. 29 (1976), no. 3, 261--291.

 \bibitem[Ma2]{Ma2} A. Majda,
A representation formula for the scattering operator and the
inverse problem for arbitrary bodies. Comm. Pure Appl. Math. 30
(1977), no. 2, 165--194.

\bibitem[M]{M} J. D. Meiss, Symplectic maps, variational principles, and transport.
Rev. Modern Phys. 64 (1992), no. 3, 795--848.

  \bibitem[MM]{MM} Sh. Marvizi and R. B. Melrose,
  Spectral invariants of convex planar regions. J. Differ. Geom. 17, 475--502
  (1982).

  \bibitem[MF]{MF} J.N. Mather and  G. Forni, Action minimizing orbits in Hamiltonian systems,
   in: Graffi (ed.): {\it Transition to Chaos in Classical and Quantum Mechanics}, Springer LNM 1589 (1992), 92--186.

  \bibitem[Me]{Me} R. B. Melrose,  The inverse spectral problem for planar domains.
 {\it  Instructional Workshop on Analysis and Geometry},
  Part I (Canberra, 1995), 137--160, Proc. Centre Math. Appl. Austral. Nat. Univ.,
  34, Austral. Nat. Univ., Canberra, 1996.

  \bibitem[Me2]{Me2} R. B. Melrose, "Geometric Scattering Theory," Cambridge Univ. Press, Cambridge/New York/Melbourne, 1995

  \bibitem[MS]{MS} R. B.  Melrose and J. Sj\"ostrand,
 Singularities of boundary value problems. I. Comm. Pure Appl. Math. 31 (1978), no. 5, 593--617.


  \bibitem[M1]{M1} R. Michel, Sur la rigidit\'e impos\'ee par la longueur des g\'eod\'esiques,
   Invent. Math. 65 (1981), 71--83.

   \bibitem[M2]{M2}
R. Michel, Restriction de la distance g\'eod\'esique \`a un arc et
rigidit\'e, Bull. Soc. Math. France 122 (1994), 435--442.

\bibitem[Mo]{Mo} J. Moser,
On invariant curves of area-preserving mappings of an annulus.
Nachr. Akad. Wiss. Göttingen Math.-Phys. Kl. II 1962 1962 1--20.

\bibitem[Mo1]{Mo1} J.Moser, Convergent series expansions for quasi-periodic motions,
Math.Annalen 169 (1967), 136-176.

\bibitem[Mo.2]{Mo.2} J.Moser, {\it Stable and Random Motions in Dynamical
Systems}, Ann.Math.Studies 77, Princeton U Press, Princeton
(1973).


\bibitem[Mu]{Mu}  R. G. Mukhometov, The reconstruction problem of a
two-dimensional Riemannian metric, and integral geometry,
(Russian), Dokl. Akad. Nauk SSSR 232 (1977), 32--35.

\bibitem[OPS]{OPS} B. Osgood, R. Phillips, and P. Sarnak,
 Moduli space, heights and
isospectral sets of plane domains. Ann. of Math. (2) 129 (1989),
no. 2, 293--362.

\bibitem[0]{O} J. P. Otal,
Le spectre marqué des longueurs des surfaces à courbure négative.
 Ann. of Math. (2) 131 (1990), no. 1, 151--162.

 \bibitem[O2]{O2} J. P. Otal,
Sur les longuer des g\'eod\'esiques d'une m\'etrique a courbure
n\'egative dans le disque, Comment. Math. Helv. 65 (1990),
334--347.


\bibitem[Pa]{Pa} V.K. Patodi,
Curvature and the fundamental solution of the heat operator. J.
Indian Math. Soc. 34 (1970), no. 3-4, 269--285 (1971).

\bibitem[PM]{PM} R. P\'erez-Marco,
Convergence or generic divergence of the Birkhoff normal form.
 Ann. of Math. (2) 157 (2003), no. 2,
557--574.

 \bibitem[PS]{PS}  V.M.Petkov and L.N.Stoyanov, {\it Geometry of Reflecting Rays
and Inverse Spectral Problems}, John Wiley and Sons, N.Y. (1992).

\bibitem[P]{P} G. Popov,  Invariants of the length spectrum and spectral invariants for convex planar domains. Commun. Math. Phys. 161, 335--364 (1994)

\bibitem[P2]{P2} G. Popov,  Length spectrum invariants of Riemannian manifolds. Math. Z. 213 (1993), no. 2, 311--351.


\bibitem[P3]{P3} G. Popov,
On the contribution of degenerate periodic trajectories to the
wave-trace.  Comm. Math. Phys. 196 (1998), no. 2, 363--383.


\bibitem[PT]{PT} G.  Popov and P. Topalov,  Liouville billiard tables and an inverse spectral
result. Ergodic Theory Dynam. Systems 23 (2003), no. 1, 225--248.

\bibitem[Rau]{Rau} J. Rauch,  Illumination of bounded domains. Amer. Math.
Monthly 85 (1978), no. 5, 359--361.

\bibitem[R]{R} M. Rouleux, Semi-classical integrability, hyperbolic flows and the Birkhoff normal form,
to appear in Canadian J. Math (see also  arXiv preprint
math.DS/0207026).







 \bibitem[SU]{SU} V. Sharafutdinov and G.  Uhlmann,
  On deformation boundary rigidity and spectral rigidity of Riemannian surfaces
   with no focal points. J. Differential Geom. 56 (2000), no. 1, 93--110.

   \bibitem[S1]{S1} K. F. Siburg, Aubry-Mather theory and the inverse spectral problem for
   planar convex domains. Israel J. Math. 113 (1999), 285--304.

   \bibitem[S2]{S2} K. F. Siburg, Symplectic invariants of elliptic fixed points.
   Comment. Math. Helv. 75 (2000), no. 4, 681--700.

   \bibitem[S3]{S3} K. F. Siburg, The minimal action in geometry and dynamics (to appear).

   \bibitem[S.M]{S.M} C.L.Siegel and J.Moser, {\it Lectures on Celestial Mechanics},
Grundlehren der math Wiss.Einz. 187, Springer-Verlag, Berlin
(1971).


   \bibitem[SjZ]{SjZ} J. Sj\"ostrand and  M. Zworski, Quantum monodromy
   and semi-classical trace formul\ae , J. Math. Pure Appl. 81 (2002), 1--33.

   \bibitem[SjZ2]{SjZ2} J. Sj\"ostrand and  M. Zworski, Quantum monodromy
   revisited, these proceedings.


   \bibitem[Sm]{Sm} L. Smith,
The asymptotics of the heat equation for a boundary value problem.
Invent. Math. 63 (1981), 467--493.

\bibitem[StU]{StU} P. Stefanov and G.  Uhlmann,  Rigidity for metrics with the same lengths of geodesics. Math. Res. Lett. 5 (1998), no. 1-2, 83--96.




\bibitem[St]{St} L. Stoyanov,  Rigidity of the scattering length spectrum. Math. Ann. 324 (2002), no. 4, 743--771.

\bibitem[St2]{St2} L. Stoyanov,  On the scattering length spectrum for real analytic obstacles. J. Funct. Anal. 177 (2000), no. 2, 459--488.

\bibitem[Su]{Su} T. Sunada,
Riemannian coverings and isospectral manifolds. Ann. of Math. (2)
121 (1985), no. 1, 169--186.


\bibitem[T]{T} S. Tanno,
A characterization of the canonical spheres by the spectrum. Math.
Z. 175 (1980), no. 3, 267--274.


\bibitem[T2]{T2} S. Tanno,
Eigenvalues of the Laplacian of Riemannian manifolds. T\v ohoku
Math. J. (2) 25 (1973), 391--403.

\bibitem[T3]{T3} S. Tanno,
A characterization of a complex projective space by the spectrum.
Kodai Math. J. 5 (1982), no. 2, 230--237

\bibitem[Wa]{Wa} C.P. Walkden, C. P. Livsic theorems for hyperbolic flows.
Trans. Amer. Math. Soc. 352 (2000), no. 3, 1299--1313

\bibitem[W]{W} K. Watanabe,
Plane domains which are spectrally determined. Ann. Global Anal.
Geom. 18 (2000), no. 5, 447--475.



 \bibitem[Z1]{Z1} S. Zelditch, Spectral determination of analytic bi-axisymmetric plane domains.
    Geom. Funct. Anal. 10 (2000), no. 3, 628--677.

\bibitem[Z2]{Z2} S.  Zelditch, The inverse spectral problem for surfaces of revolution.
   J. Differential Geom. 49 (1998), no. 2, 207--264.



\bibitem[Z3]{Z3} S.  Zelditch,  Wave invariants at elliptic closed geodesics.
Geom. Funct. Anal. 7 (1997), no. 1, 145--213.

\bibitem[Z4]{Z4} S. Zelditch,  Wave invariants for non-degenerate closed
geodesics. Geom. Funct. Anal. 8 (1998), no. 1, 179--217.


\bibitem[Z5]{Z5} S. Zelditch, The inverse spectral problem for analytic plane domains. I:
Balian-Bloch trace formula (http://arXiv.org/abs/math.SP/0111077,
to appear in Comm. Math. Phys.); and II: domains with symmetry,
(http://arXiv.org/abs/math.SP/0111078.)

\bibitem[Z6]{Z6} S. Zelditch,
Isospectrality in the FIO category. J. Differential Geom. 35
(1992), no. 3, 689--710.

\bibitem[Z7]{Z7} S. Zelditch, Normal forms and inverse spectral theory,
 Journées \'Equations aux d\'eriv\'ees partielles, 2--5 June 1998, GDR 1151 (CNRS).

 \bibitem[Z8]{Z8} S. Zelditch,
Maximally degenerate Laplacians. Ann. Inst. Fourier (Grenoble) 46
(1996), no. 2, 547--587.

\bibitem[Z9]{Z9} S. Zelditch,  Lectures on wave invariants. {\it Spectral theory and geometry}
 (Edinburgh, 1998), 284--328, London Math. Soc. Lecture Note Ser., 273, Cambridge Univ. Press, Cambridge, 1999.

 \bibitem[Z10]{Z10} S. Zelditch Inverse resonance problem for $Z_2$ symmetric
analytic obstacles in the plane,   IMA Volume 137: {\it Geometric
Methods in Inverse Problems and PDE Control}. C.B. Croke, I.
Lasiecka, G. Uhlmann, and M. S.Vogelius eds.


\bibitem[Z11]{Z11} S. Zelditch, Kuznecov sum formulae and Szeg\"o limit formulas on manifolds. Comm. P.D.E. 17, 221--260 (1992)


\bibitem[Zw2]{Zw2} M. Zworski,  Resonances in physics and geometry. Notices Amer. Math. Soc. 46 (1999), no. 3, 319--328.

\bibitem[Zw3]{Zw3} M. Zworski, Quantum resonances and partial differential equations. Proceedings of the International Congress of Mathematicians, Vol. III (Beijing, 2002), 243--252, Higher Ed. Press, Beijing, 2002.
\end{thebibliography}
\end{document}